\documentclass[12pt]{elsarticle}

\usepackage[utf8]{inputenc}

\usepackage{fullpage}

\usepackage{caption}
\usepackage{subcaption}

\usepackage{amsmath}
\usepackage{amssymb}
\usepackage{latexsym}

\usepackage{hyperref}

\usepackage{listings}

\usepackage{ifpdf}

\usepackage[mathlines]{lineno}

\usepackage{xspace}

\newcommand{\pd}[2]{\ensuremath{\frac{\partial #1}{\partial #2}}} % partial
\newcommand{\dx}{\ensuremath{\Delta x}}                     % Delta x
\newcommand{\dy}{\ensuremath{\Delta y}}                     % Delta y
\newcommand{\dt}{\ensuremath{\Delta t}}                     % Delta t
\newcommand{\R}{\ensuremath{\mathbb{R}}}                    % Real field
\newcommand{\wave}{\ensuremath{\mathcal{W}}\xspace}             % Wave
\newcommand{\fwave}{\ensuremath{\mathcal{Z}}\xspace}            % F-Waves
\newcommand{\cell}{\ensuremath{\mathcal{C}}\xspace}             % FV grid cell
\newcommand{\apdq}{\ensuremath{\mathcal{A}^+ \Delta Q}\xspace}      % A+dq
\newcommand{\amdq}{\ensuremath{\mathcal{A}^- \Delta Q}\xspace}      % A-dq
\newcommand{\geoclaw}{{\sc GeoClaw}\xspace}

\newcommand{\adcirc}{ADCIRC\xspace}

\begin{document}

\ifpdf
\DeclareGraphicsExtensions{.pdf, .png, .jpg, .tif}
\else
\DeclareGraphicsExtensions{.png, .jpg, .tif, .eps}
\fi

\begin{frontmatter}

\title{Adaptive Mesh Refinement for Storm Surge}

\author[ut]{Kyle T. Mandli}

\ead{kyle@ices.utexas.edu}
\ead[http://users.ices.utexas.edu/~kyle]{http://users.ices.utexas.edu/~kyle}

\author[ut]{Clint N. Dawson}
\date{\today}
\address[ut]{Institute for Computational Engineering and Science, University of Texas at Austin, 201 E 24th ST. Stop C0200, Austin, TX 78712-1229, USA}

\date{\today}

\begin{abstract}
An approach to utilizing adaptive mesh refinement algorithms for storm surge modeling is proposed.  Currently numerical models exist that can resolve the details of coastal regions but are often too costly to be run in an ensemble forecasting framework without significant computing resources.  The application of adaptive mesh refinement algorithms substantially lowers the computational cost of a storm surge model run while retaining much of the desired coastal resolution.  The approach presented is implemented in the \geoclaw framework and compared to \adcirc for Hurricane Ike along with observed tide gauge data and the computational cost of each model run.
\end{abstract}


\end{frontmatter}
% \linenumbers

% ==============================================================================
\section{Introduction} \label{sec:intro}

As computer technology advances, scientists continually attempt to use numerical modeling to better predict a growing number of high-impact geophysical events.  In particular, coastal hazards have become an increasing concern as the world's population continues to grow and move towards the coastline, in fact 44\% of the world's population lives within 150 km of the coast and 8 of the 10 largest cities in the world lie in that range \cite{Anonymous:td}. As a consequence, loss of life and property is becoming a larger concern than ever before.  One of the most recurring and wide spread hazards to many coastal communities is the inundation of coastlines that is associated with strong storms, one part of which is known as storm surge.  A storm surge is a rise in the sea accompanying extratropical or tropical cyclones, the strongest examples of which are hurricanes and typhoons.  Storm surges can cause massive amounts of damage, as was demonstrated by Hurricane Katrina, which caused an estimated \$81 billion of damage \cite{Blake:7uq}.  Of the world's largest cities, 4 lie within threat zones from tropical cyclones.  With the mounting evidence that severe storms may be increasingly common \cite{Anonymous:ws}, the task of modeling these events becomes even more critical to communities along the coasts.

Modeling of storm surges was first carried out by local empirical observations.  Unfortunately, for more severe storms such as Katrina, these types of prediction can grossly under-predict storm surge size and effect.  By the 1960's, scientists started using computer simulations to predict storm surge but, because these simulations were limited in resolution and size, these models had the same short-comings as the empirically-based models.  It was not until recently with increased observational evidence, improved efforts in modeling underlying physical processes, and increases in available computational power that substantial progress has been made simulating large-scale storm surge for use in hazard planning.

The current state-of-the-art numerical models for storm surge simulations rely on single-layer depth-averaged equations for the ocean and make assumptions about the ocean's response to a storm.  The National Weather Service (NWS) utilizes a storm surge model called ``Sea, Lake and Overland Surges from Hurricanes'', or SLOSH, which uses local grids defined for many regions of the United States coastline, to make predictions \cite{Jelesnianski:1992aa}.  These simulations are efficient enough that ensembles of runs can be made quickly for multiple different hurricane paths and intensities.  This capability can be critical for effective forecasting due to the uncertainty in the storm forecast.  The primary drawback to using the SLOSH model is the limited domain size and extents allowed due to the grid mapping used and formulation of the equations.

Another model currently in use is the Advanced Circulation Model (\adcirc), a finite element model which has been applied to southern Louisiana \cite{westerink:2008} and recently to Hurricane Ike \cite{Hope:2013jt}.  One of the key advantages \adcirc has it its use of an unstructured grid.  Unstructured models allow easy application of variable resolution, especially at the coastline where fine scale features need to be resolved.  They can also map to coastlines in a way even a cleverly mapped structured grid cannot.  Another advantage of unstructured grids relates to the importance of including entire ocean basins for surge predictions \cite{Blain:1994aa,Li:2013bf}.  Unstructured grids can allow the domain of the numerical model to stretch well away from coastlines to include ocean basins while reducing the cost of the model by substantially decreasing resolution in the basin compared to the coastal regions.  Unfortunately these models, even with the above advantages, can still be computationally costly and require a large amount of computing resources in order to compute ensemble forecasts without the degradation of their resolution benefits.

In this paper we present an alternative computational framework and methodology to bridge the gap between the numerical cost of the unstructured grid storm surge models and the efficient but unresolved models currently in use at the NWS.  The approach leverages adaptive mesh refinement (AMR) algorithms to retain the resolution required to resolve coastal inundation but only when necessary so that ensemble calculations are still feasible.  This is accomplished by allowing nested structured grids of variable resolution to vary in time and space thereby capturing the spatial advantages of the unstructured grid approach but only when needed, and therefore decreasing the computational cost substantially.  The framework in question, \geoclaw, has successfully been used previously for tsunami modeling where similar computational requirements are present \cite{Berger:2011vi}.

% ==============================================================================
\section{Numerical Approach} \label{sec:numerical_approach}

The mathematical model for storm surge we will consider uses the classical shallow water equations with the addition of appropriate source terms for bathymetry, bottom friction, wind friction, non-constant surface pressure and Coriolis forcing which can be written as
\begin{equation} \label{eq:swe}
    \begin{aligned}
    &\pd{}{t} h + \pd{}{x} (hu) + \pd{}{y} (hv) = 0, \\
    &\pd{}{t}(hu) + \pd{}{x} \left(hu^2 + \frac{1}{2} g h^2 \right ) + \pd{}{y} (huv) = ~~ fhv - gh \pd{}{x} b + \frac{h}{\rho} \left(-\pd{}{x} P_A + \rho_{\text{air}} C_w |W| W_x - C_f |\vec{u}| u \right )\\
    &\pd{}{t} (hv) + \pd{}{x} (huv) + \pd{}{y} \left (hv^2 + \frac{1}{2} gh^2 \right) = -fhu - gh \pd{}{y} b + \frac{h}{\rho} \left(-\pd{}{y} P_A + \rho_{\text{air}} C_w |W| W_y - C_f |\vec{u}| v \right )
    \end{aligned}
\end{equation}
where $h$ is the fluid depth, $u$ and $v$ the depth-averaged horizontal velocity components, $g$ the acceleration due to gravity, $\rho$ the density of water, $\rho_\text{air}$ the density of air, $b$ the bathymetry, $f$ the Coriolis parameter, $W = [W_x, W_y]$ is the wind velocity at 10 meters above the sea surface, $C_w$ the wind friction coefficient, and $C_f$ the bottom friction coefficient.  The value of $C_w$ is defined by Garratt's drag formula \cite{Garratt:1977aa} as
\begin{equation} \label{eq:wind_drag_coefficient}
    C_w = \text{min}(2\times10^{-3}, (0.75 + 0.067 |W|) \times 10^{-3})
\end{equation}
and the value of the friction coefficient $C_f$ is determined using a hybrid Chezy-Manning's $n$ type friction law
\begin{equation} \label{eq:friction_term}
    C_f = \frac{g n^2}{h^{4/3}} \left[1-\left(\frac{h_{\text{break}}}{h}\right)^{\theta_f} \right]^{\gamma_f / \theta_f}
\end{equation}
where $n$ is the Manning's $n$ coefficient and $h_{\text{break}} = 2$, $\theta_f = 10$ and $\gamma_f = 4/3$ parameters control the form of the friction law.

The numerical approach proposed to solve \eqref{eq:swe} falls under a general class of high resolution finite volume methods known as wave-propagation methods, described in detail in \cite{LeVeque:2002aa}.  These methods are Godunov-type finite volume methods requiring the specification of a Riemann solver to update each grid cell in the domain.  On top of these methods adaptive mesh refinement is employed to allow for variable spatial and temporal resolution as the simulation progresses.  These methods have been implemented together in \geoclaw, a package that was originally designed to model tsunamis \cite{Berger:2011du} and other depth-averaged flows \cite{Berger:2011vi,George:2010p2275,Mandli:2011te}.  The rest of this section is dedicated to describing the salient points of the AMR approach and how storm surge physics are represented in \geoclaw.  A brief review of wave-propagation methods can be found in \ref{app:wave_propagation_methods} along with a basic outline of the Riemann solver employed in \ref{app:riemann_solver}.

\subsection{Adaptive Mesh Refinement} \label{sub:amr}

\begin{figure}[tb]
    \centering
    \begin{subfigure}[b]{0.40\textwidth}
        \includegraphics[width=\textwidth]{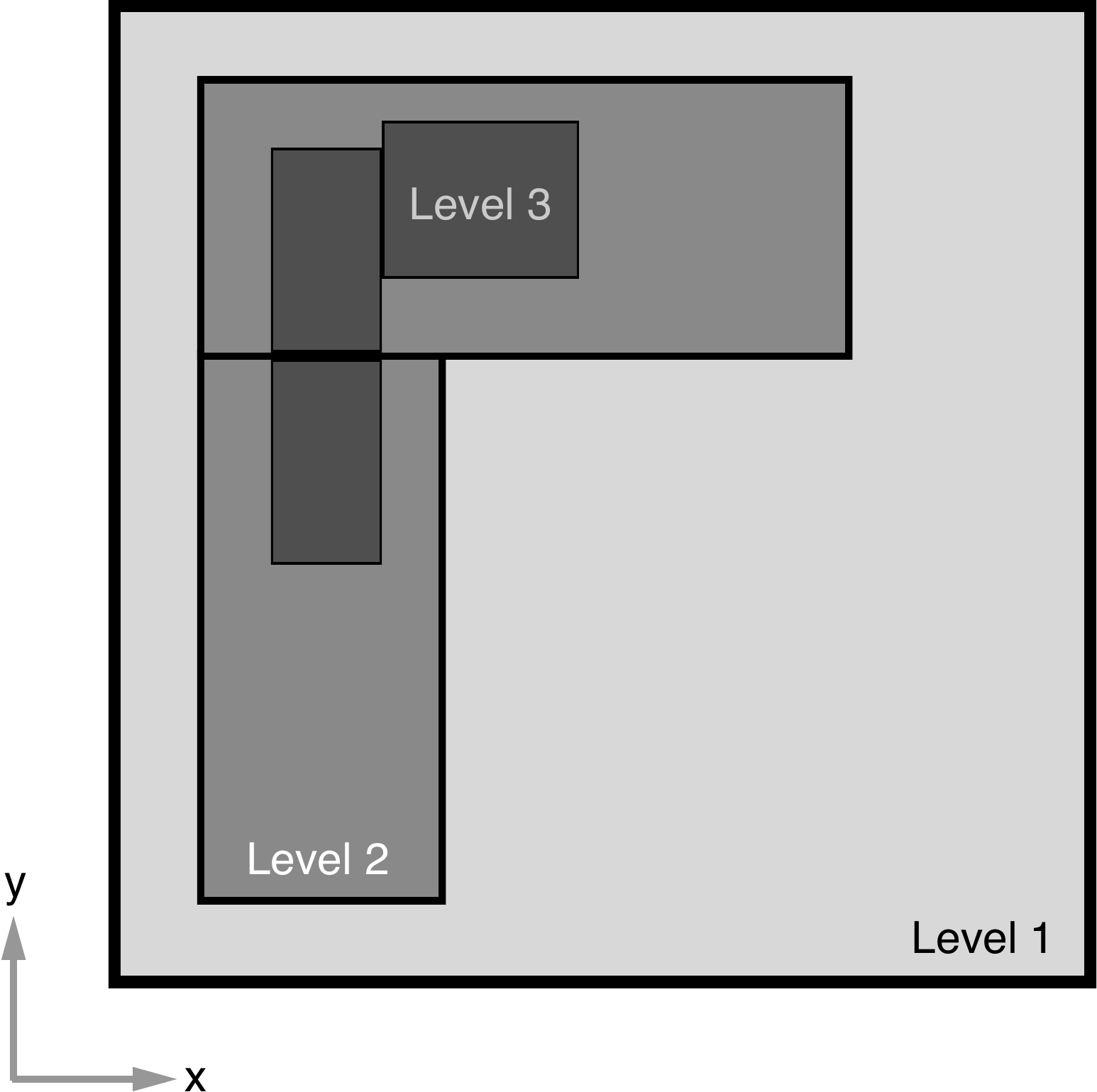}
        \caption{}
        \label{fig:planar_amr_view}
    \end{subfigure}
    \begin{subfigure}[b]{0.58\textwidth}
        \includegraphics[width=\textwidth]{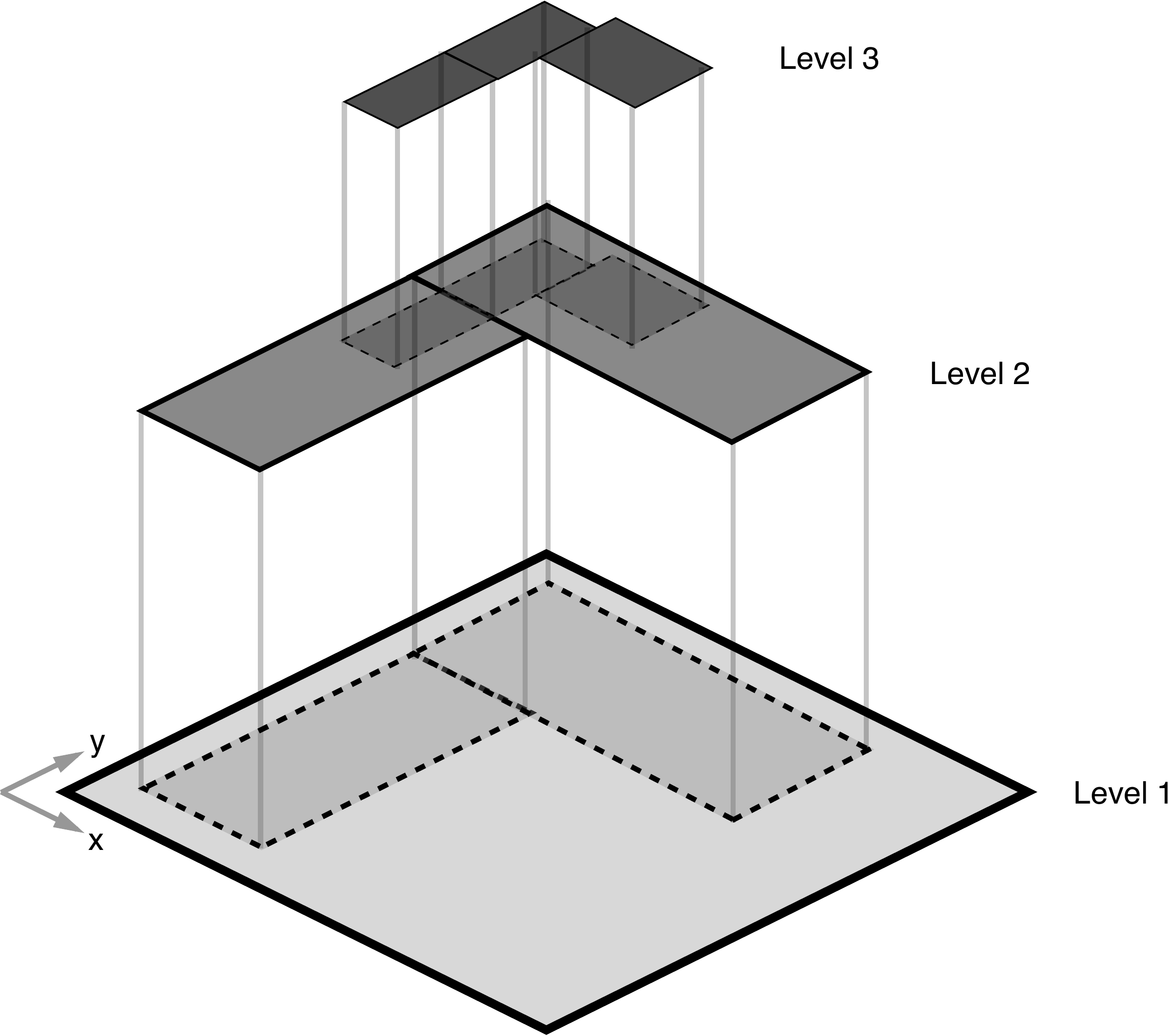}
        \caption{}
        \label{fig:isometric_amr_view}
    \end{subfigure}
    \caption{Views of a sample patch-based AMR grid.  Figure~\ref{fig:planar_amr_view} shows a planar view with three levels of refinement.  In this instance, there are 3 level 3 grids, 2 level 2 grids, and 1 level 1 grid.  Figure~\ref{fig:isometric_amr_view} shows an isometric view of the same grid layout with the z-direction representing different levels.}
    \label{fig:amr_diagrams}
\end{figure}

Adaptive mesh refinement is a core capability of \geoclaw as it allows the resolution of disparate spatial and temporal scales common to geophysical applications such as storm surge and tsunamis.  The patch-based AMR  approach used in \geoclaw employs a set of overlapping logically rectangular grids that correspond to one of many levels of refinement.  The first of these levels, enumerated starting at $\ell=1$, contains grids that cover the entire domain at the coarsest resolution.  The subsequent levels $\ell \ge 2$ represent progressively finer resolutions by a set of prescribed ratios $r^\ell$ in time and space such that
\[
    \dx^{(\ell+1)} = \dx^{(\ell)} / r^{(\ell)}_x, ~~~~ \dy^{(\ell+1)} = \dy^{(\ell)} / r^{(\ell)}_y, ~~~~~\text{and}~~~~~ \dt^{(\ell+1)} = \dt^{(\ell)} / r^{(\ell)}_t.
\]
Each subsequent level is properly nested within the union of grids in the next level coarser (see Figure~\ref{fig:amr_diagrams} for an example of the structure of these nested grids).  With this hierarchy of grids, the evolution of the solution follows as:
\begin{enumerate}
    \item Evolve the level 1 (coarsest) grids one time step to $t^{n+1}$.
    \item Fill the ghost cells of all level 2 grids by temporal and spatial interpolation.
    \item Evolve the level 2 grids the number of time steps determined by $\dt^{(2)}$ needed to reach $t^{n+1}$.
    \item Recursively continue to fill in ghost cells, evolving each level $\ell$ after the coarser level $\ell-1$ has been evolved until all levels are at time $t^{n+1}$.
    \item Fill in regions where the grids overlap with the best available data using interpolation.
    \item Adjust coarse cell values adjacent to finer cells to preserve conservation of mass (see \cite{Berger:1998aa} for a discussion on this topic).
\end{enumerate}

The key benefit of adaptive mesh refinement is the ability to change resolution as the simulation progresses.  This is done in a process that involves using a local criteria to flag each cell that requires refinement to the next level.  The algorithm then clusters the flagged cells into new rectangular patches which attempts to minimize the number of grids created and the number of grid cells unnecessarily refined \cite{Berger:1991ty}.  After the new grid structure is created, the previous solution values are copied into the new grid cells if no change in resolution was required or an interpolation and averaging is performed to either coarsen or refine the data available from a different level of refinement.

For the shallow water equations there are three primary areas where care must be taken when implementing adaptive mesh refinement.  The first involves the interpolation of the solution and bathymetry.  When interpolating where the water column is at rest but bathymetry is varying, using the depth $h$ as the interpolation field will lead to a sea surface $\eta$ that is not at rest and will result in the creation of spurious waves.  To avoid this, interpolation is done with the sea surface instead.  From this the depth is computed and the resulting momentum interpolated.  This process can be perfomed in wet cells while conserving mass and momentum.  This is not the case in the near shore where cells may change from their wet (or dry) state.  In order to avoid spurious waves being created mass cannot be conserved and is either lost or gained.  This is a result of different grid resolutions of bathymetry being used through the computation.  For these cases it is often desirable to refine coast lines before a wave arrives such that the loss or gain in mass does not effect the wave itself and overall leads to a negligible increase in mass overall in the simulation.  In the end, the following properties are true for the interpolation approach used.
\begin{itemize}
    \item Mass is conserved except possibly near coast lines.
    \item Momentum is conserved if mass is conserved or gained and lost if mass is lost.
    \item New extrema in the sea surface and velocity are not created.
\end{itemize}
More discussion and the derivation of the rules for coarsening and refinement can be found in \cite{Berger:2011du}.

The second area of concern is the form and types of refinement criteria that should be used.  In \geoclaw, a number of different solution based criteria are used to determine refinement.  The first criteria is triggered if the absolute difference between the initial sea-level $\eta_{\text{sea-level}}$ and the calculated sea-surface $\eta$ in cell $i$,$j$ is greater than some tolerance $T_{\text{wave}}$, \emph{i.e.}
\[
    T_{\text{wave}} < |\eta_{i,j} - \eta_{\text{sea-level}}|
\]
will mark the $i^{\text{th}},j^{\text{th}}$ cell as needing refinement.  Additionally water currents can be used as a refinement criteria.  This criteria uses a set of speed tolerances specified at each level, only marking a cell as needing refinement if the tolerance $T_\text{speed}$ for the corresponding level is less than the current speed in that cell.  For example, if tolerances for level 1, 2, and 3 were set to $2$, $2.5$, and $3$ m/s respectively, a cell with a speed of $2.6$ m/s currently at level 1 then would be marked as needing refinement.  Conversely, if the cell is already at level 2 then no refinement is necessary.  In addition to these physics based criteria, a user can specify regions with minimum or maximum refinement constraints.  There are also constraints based on initial depth of the water allowing restriction of refinement to only shallower depths.  In an instance where a conflict between the physics based criteria and the region constraints occurs the region constraints take precedence.

The last concern of note involves the relative spatial and temporal refinement ratios.  Isotropic refinement ratios, when the spatial and temporal refinement ratios are identical, are generally needed in AMR to satisfy the CFL condition.  In the case of the shallow water equations and inundation modeling the desired spatial refinement is usually highest at the coastline where the wave speeds are lowest.  This allows for the use of anisotropic refinement between the spatial and temporal directions to take advantage of the slower wave speeds located near shore.  Along with anisotropic refinement, \geoclaw allows for the automatic determination of the optimal temporal refinement ratios since the wave speed estimate for the shallow water equations is inexpensive to evaluate.

\subsubsection{Storm Based AMR Criteria} \label{ssub:ss_amr}

Additional criteria related to the location and strength of the storm can also be important to resolving storm surge.  The first available criteria is based on the distance of a cell from the eye of the storm.  As was the case with the speed criteria, these tolerances are specified by level such that a cell will be flagged as needing refinement if it is within the specified radius of the storm's eye, i.e.
\[
    T_{r}^\ell < r_\text{eye}
\]
for level $\ell$.  The second criteria flags based on storm intensity, wind speed tolerances are used in a similar way.  Often the wind speed tolerances and eye distance tolerances give similar but slightly different regions of refinement depending, for instance, on the strength of the storm.

\subsection{Source Term Evaluation} \label{ssub:source_term_numerics}

Since the storm's impact on the ocean depends only on the momentum source terms for the wind-stress and pressure gradients, it is important that all the momentum source terms are handled carefully.  Apart from the storm's impact, the rest of the momentum source terms are bathymetric, friction and Coriolis source terms.  The bathymetric source terms are handled in the Riemann solver so that well-balancedness can be maintained along with other useful properties (see~\ref{app:riemann_solver} and \cite{George:2008aa} for these details).  For the remaining source terms, \geoclaw uses an operator-splitting approach to evaluate the remaining source terms.  This involves solving two simpler problems, $q_t + A(q) q_z + B(q) q_y = S_1(q)$ and $q_t = S_2(q)$, and combining updates of the vector $q$ in a consistent manner.  The method employed uses Godunov splitting which involves solving each simpler system alternately for the full time step $\dt$.  Although this approach is formally only first-order accurate, observed errors due to the splitting do not dominate the overall error in practice (see \cite{LeVeque:2002aa} for a more thorough discussion of these issues).  

Starting with the friction and Coriolis, the bottom friction terms can be evaluated as
\[
    (hu)_t = C_f hu ~~~~\text{and}~~~~ (hv)_t = C_f hv
\]
and are solved using a backwards Euler method for computing the loss of momentum so that
\[
    (hu)^{n+1}_{ij} = \frac{(hu)^n_{ij}}{1 + (C_f)_{ij} \dt} ~~~~\text{and}~~~~
    (hv)^{n+1}_{ij} = \frac{(hv)^n_{ij}}{1 + (C_f)_{ij} \dt}
\]
where the drag coefficient is computed using the previous time step's state as
\[
    (C_f)_{ij} = \frac{gn_{ij}^2}{(h^n_{ij})^{-7/3}} \sqrt{[(hu)^n_{ij}]^2 + [(hv)^n_{ij}]^2} \left[1-\left(\frac{h_{\text{break}}}{h^n_{ij}}\right)^{\theta_f} \right]^{\gamma_f / \theta_f}.
\]
The values $h_{\text{break}}$, $\theta_f$, and $\gamma_f$ are constant parameters that define the hybridization between the Chezy and Manning's $n$ formulation from \eqref{eq:friction_term}.  The Coriolis terms
\[
    (hu)_t = -f hv ~~~~\text{and}~~~~ (hv)_t = f hu
\]
are evaluated using a matrix exponential up to the 4th term in the series such that the update becomes
\[
    \begin{bmatrix}
        hu \\ hv     
    \end{bmatrix}^{n+1}_{ij} = 
    \begin{bmatrix}
        hu \\ hv     
    \end{bmatrix}^{n}_{ij} \cdot 
    \begin{bmatrix}
        1 - \frac{1}{2} (f \dt)^2 + \frac{1}{24} (f \dt)^4 & f \dt - \frac{1}{6} (f \dt)^3 \\
        - f \dt + \frac{1}{6} (f \dt)^3 & 1 - \frac{1}{2} (f \dt)^2 + \frac{1}{24} (f \dt)^4
    \end{bmatrix}
\]
where $f = 2 \Omega \sin y$ with $\Omega = 2 \pi / 8.61642\times10^4$ and $y$ the longitudinal coordinate in radians. 

\subsubsection{Storm Representation and Source Term Evaluation} \label{ssub:storm_source_terms}

One of the significant additions to \geoclaw needed to model storm surge is the representation of the storm fields and their source terms.  In the forecasting setting that we are targeting, pressure and wind fields of a storm are often evaluated based on the empirically derived Holland model which provides profiles of the wind speed and pressure in a hurricane based on storm characteristics  \cite{Holland:1980aa}.  These profiles are then rotated to produce a two-dimensional, rotationally symmetric field.  The wind speed profile $W(r,t)$ takes the form
\begin{equation} \label{eq:holland_wind}
    W(r,t) = \sqrt{\left (\frac{r_{\text{w}}}{r}\right)^B W_{\text{max}}^2 e^{(1 - (r_{\text{w}} / r)^B)}  + \frac{(r f)^2}{4}} - \frac{r f}{2}
\end{equation}
where $r$ is the radial coordinate centered at the eye of the storm and the parameters $r_\text{w}$ and $W_\text{max}$ are given by the storm forecast and represent the radius of maximum winds, and the maximum wind speed, respectively.  The Holland $B$ parameter, an empirical fitting parameter from \cite{Holland:1980aa}, takes the form
\[
    B=\frac{\rho_\text{air} (W_\text{max}')^2 e}{P_n - P_c}
\] 
where $P_n$ and $P_c$ are the background and central storm atmospheric pressures.  The value $W_\text{max}'$ is a correction form of $W_\text{max}$ that accounts for storm translation and movement in spherical coordinates.  The pressure profile is similarly
\begin{equation} \label{eq:holland_pressure}
      P_A = P_c + (P_n - P_c) e^{-(r_\text{w} / r)^B}.
\end{equation}
The Holland wind field velocities from \eqref{eq:holland_wind} are converted to a cyclonic wind field via
\[
    \vec{W} = \begin{bmatrix}
        - W \sin \theta \\ W \cos \theta
    \end{bmatrix}
\]
where $\theta$ is the azimuthal angle with respect to the storm's eye.  Furthermore, the storm translation speed is added to take into account the relative motion of the storm and atmosphere.  The last modification to the fields is a ramping function, $\mathcal{R}(r)$, applied at a radial distance from the storm to smoothly decrease the wind velocity to zero and the pressure to the background pressure via
\[
    \mathcal{R}(r) = \frac{1}{2} \left(1 - \tanh\left( \frac{r - R_p}{R_w}\right)\right)
\]
where $r$ is the radial distance from the eye of the storm, $R_w$ the width of the ramping function, and $R_p$ the storm parameter recording the radius of the last closed iso-bar.  $R(r,t)$ is then applied to the wind and pressure fields as
\[
    W(r) = W(r) \cdot \mathcal{R}(r)~~~~\text{and}~~~~ P(r) = P_n + (P(r) - P_n) \cdot \mathcal{R}(r).
\]

A storm's intensity, speed, and location is allowed to vary in time via forecasted and best-track data at the times specified.  \geoclaw then reconstructs the wind and pressure fields using linear interpolation between time points so that the storm forecast evolves continuously rather than only at the time points available.  In addition to the intensity of the storm, the velocity and location are linearly interpolated taking into account that the storm is traveling on the surface of a sphere.  Once the simulation passes the end of the forecasted data, these parameters are extrapolated in time and kept constant.

Turning now to the evaluation of the source terms due to the wind and pressure fields, we utilize source term splitting as described earlier.  The wind source term $\rho_{\text{air}} C_w |W| W$ is evaluated using a single forward Euler time-step with
\[
    (hu)^{n+1}_{ij} = \frac{h^n_{ij}}{\rho} \rho_\text{air} (C_w)^n_{ij} |W^n_{ij}| (W_x)^n_{ij} ~~~~\text{and}~~~~ 
    (hv)^{n+1}_{ij} = \frac{h^n_{ij}}{\rho} \rho_\text{air} (C_w)^n_{ij} |W^n_{ij}| (W_y)^n_{ij}
\]
where the coefficient of drag is defined in \eqref{eq:wind_drag_coefficient}.  The pressure term is also handled via a forward Euler approach such that
\[
    (hu)^{n+1}_{ij} = \frac{h^n_{ij}}{\rho} \frac{P_{i-1j} - P_{i+1 j}}{2 \dx} ~~~~\text{and}~~~~ 
    (hv)^{n+1}_{ij} = \frac{h^n_{ij}}{\rho} \frac{P_{ij-1} - P_{i j+1}}{2 \dy}.
\]
The values of $\dx$ and $\dy$ for the second order-centered finite difference used for the pressure gradient were also converted from the latitude-longitude coordinate system to meters.

% ==============================================================================
\section{Comparisons} \label{sec:comparisons}

As a demonstration of the advantages of using adaptive mesh refinement for storm surge, \geoclaw was used to simulate Hurricane Ike.  \geoclaw's results were then compared to gauge data taken during the storm from \cite{Kennedy:2011kt}, and to \adcirc results, previously validated for Hurricane Ike in \cite{Hope:2013jt}.  The intention was to do the comparison in a forecasting type of scenario and as a consequence some forcing terms and resolution were sacrificed. Each simulation was computed 3 days ahead of the approximate landfall time and run to 18 hours after landfall.

\subsection{Hurricane Ike} \label{sub:hurricane_ike}

Hurricane Ike was a storm that caused significant devastation with maximum sustained winds of 54 m/s (10 minutes average) and a minimum pressure of 935 mbar, making landfall in the United States near Galveston, TX on September 13th, 2008 \cite{Hope:2013jt}.  One of the key features of the U.S. landfall was a prominent forerunner, a rise in water elevation that arrived prior to the arrival of the storm, leading to an increase of the overall surge throughout the storm.  The source of this prominent forerunner is thought to be Ekman setup on the Louisiana-Texas shelf \cite{Kennedy:2011kt} as the storm moved parallel to the shelf.  This also lead to significant setup in the relatively low lying coastal areas of Louisiana and Texas east of Galveston as the storm winds pushed the surge westward.  For a more complete description of the processes of interest during Hurricane Ike see \cite{Hope:2013jt}.

\subsection{\adcirc Run Description} \label{sub:adcirc_ike}

ADCIRC uses a continuous Galerkin finite element method that solves the modified shallow water equations \cite{westerink:2008, LuettichJr:2004ws, Dawson:2006cr}.  It has been validated using hindcasted storm data for a variety of hurricanes, in particular for Louisiana and Texas.  A key advantage of ADCIRC over many other available models is the unstructured grids it employs to model storm surge allowing for high levels of resolution near the coasts and less in the Gulf of Mexico and Atlantic Ocean.  Even with an unstructured grid, however, the number of nodes and elements required to run a storm surge simulation is massive due to the required domain size and resolution.  Because of this, ADCIRC employs multi-process parallelism via the message passing interface (MPI) standard to produce results in a reasonable amount of time.  Work has been done to make this parallelization as efficient and scalable (on the order of 10,000 processes) as possible \cite{Tanaka:2011cq}.

The \adcirc simulation compared here is based on a forecasting mode commonly used to model surge along the Texas coastline.  The grid used contains 3,331,560 nodes and 6,633,623 elements concentrated along the Texas coastline and identical in that region to the grid used to hindcast Hurricane Ike in \cite{Hope:2013jt}.  The grid cell sizes range from 30 kilometers to 37 meters with a mean size of 700 meters.  Along with the bathymetry, the bottom friction coefficient is allowed to vary at every node of the grid.  Additionally, because the simulation was run in forecast mode, a Holland \cite{Holland:1980aa} based storm model was used and tides, wave-stress, and riverine input were not included in the computation.

\begin{figure}[htb]
    \centering
    \begin{subfigure}[b]{0.48\textwidth}
        \includegraphics[width=\textwidth]{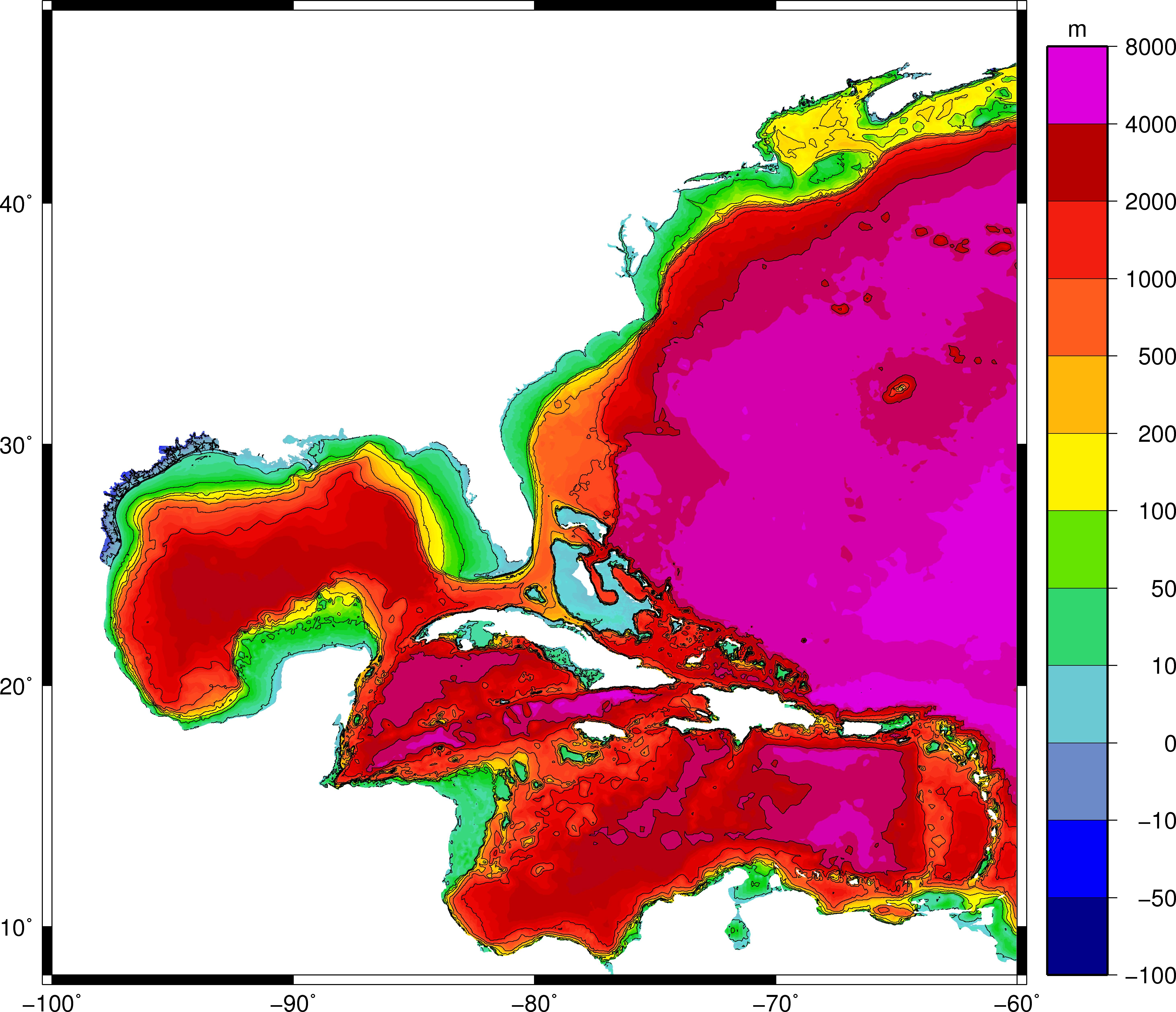}
        \caption{}
        \label{fig:adcirc_bathy}
    \end{subfigure}
    \begin{subfigure}[b]{0.48\textwidth}
        \includegraphics[width=\textwidth]{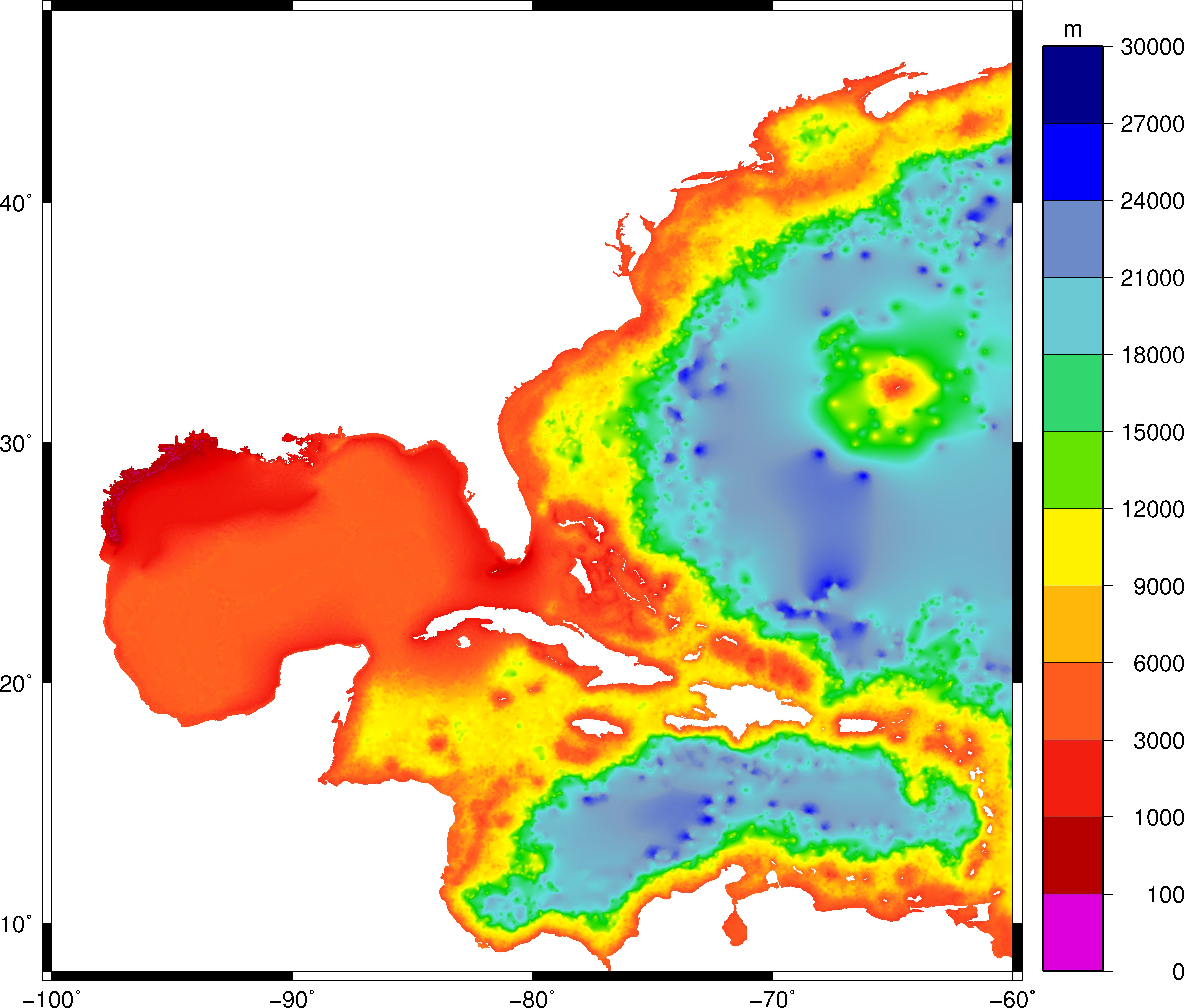}
        \caption{}
        \label{fig:adcirc_grid_resolution}
    \end{subfigure}
    \caption{Characteristics of the grid used in the \adcirc comparison run.  Figure~\ref{fig:adcirc_bathy} shows the bathymetry used in the \adcirc grid, Figure~\ref{fig:adcirc_grid_resolution} shows the resolution of the elements used.}
    \label{fig:adcirc_domain}
\end{figure}

\subsection{\geoclaw Run Description} \label{sub:geoclaw_ike}

In this section we describe the setup for the \geoclaw simulation.  The data and settings used for this simulation can be obtained at \url{http://github.com/clawpack/apps/tree/master/storm_surge/gulf/ike} along with the \geoclaw software itself at \url{http://www.clawpack.org}.

\subsubsection{Bathymetry Sources} \label{ssub:bathymetry}

The bathymetry data used is from the NOAA NGDC US coastal relief model grids with 3 seconds of resolution in the Louisiana-Texas shelf and coastal regions \cite{USCoastalRelief:qiDPbJY4}.  Outside of these regions, the ETOPO1 global relief bathymetry was used \cite{Amante:2009ud}.  Sea-level was set to 0.28 meters above the sea-level in the bathymetric sources.  This is to account primarily for the swelling of the Gulf of Mexico during the summer along with other effects such as upper layer warming, seasonal riverine discharges, and the measured sea level rise \cite{Hope:2013jt}.

\subsubsection{Variable Friction Field} \label{ssub:variable_friction}

\begin{figure}[htb]
    \centering
    \begin{subfigure}[b]{0.46\textwidth}
        \includegraphics[width=\textwidth]{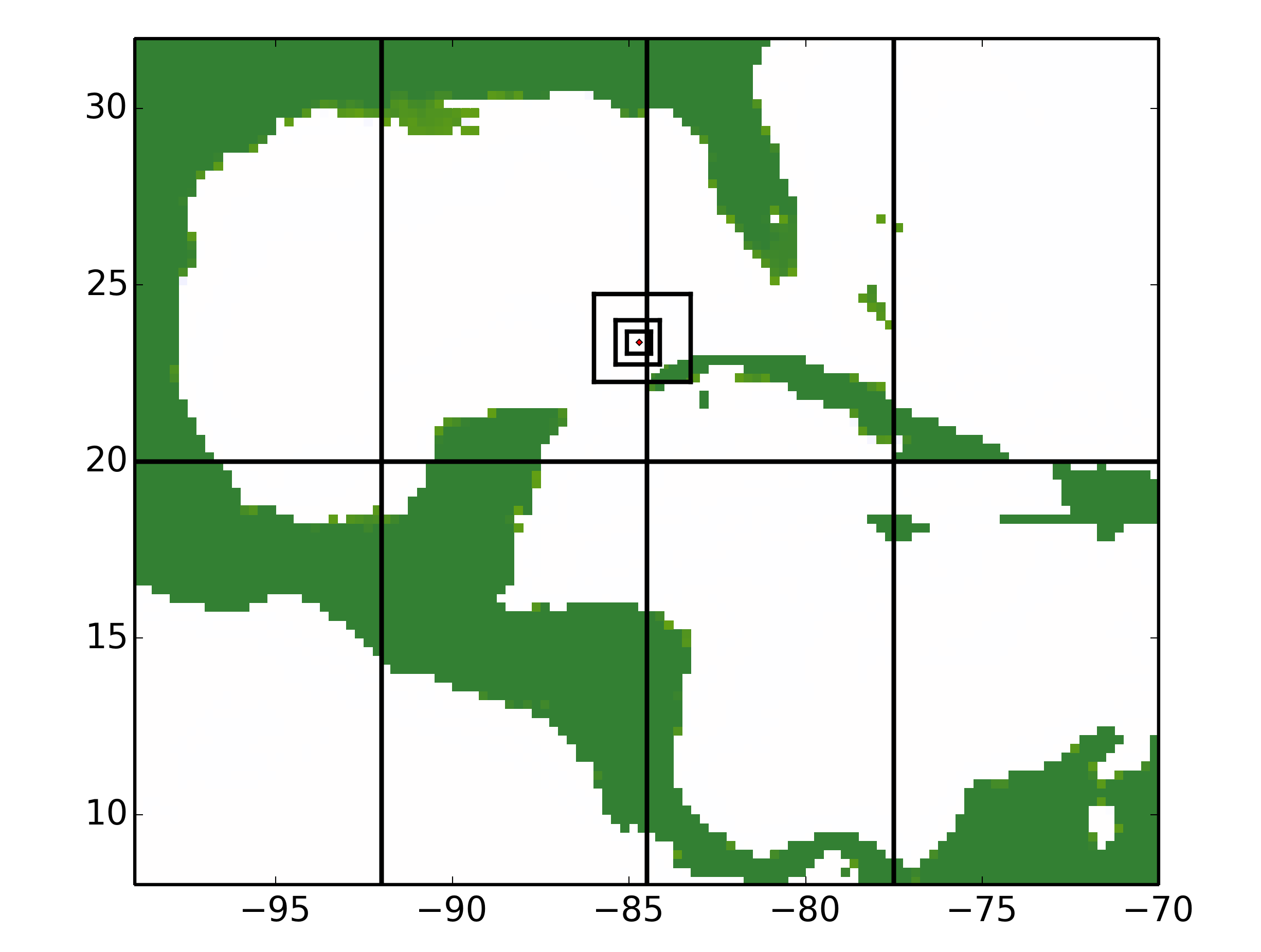}
        \caption{}
        \label{fig:full_domain}
    \end{subfigure}
    \begin{subfigure}[b]{0.5\textwidth}
        \includegraphics[width=\textwidth]{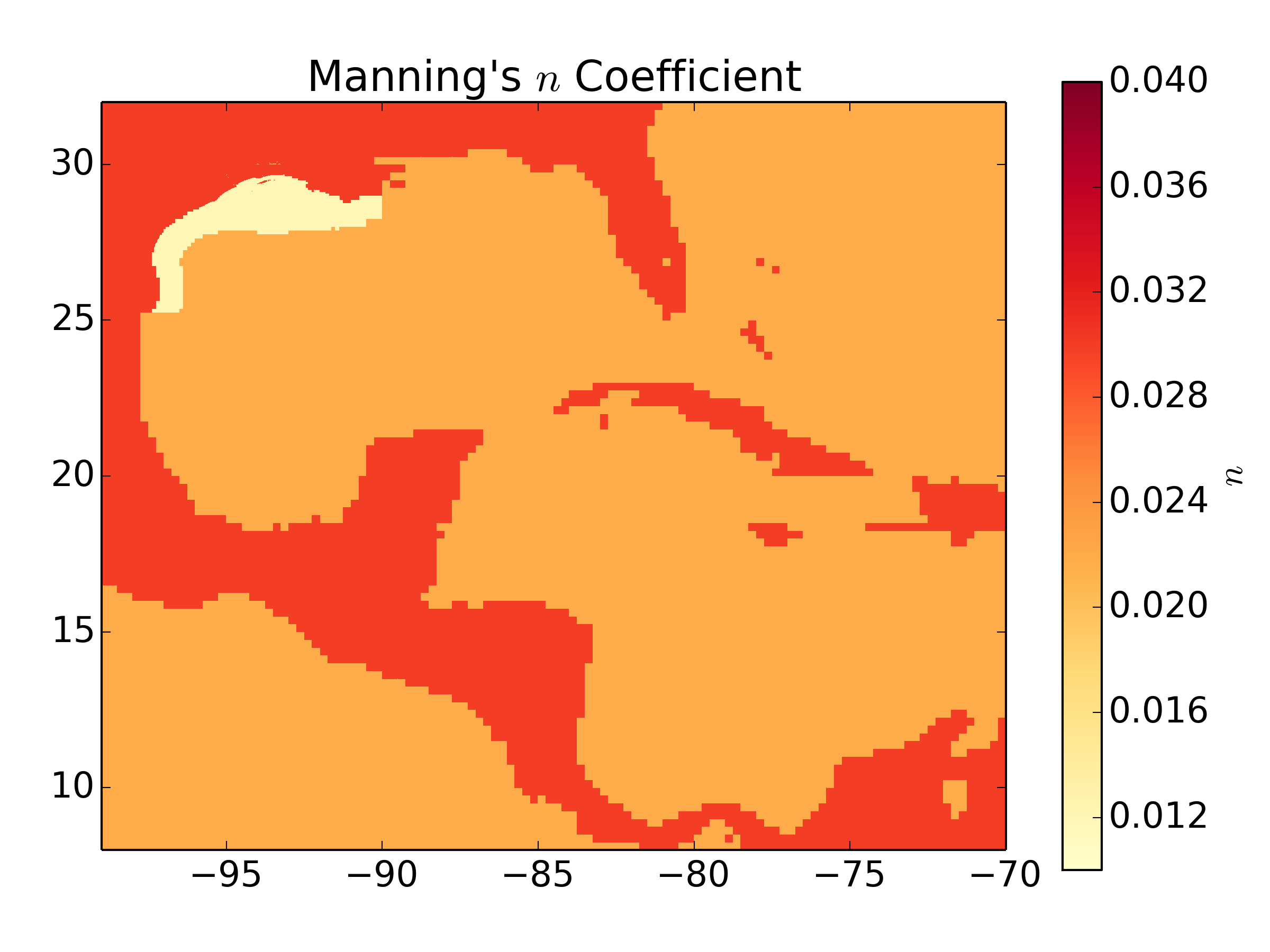}
        \caption{}
        \label{fig:friction}
    \end{subfigure}
    \caption{Figure~\ref{fig:full_domain} represents the entire simulation domain including the starting position of the hurricane and refinement patches.  Figure~\ref{fig:friction} shows the Manning's $n$ coefficients used throughout the domain.}
    \label{fig:geoclaw_domain_characterisitics}
\end{figure}

As was mentioned previously, variable friction can be important to take into account in storm surge simulations.  For the \geoclaw results presented, a simple variable friction field was specified based on bathymetry contours relative to sea-level.  For regions which were initially above sea-level, a Manning's $n$ coefficient of $0.030$ was used and $0.022$ for regions initially below sea-level.  The one exception to this is in the Louisiana-Texas shelf region (defined as $(25.25^\circ \text{N}, 98^\circ \text{W}) \times (30^\circ \text{N}, 90^\circ \text{W})$ in the simulation) where additional contours at $5$ and $200$ meters depth were used such that
\[
    n = \left \{ \begin{aligned}
        &0.030 & ~~~~~& \text{if}~b > \eta_\text{sea-level} - 10 \\
        &0.012 & ~~~~~& \text{if}~\eta_\text{sea-level} - 5 > b > \eta_\text{sea-level}-200 \\
        &0.022 & ~~~~~& \text{if}~ \eta_\text{sea-level}-200 > b
    \end{aligned} \right .
\]
where $\eta_\text{sea-level} = 0.28$ meters as mentioned before.  Figure~\ref{fig:friction} illustrates the resulting values of $n$ through the entire domain.  The seemingly counter-intuitive change in friction on the shelf was suggested in \cite{Kennedy:2011kt} as a means to explain the anomalously large forerunner observed from Ike.  

It should be noted that the value of $n$ can change due to refinement.  This is currently handled by evaluating the cell depth at the current resolution.  For example, in the near shore region a grid cell could be refined to 4 cells which may not have all been above sea-level independently.  When this is the case, the depth of each new cell would determine its value of $n$.  If the fields were set based on a given resolution grid (for instance from NASA land use data), a consistent averaging or interpolation is needed.

\subsection{Adaptive Grids} \label{sub:amr_ike}

As there was significant setup in the low lying coastal areas east of Galveston, the region of greatest interest, higher refinement was forced in these areas before the storm arrived to capture this effect.  The original domain is defined as the region $8^\circ$ N to $32^\circ$ N and $99^\circ$ W to $70^\circ$ W at $1/4$ degree resolution.  From this coarsest grid 6 more levels of refinement were used with ratios and resolutions detailed in Table~\ref{tab:refinement_ratios}.  As mentioned in Section~\ref{sub:amr}, an optimal refinement ratio in time is used based on the CFL condition on each level.  The refinement criteria included sea-level, speed, and storm based criteria and are specified in Table~\ref{tab:refinement_criteria}.  These values were chosen after doing multiple simulations of both Hurricanes Ike and Irene to determine qualitatively the best values in Table~\ref{tab:refinement_criteria}.   Maximum refinement over the entire domain was limited to level 5 away from the coasts of Louisiana and Texas where the full resolution was allowed.  As mentioned in Section~\ref{sub:amr}, there was an increase in total mass of approximately $ 0.0015 \%$ compared to the original total mass of the domain.

\begin{table}[htb]
    \begin{center}
    \begin{tabular}{c|c|rr}
    \textbf{Level} & \textbf{$r^\ell_{\dx,\dy}$} & \multicolumn{2}{c}{\textbf{Resolution (m)}} \\
    ~ & ~ & Latitude & Longitude \\
    \hline
    1 & ~ & 25250 & 27700 \\
    2 & 2 & 12600 & 13850 \\
    3 & 2 & 6300 & 6925 \\
    4 & 2 & 3150 & 3460 \\
    5 & 6 & 525 & 575 \\
    6 & 4 & 130 & 144 \\
    7 & 4 & 32.9 & 36.1
    \end{tabular}
    \end{center}
    \caption{Refinement ratios $r^\ell_{\dx,\dy}$ and effective resolutions at each level.  Since the grid is defined in latitude-longitude coordinates, the meters are an approximate resolution in each direction.}
    \label{tab:refinement_ratios}
\end{table}

\begin{table}[htb]
    \begin{center}
    \begin{tabular}{l|c}
    \textbf{Criteria} & \textbf{Tolerances} \\
    \hline
       $T_{\text{wave}}$  & 1 m\\
       $T_{\text{speed}}$  & $[1, 2, 3, 4]$ m/s\\
       $T_{\text{r}}$  & $[60, 40, 20]$ km\\
       $T_{\text{wind}}$ & $[20, 40, 60]$ m/s\\
    \end{tabular}
    \end{center}
    \caption{Refinement criteria tolerances for the sea-surface height $T_{\text{wave}}$, the water speed $T_{\text{speed}}$, the radial distance from the eye $T_{\text{r}}$, and the wind speed $T_{\text{wind}}$.  The lists of tolerances correspond to level criteria, i.e. the first entry is the tolerance for moving from level 1 to level 2.}
    \label{tab:refinement_criteria}
\end{table}

\subsection{Simulation Results} \label{sub:results}

Comparisons of the sea-surface and currents produced by \geoclaw and \adcirc at times before and after land-fall are shown in some of the relevant regions, in particular on the Texas-Louisiana shelf and near Galveston Island and Galveston Bay.  For clarity of the comparison the \geoclaw grid patches were not drawn on the plots.  

\begin{figure}[htb]
    \centering
    \includegraphics[width=0.49\textwidth]{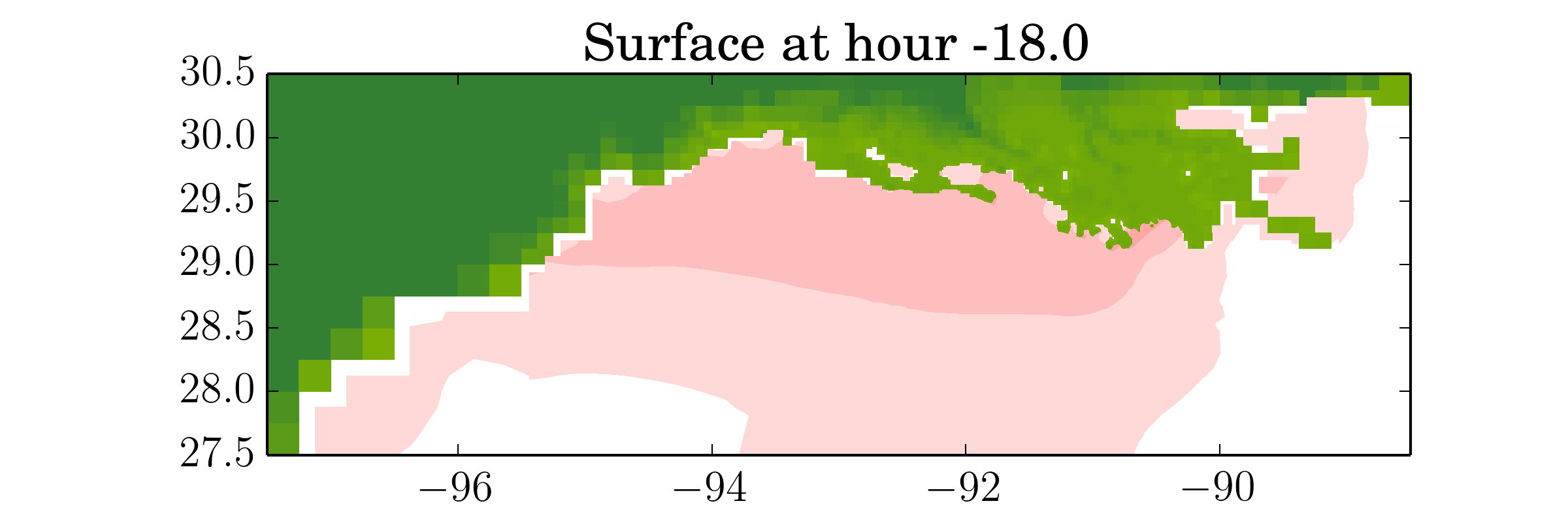} 
    \includegraphics[width=0.49\textwidth]{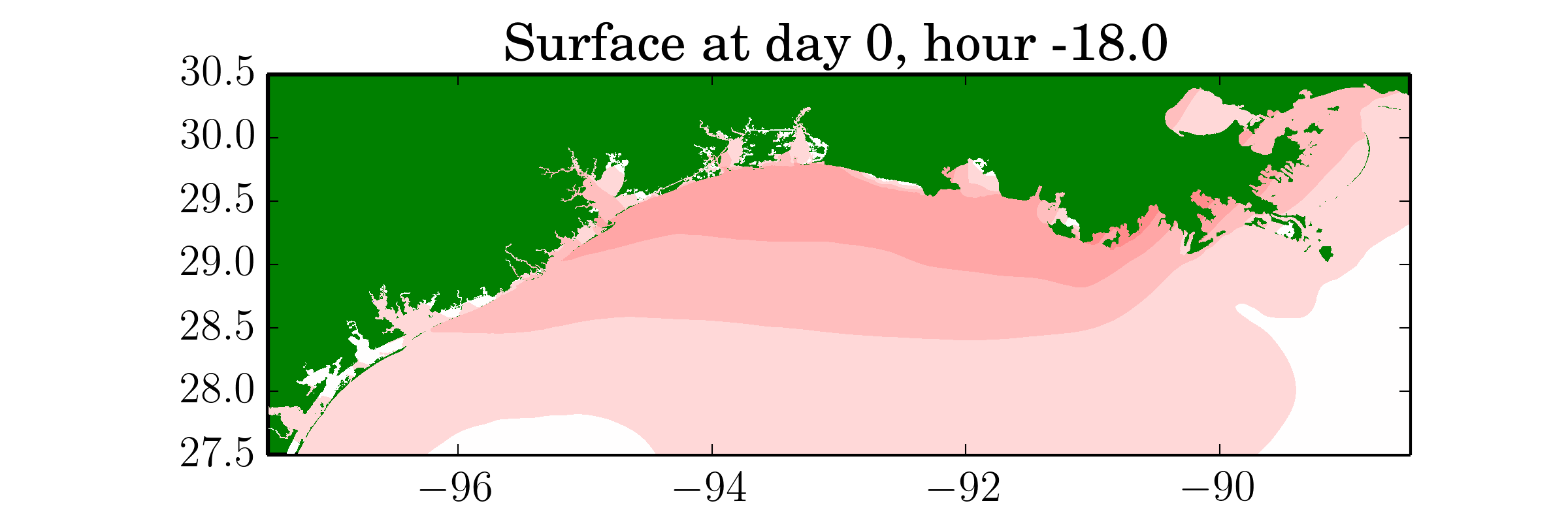} \\
    \includegraphics[width=0.49\textwidth]{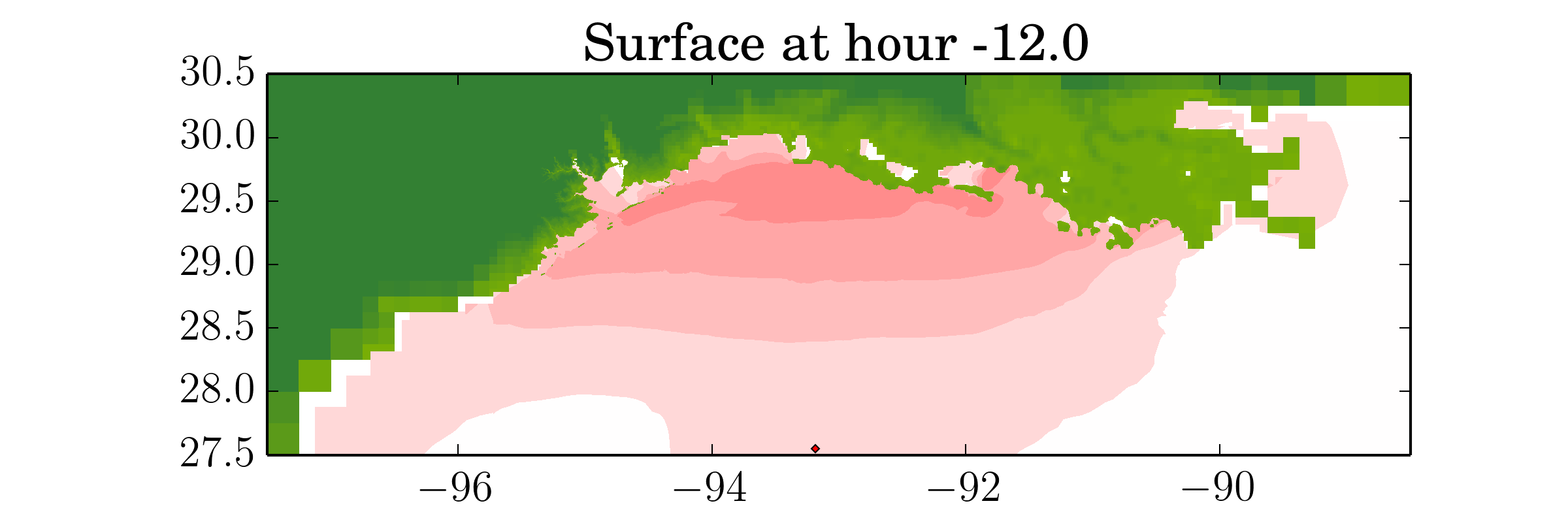} 
    \includegraphics[width=0.49\textwidth]{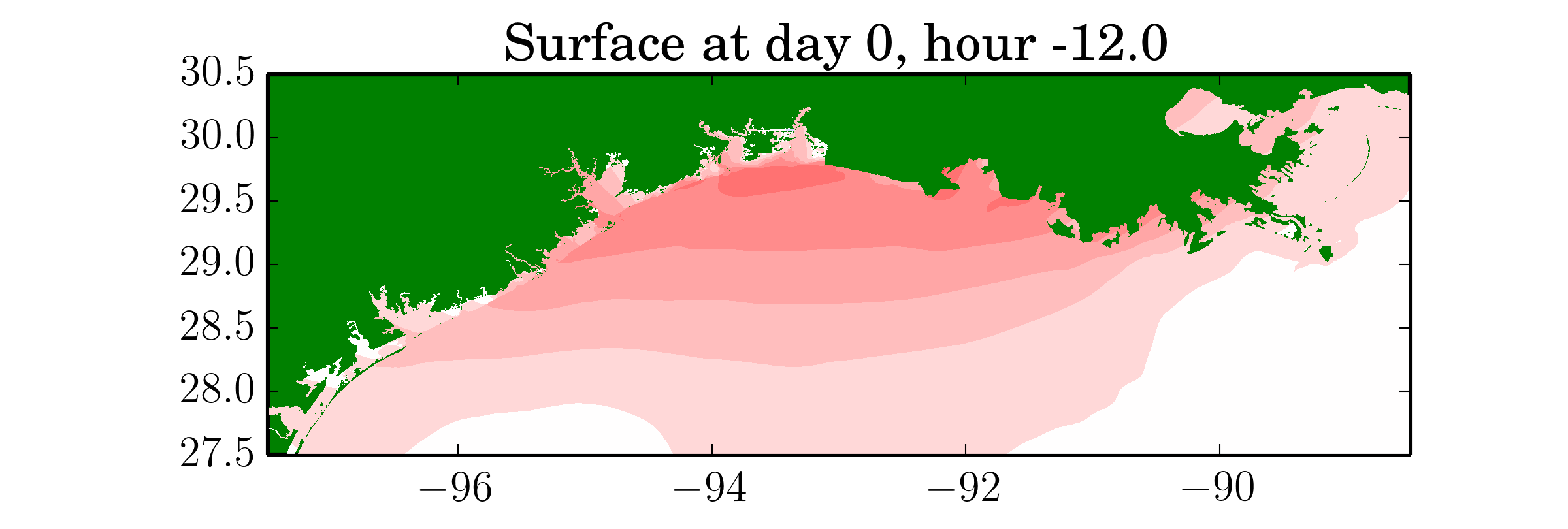} \\
    \includegraphics[width=0.49\textwidth]{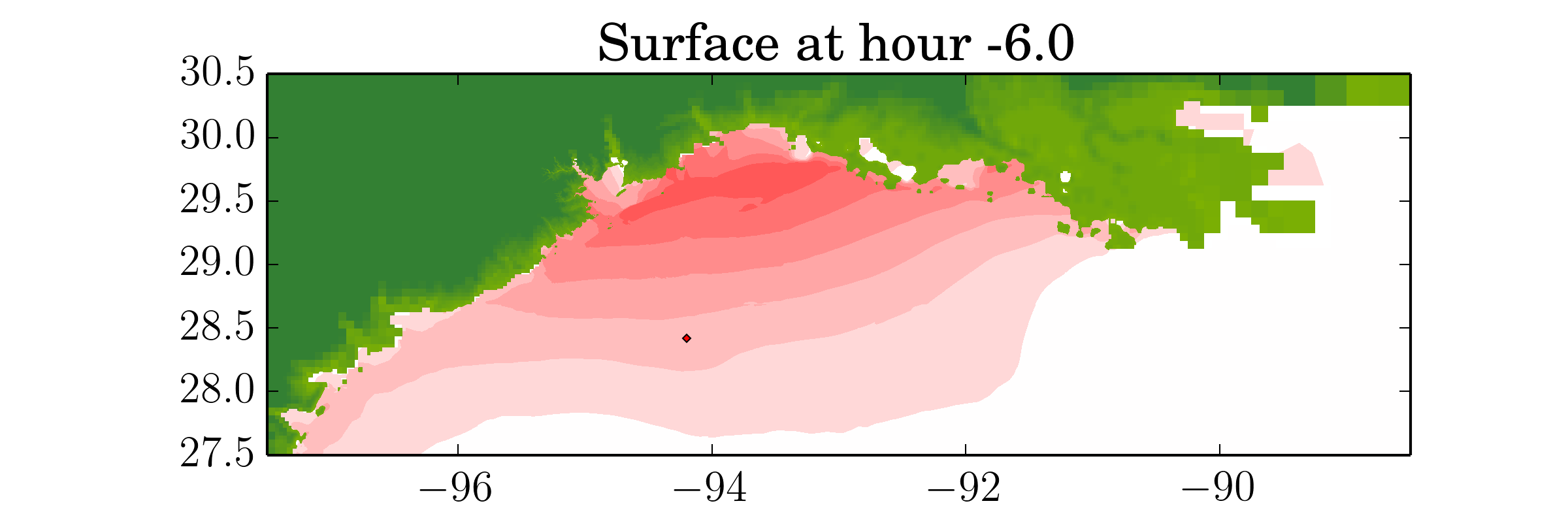} 
    \includegraphics[width=0.49\textwidth]{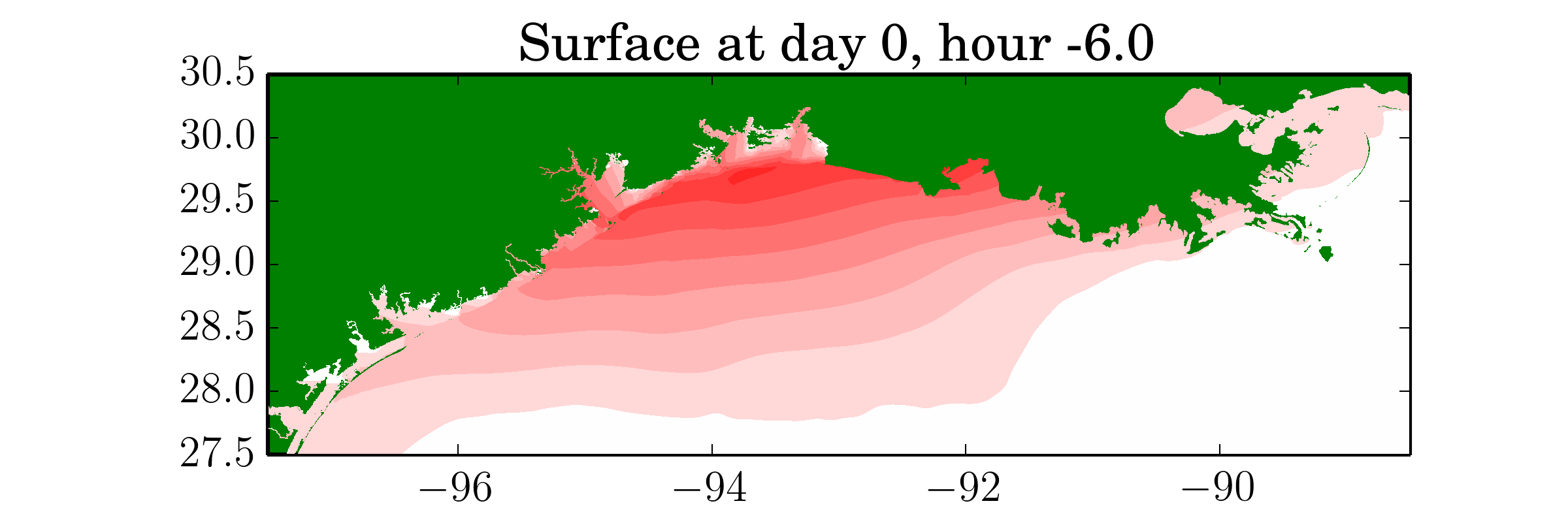} \\
    \includegraphics[width=0.49\textwidth]{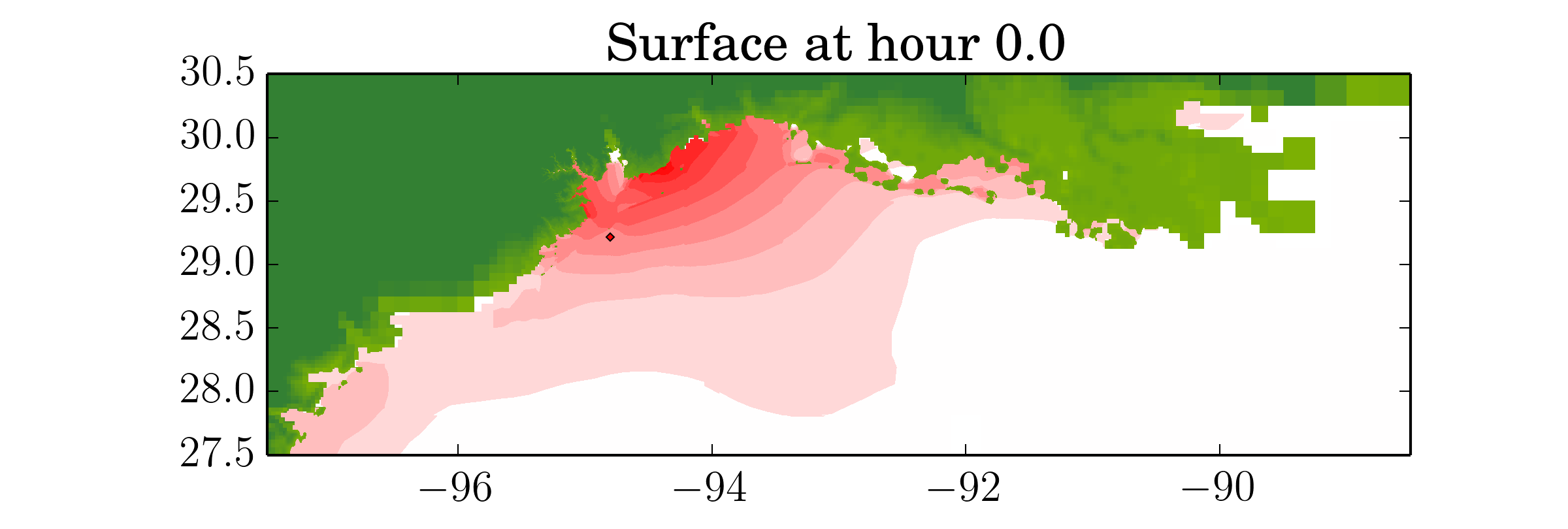} 
    \includegraphics[width=0.49\textwidth]{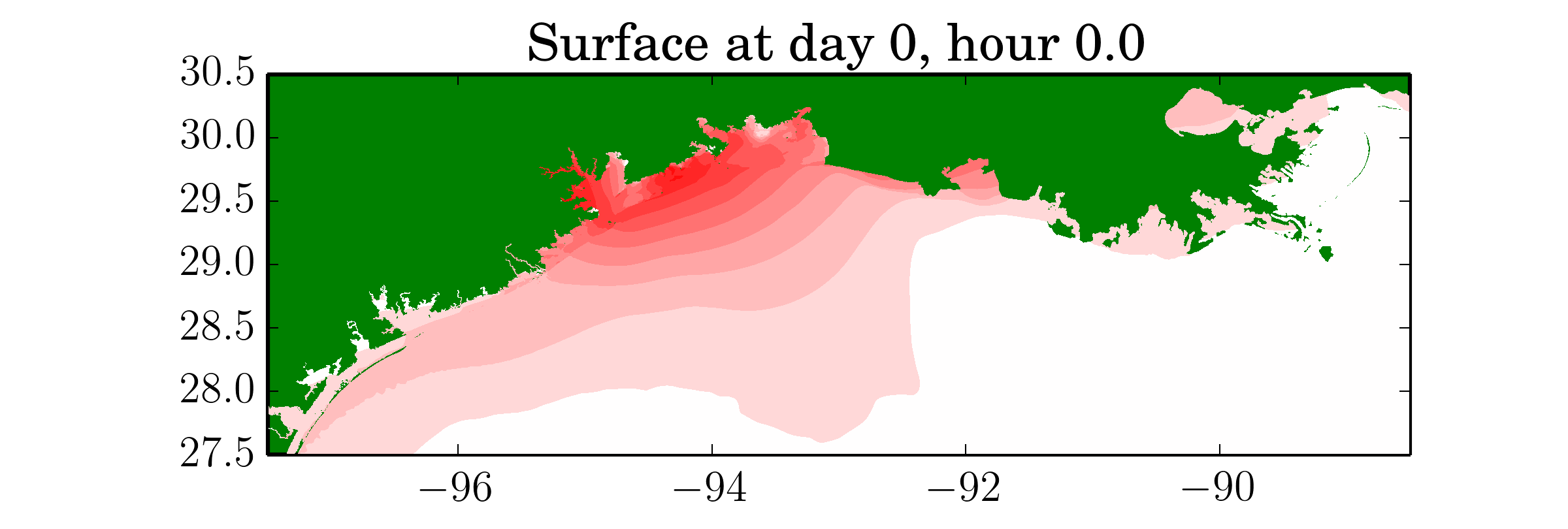} \\
    \includegraphics[width=0.49\textwidth]{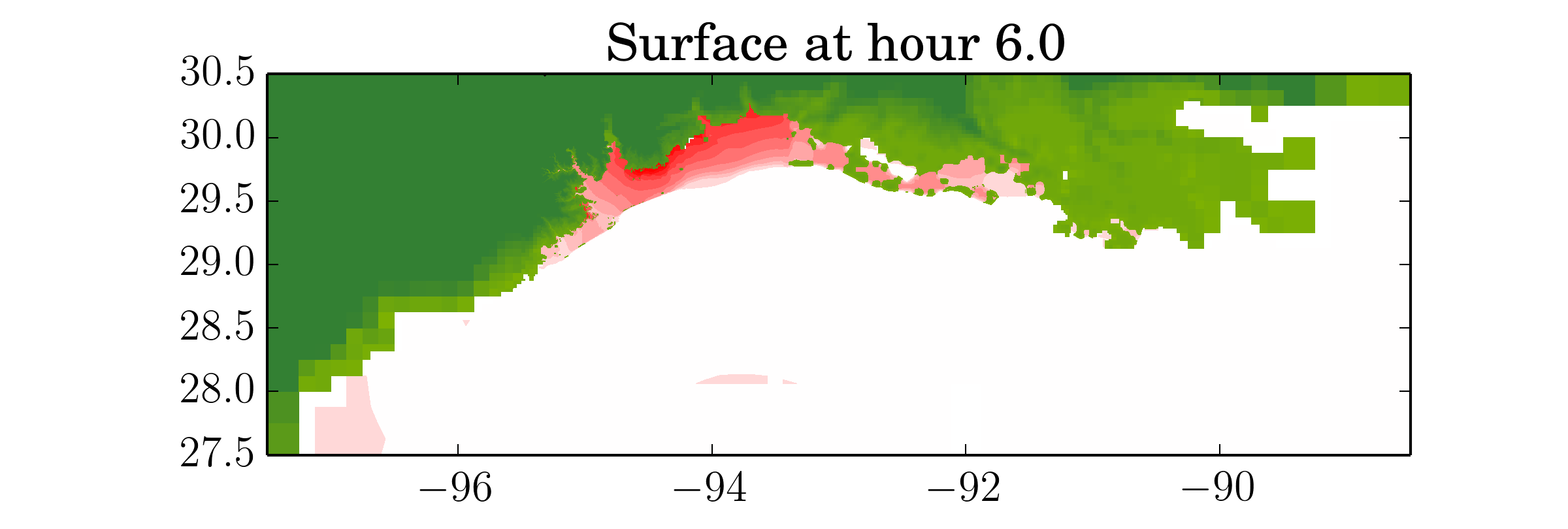} 
    \includegraphics[width=0.49\textwidth]{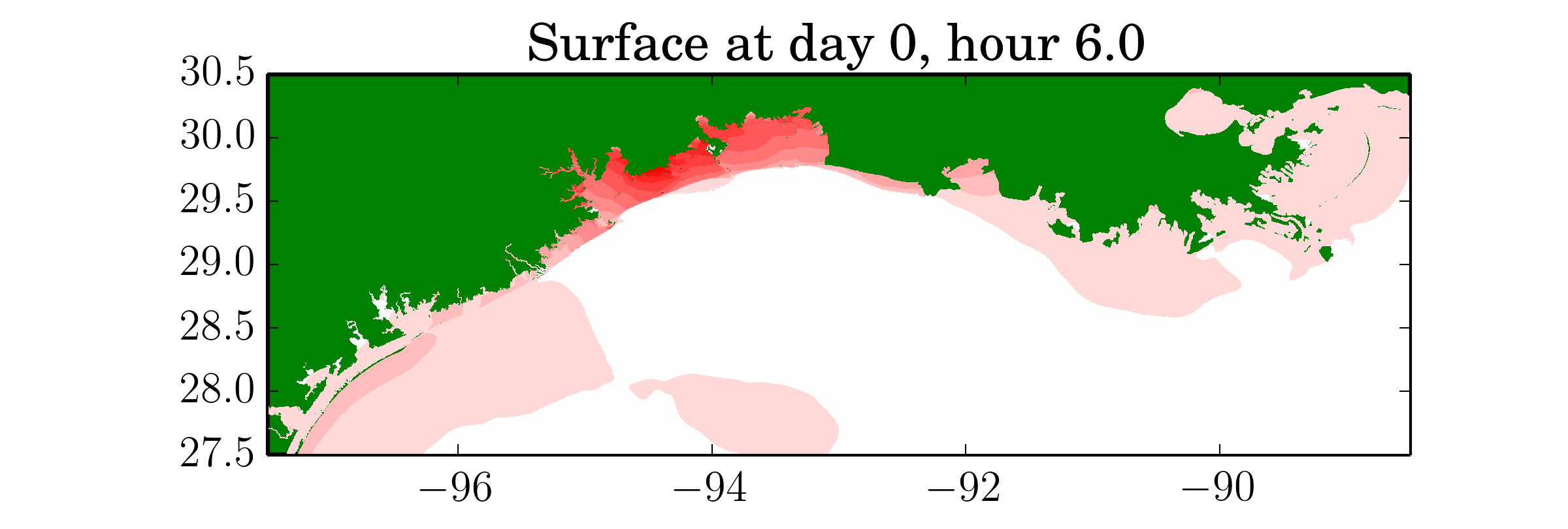} \\
    \includegraphics[width=0.49\textwidth]{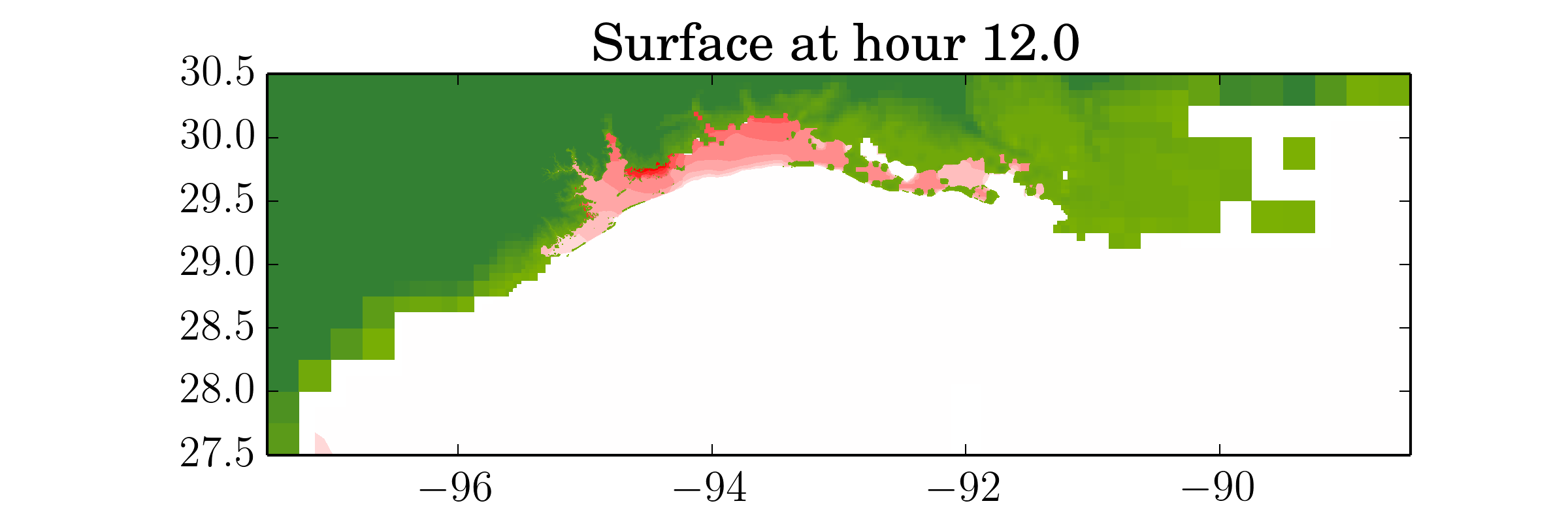} 
    \includegraphics[width=0.49\textwidth]{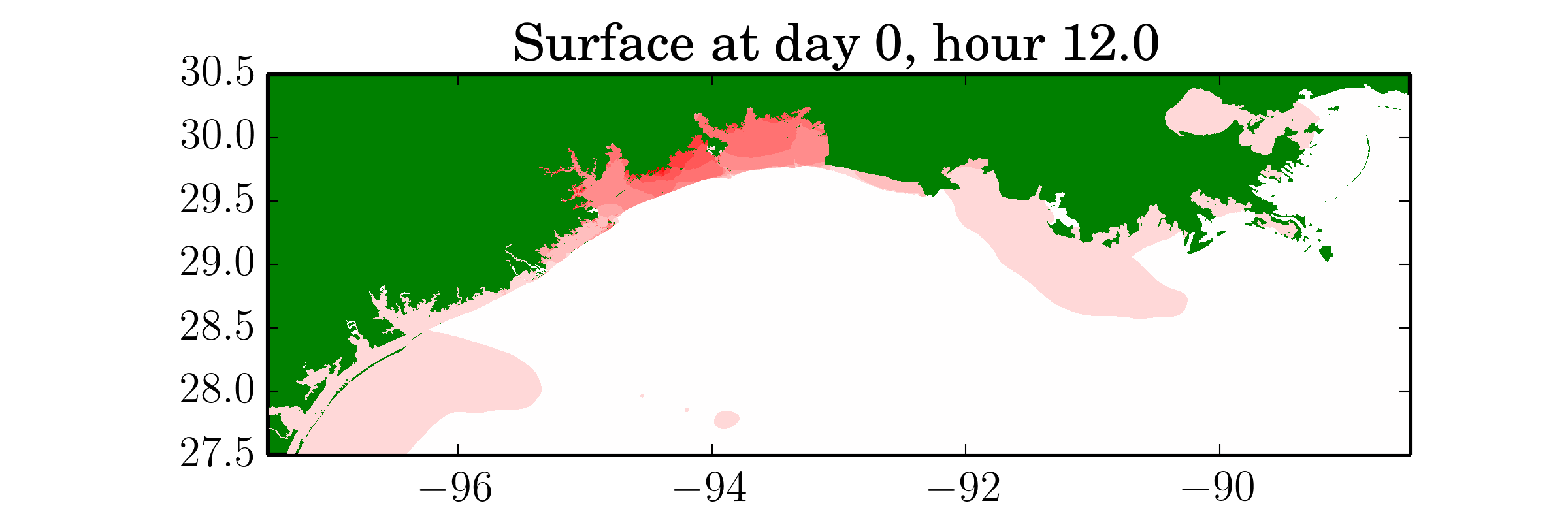} \\
    \includegraphics[width=0.90\textwidth]{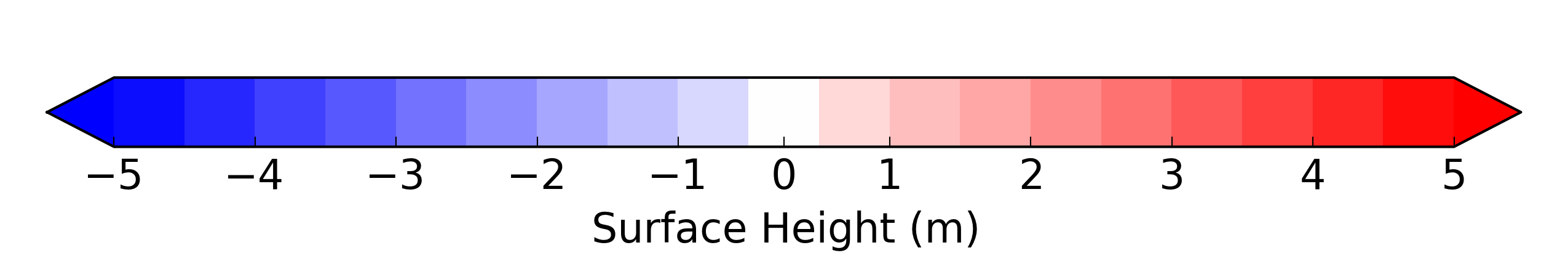}
    \caption{Surface elevations on the Texas-Louisiana shelf produced by \geoclaw (on the left) and \adcirc (on the right) before and after Hurricane Ike makes landfall.}
    \label{fig:latex_shelf_surface}
\end{figure}

\begin{figure}[htb]
    \centering
    \includegraphics[width=0.49\textwidth]{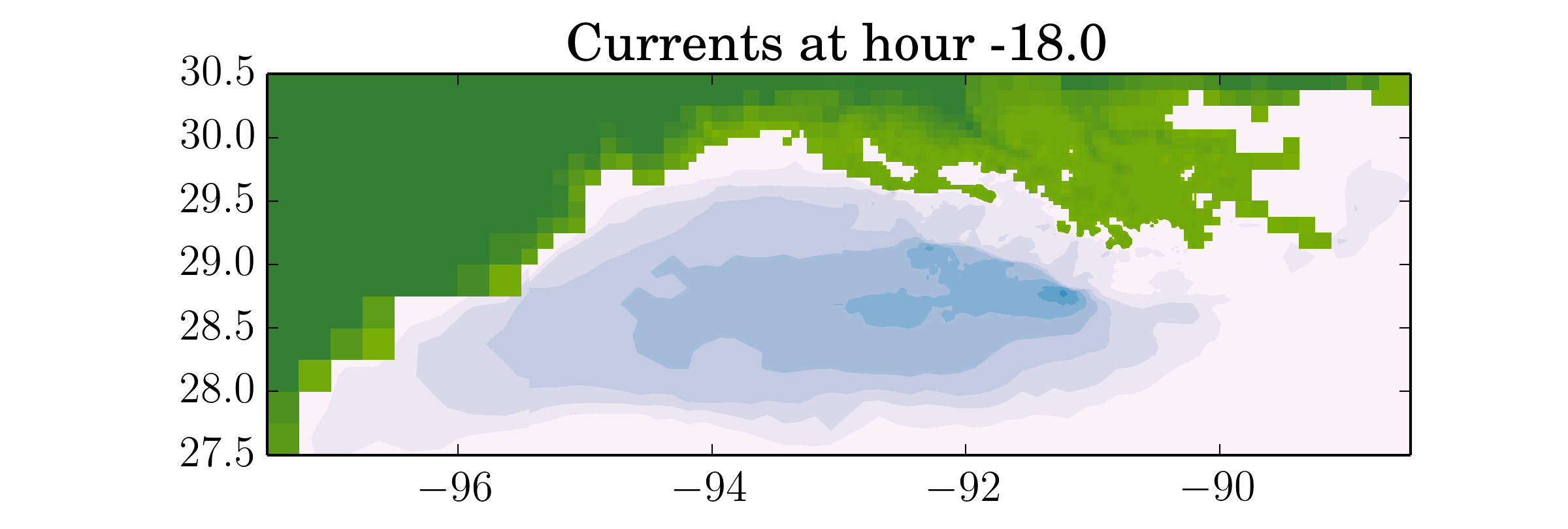}
    \includegraphics[width=0.49\textwidth]{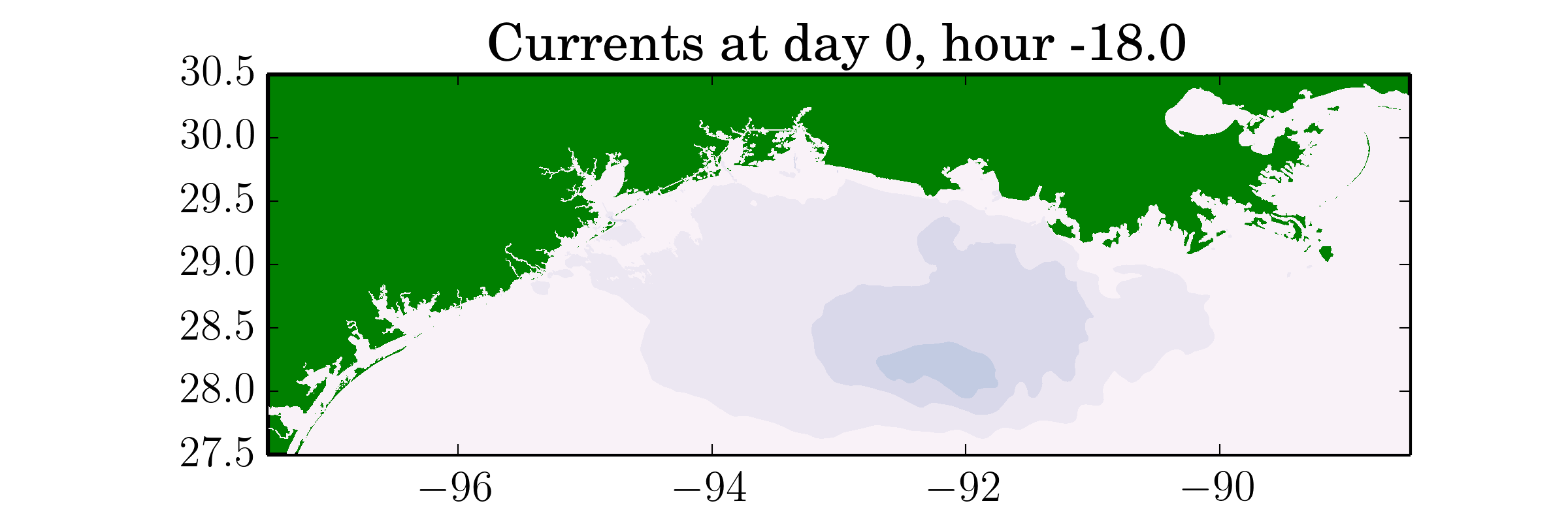} \\ 
    \includegraphics[width=0.49\textwidth]{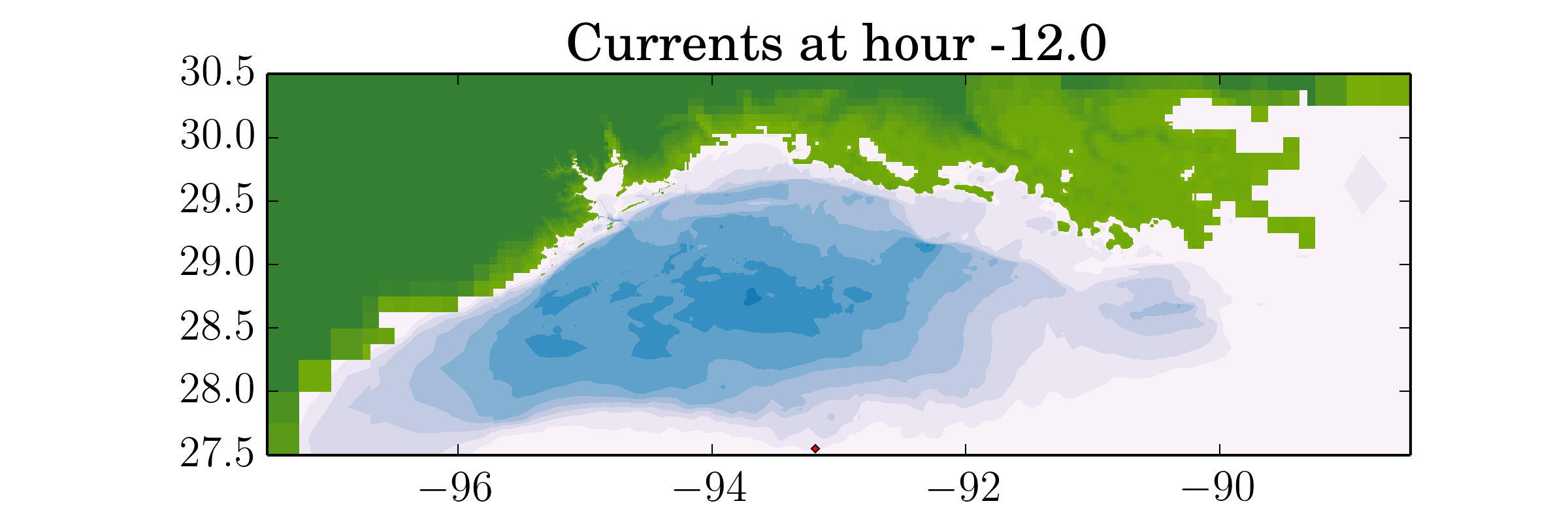}
    \includegraphics[width=0.49\textwidth]{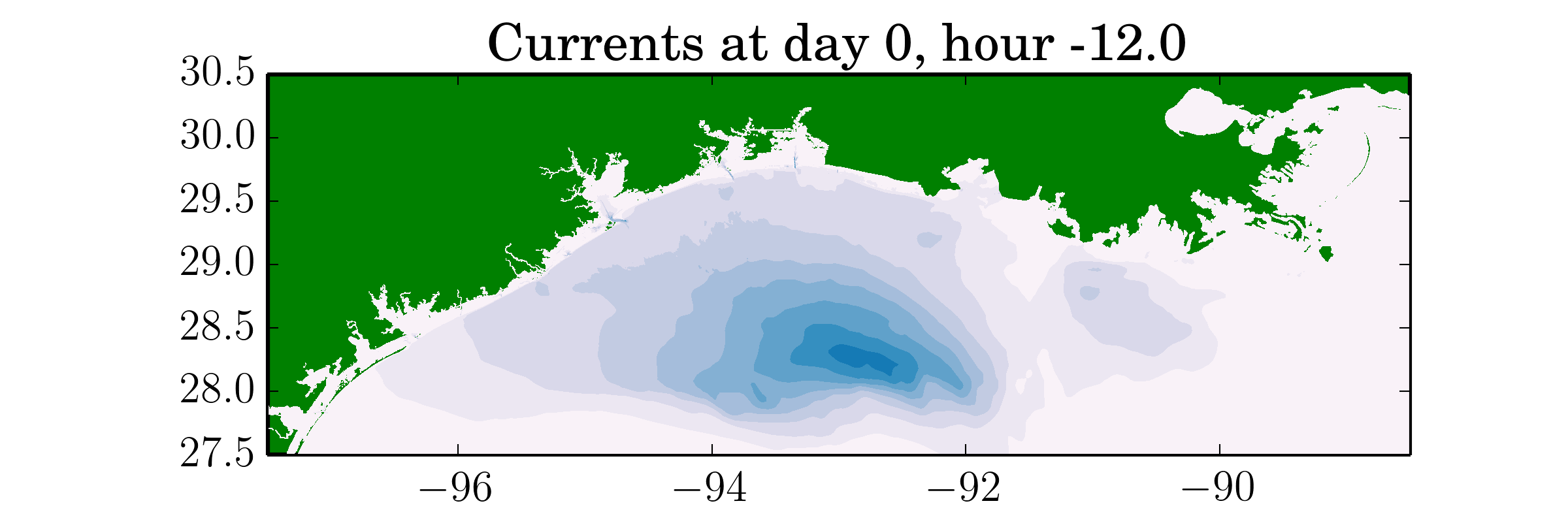} \\
    \includegraphics[width=0.49\textwidth]{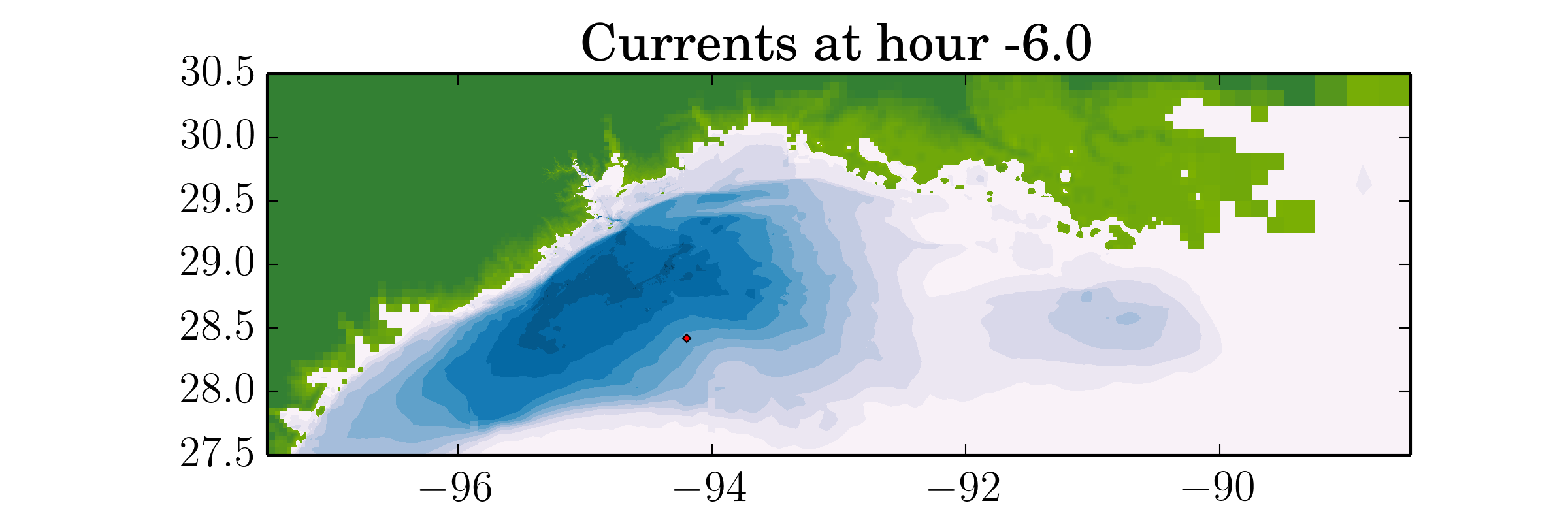}
    \includegraphics[width=0.49\textwidth]{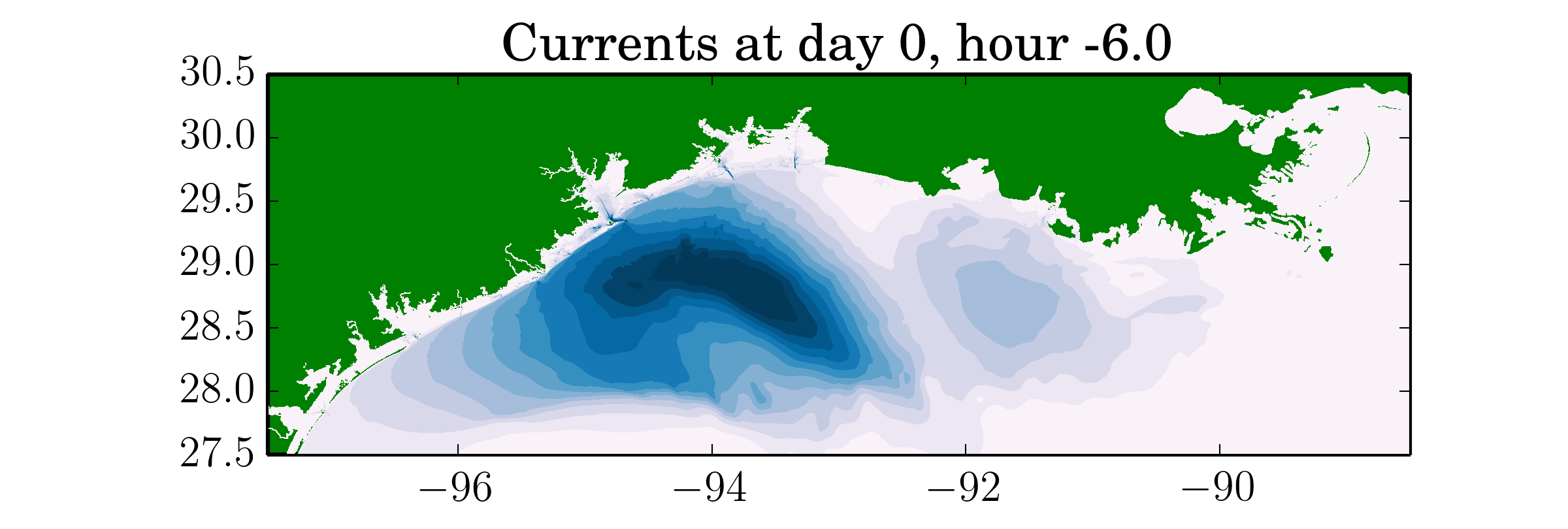} \\
    \includegraphics[width=0.49\textwidth]{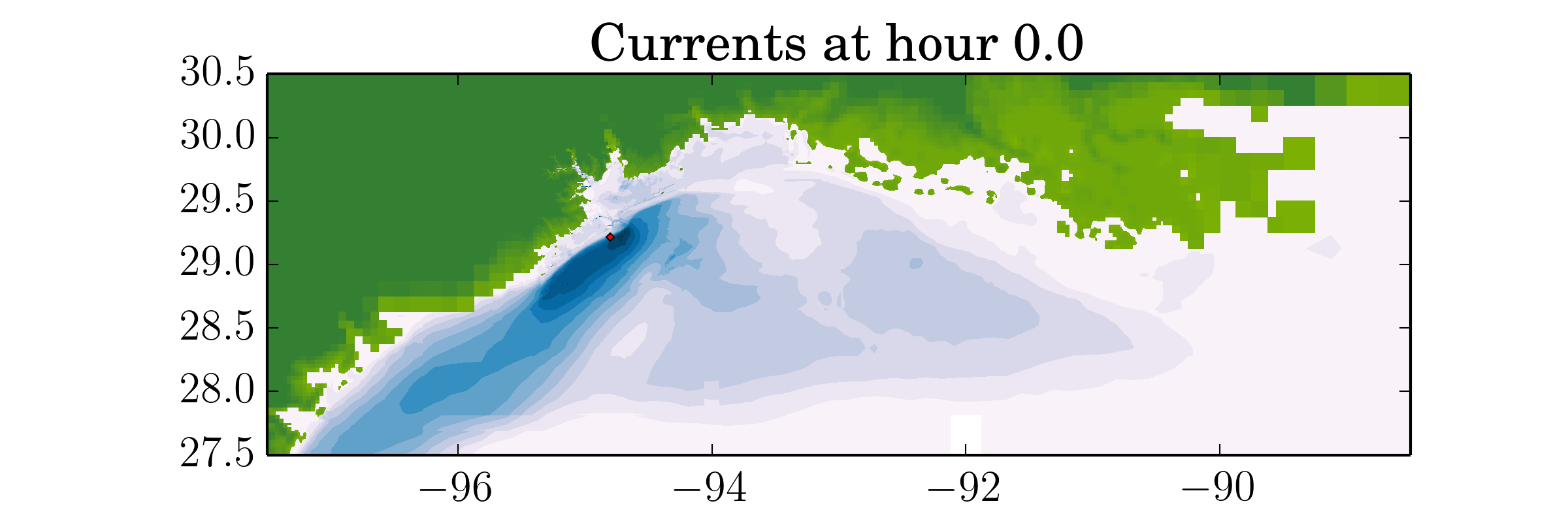}
    \includegraphics[width=0.49\textwidth]{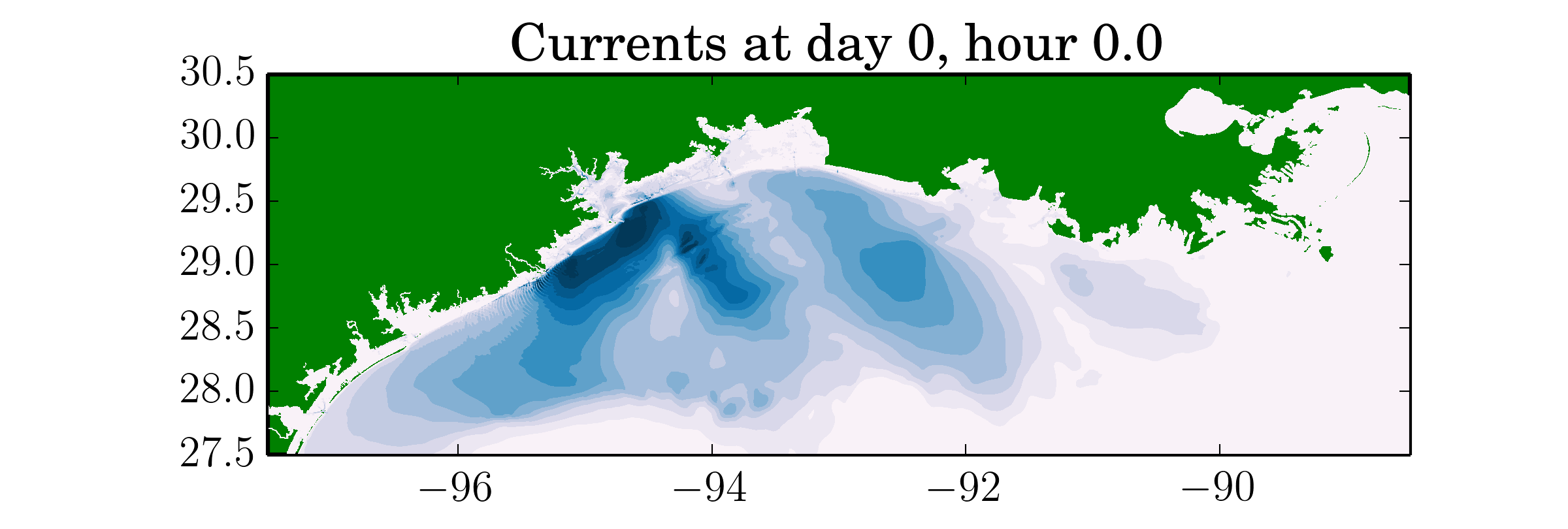} \\
    \includegraphics[width=0.49\textwidth]{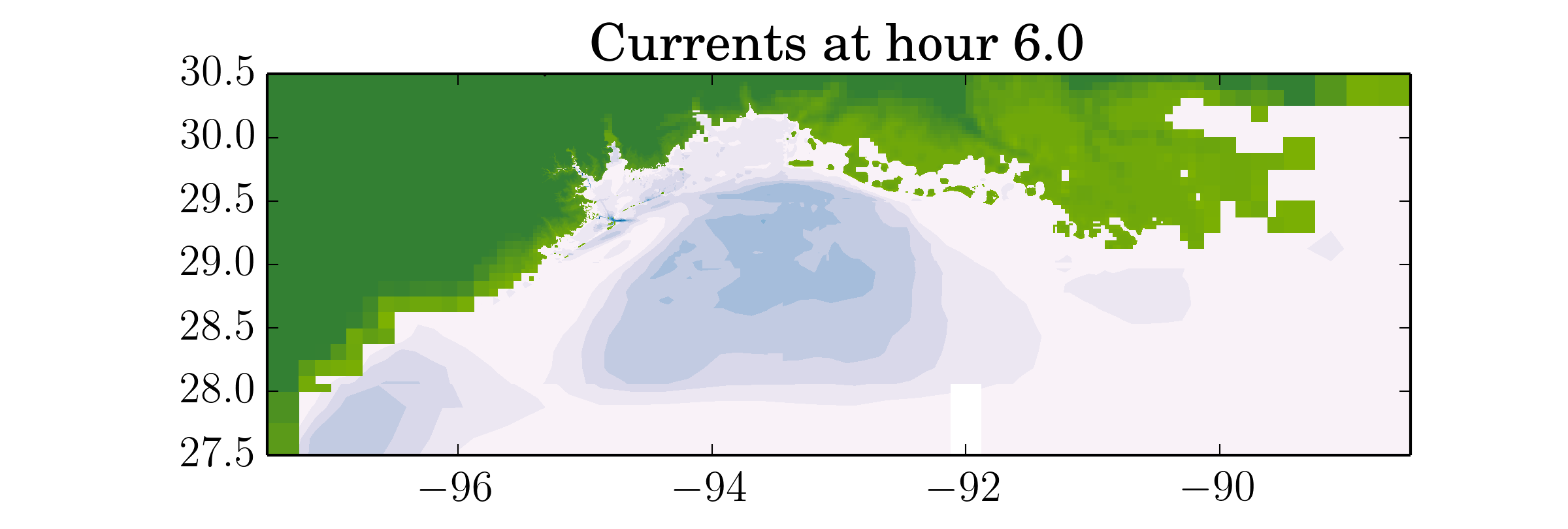}
    \includegraphics[width=0.49\textwidth]{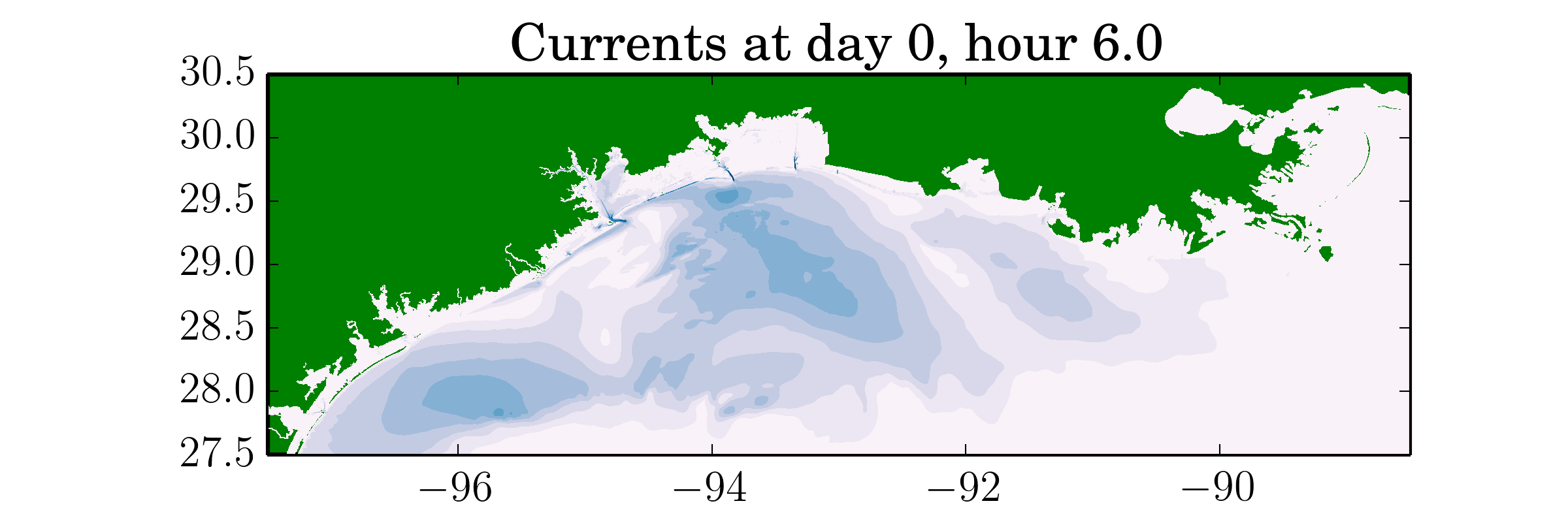} \\
    \includegraphics[width=0.49\textwidth]{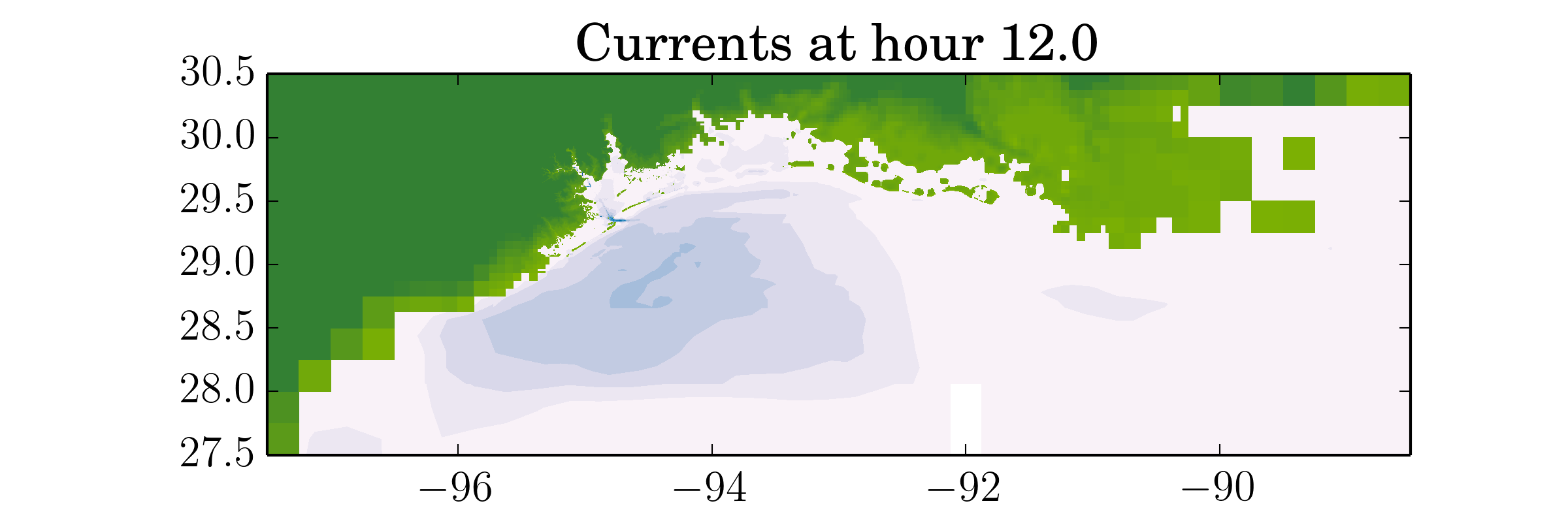}
    \includegraphics[width=0.49\textwidth]{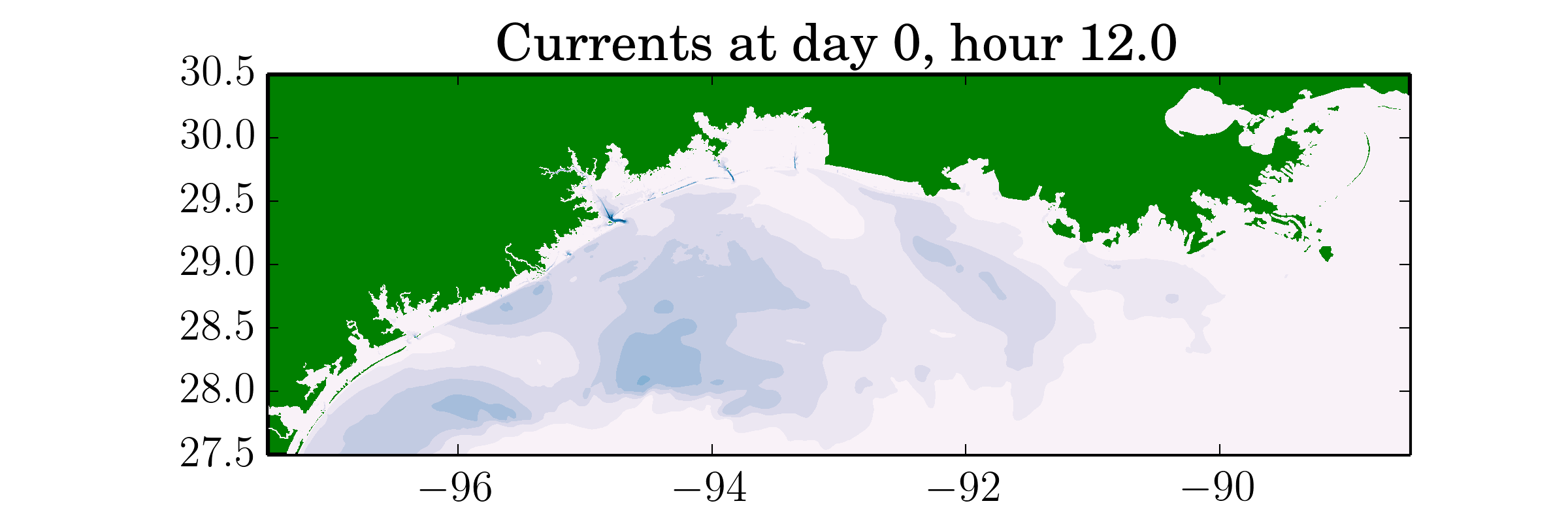} \\
    \includegraphics[width=0.9\textwidth]{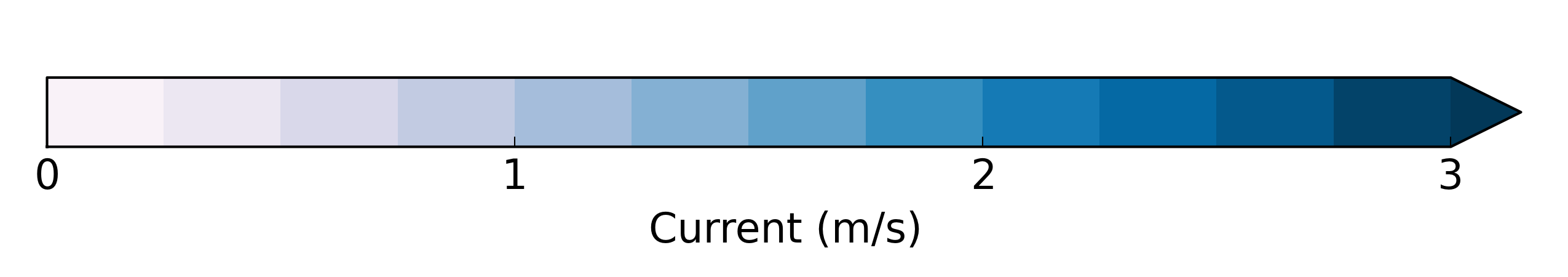}
    \caption{Currents on the Texas-Louisiana shelf produced by \geoclaw (on the left) and \adcirc (on the right) before and after Hurricane Ike makes landfall.}
    \label{fig:latex_shelf_currents}
\end{figure}

Figure~\ref{fig:latex_shelf_surface} shows the surface displacement in the Texas-Louisiana shelf region produced by \geoclaw and \adcirc.  The simulations match well over most of the shelf and only show differences in regions where the \geoclaw simulation did not fully resolve the features of the surge.  Similarly, Figure~\ref{fig:latex_shelf_currents} shows the comparison for the currents in the same region.  Again, in the primary surge region, the simulations match well but away from the regions of interest we see additional structure in the \adcirc simulations, probably due to greater resolution in these areas.  These differences may also be due to the more complex frictional field included in the \adcirc model.

\begin{figure}[htb]
    \centering
    \includegraphics[width=0.46\textwidth]{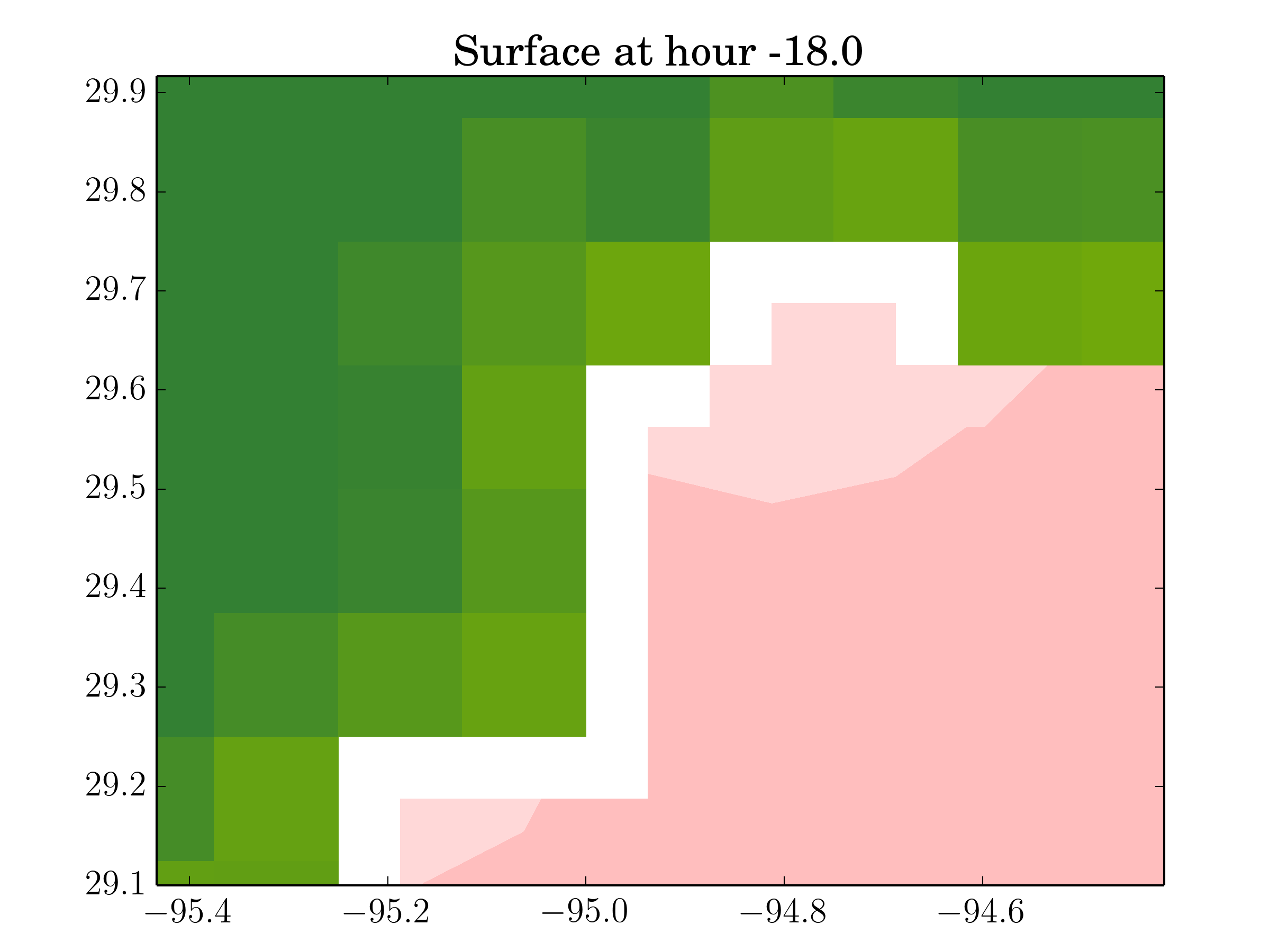} 
    \includegraphics[width=0.46\textwidth]{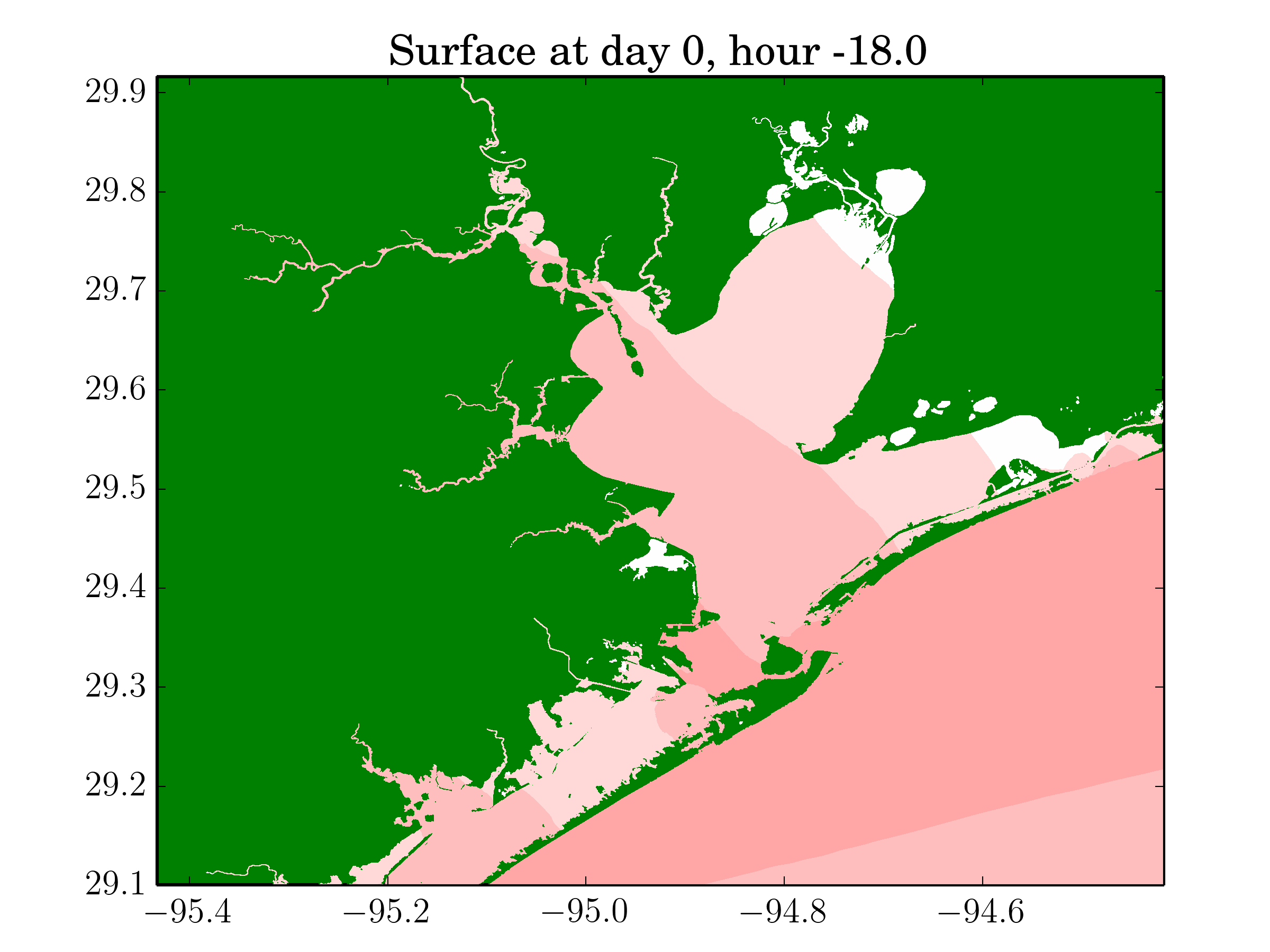} \\
    \includegraphics[width=0.46\textwidth]{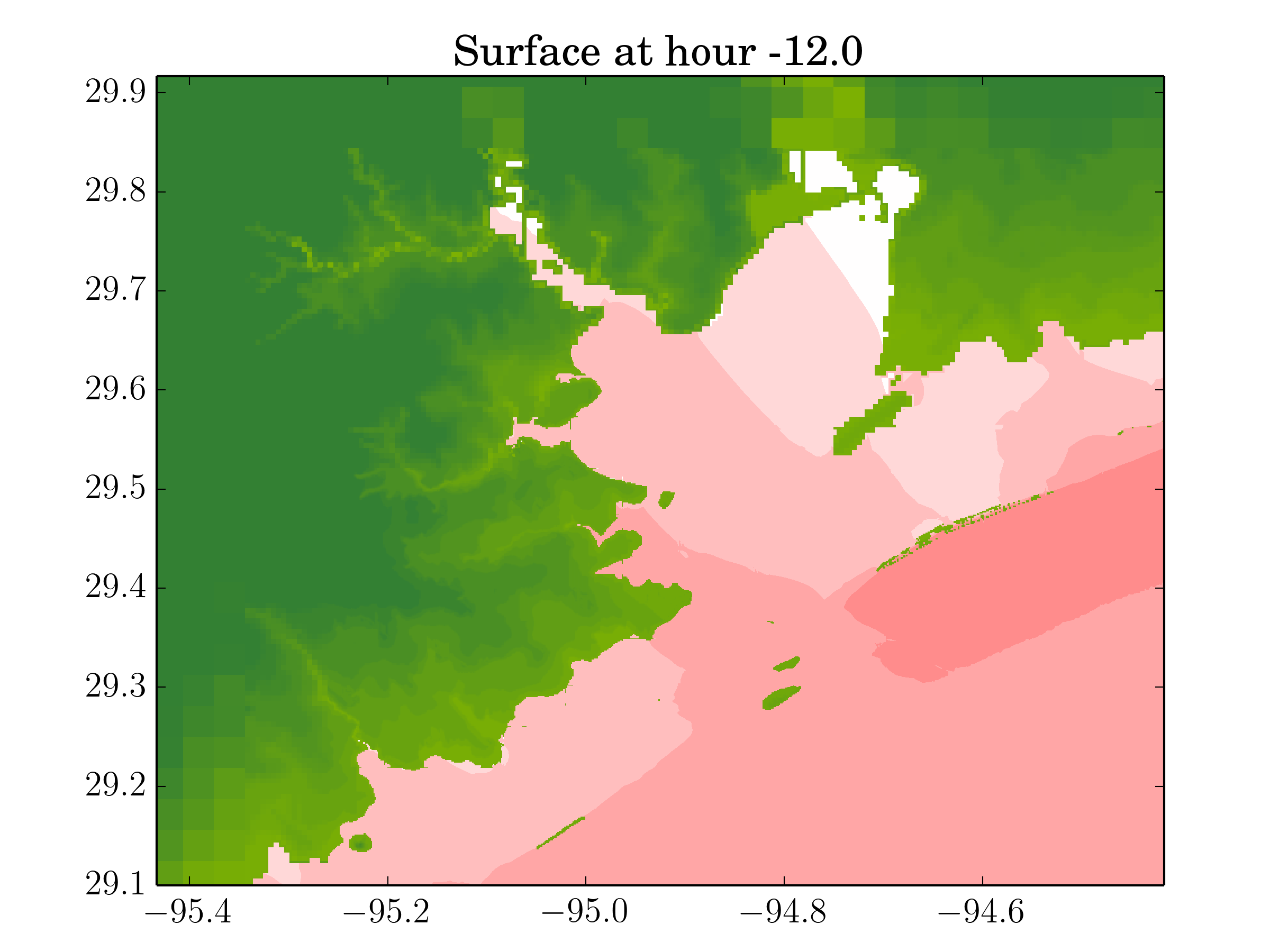} 
    \includegraphics[width=0.46\textwidth]{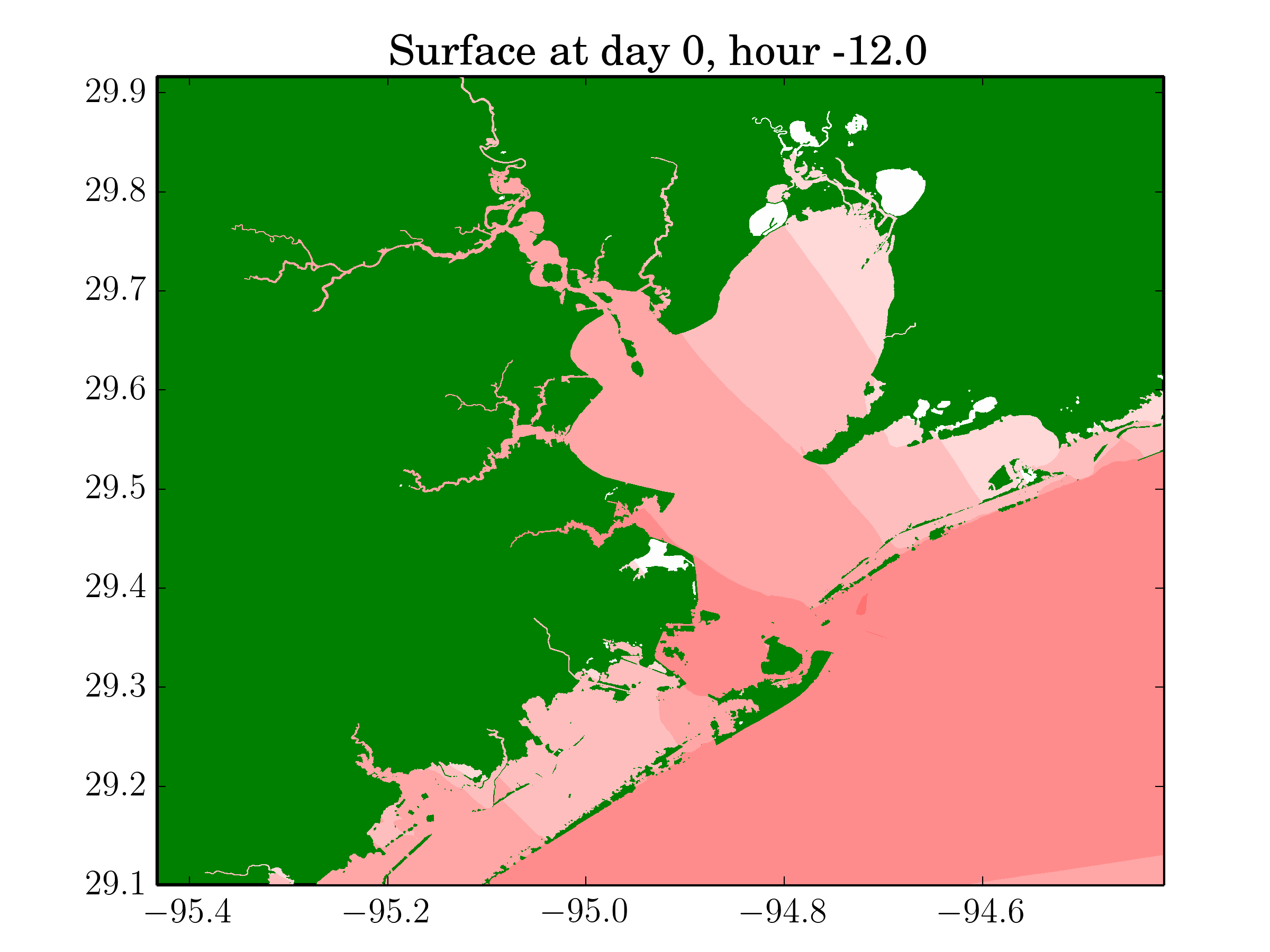} \\
    \includegraphics[width=0.46\textwidth]{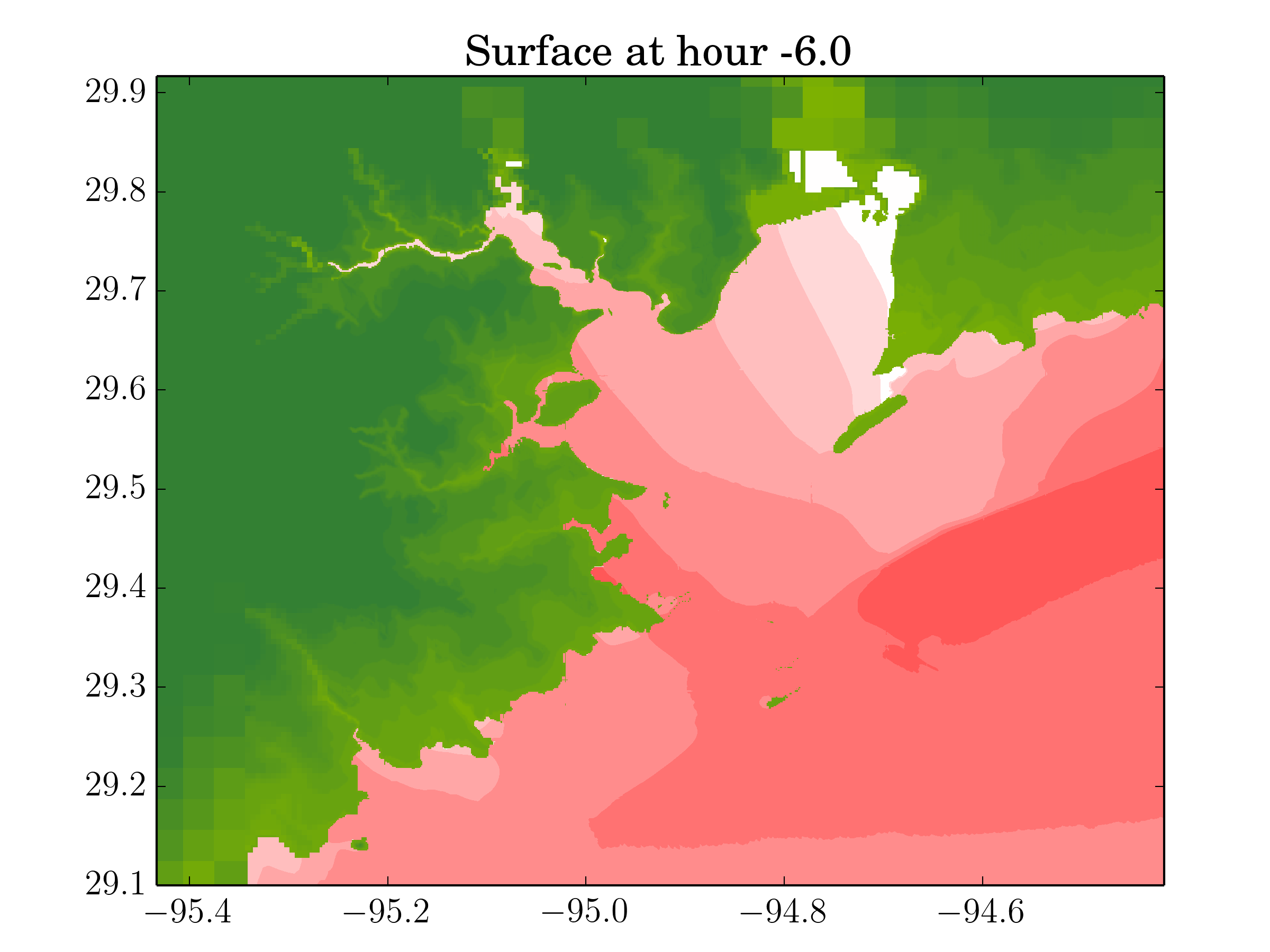} 
    \includegraphics[width=0.46\textwidth]{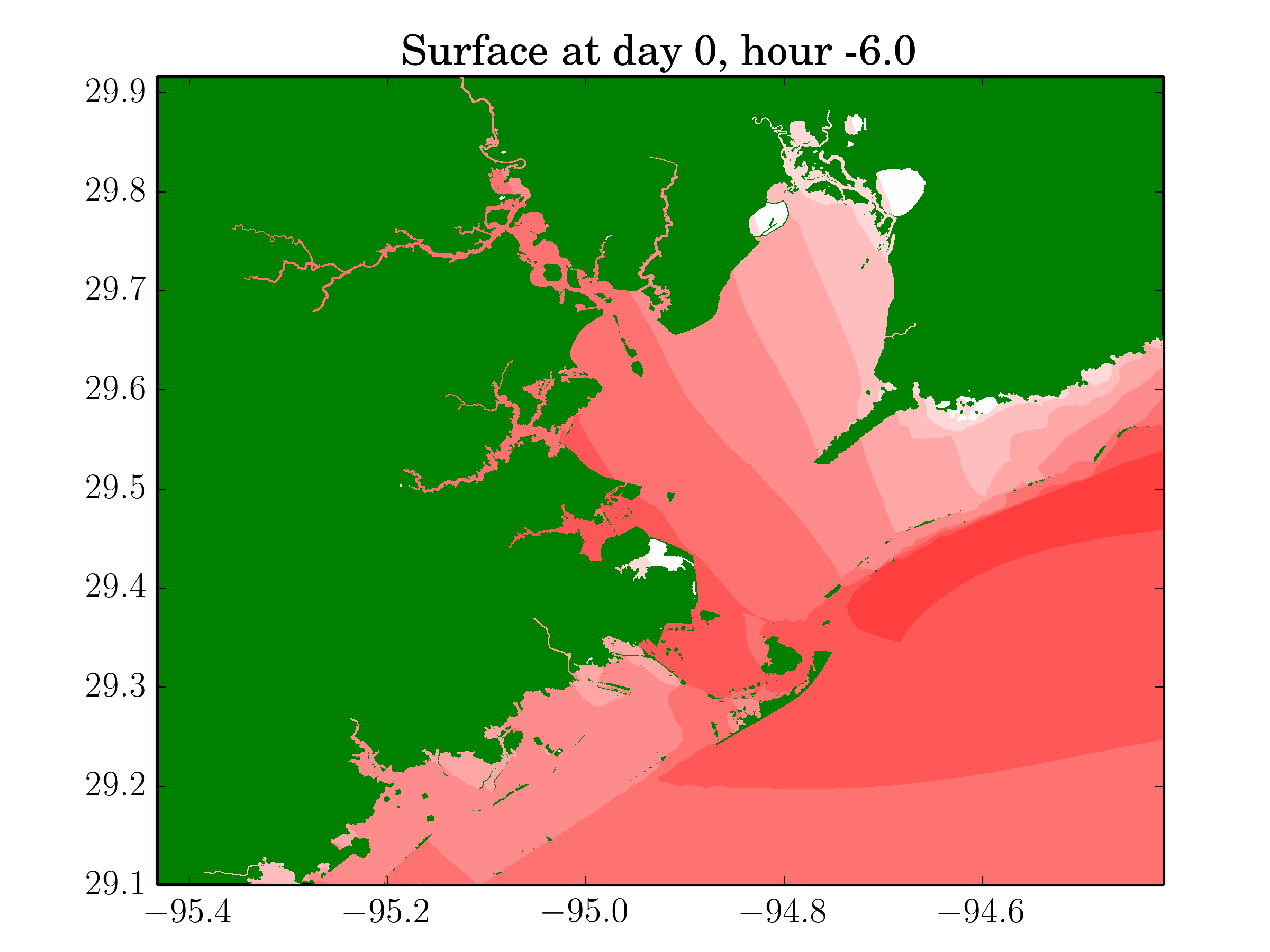} \\
    \includegraphics[width=0.9\textwidth]{surface_cmap.png}
    \caption{Surface elevations in the Galveston Bay region produced by \geoclaw (on the left) and \adcirc (on the right) before Hurricane Ike made landfall.}
    \label{fig:galveston_surface_before}
\end{figure}

\begin{figure}[htb]
    \centering
    \includegraphics[width=0.46\textwidth]{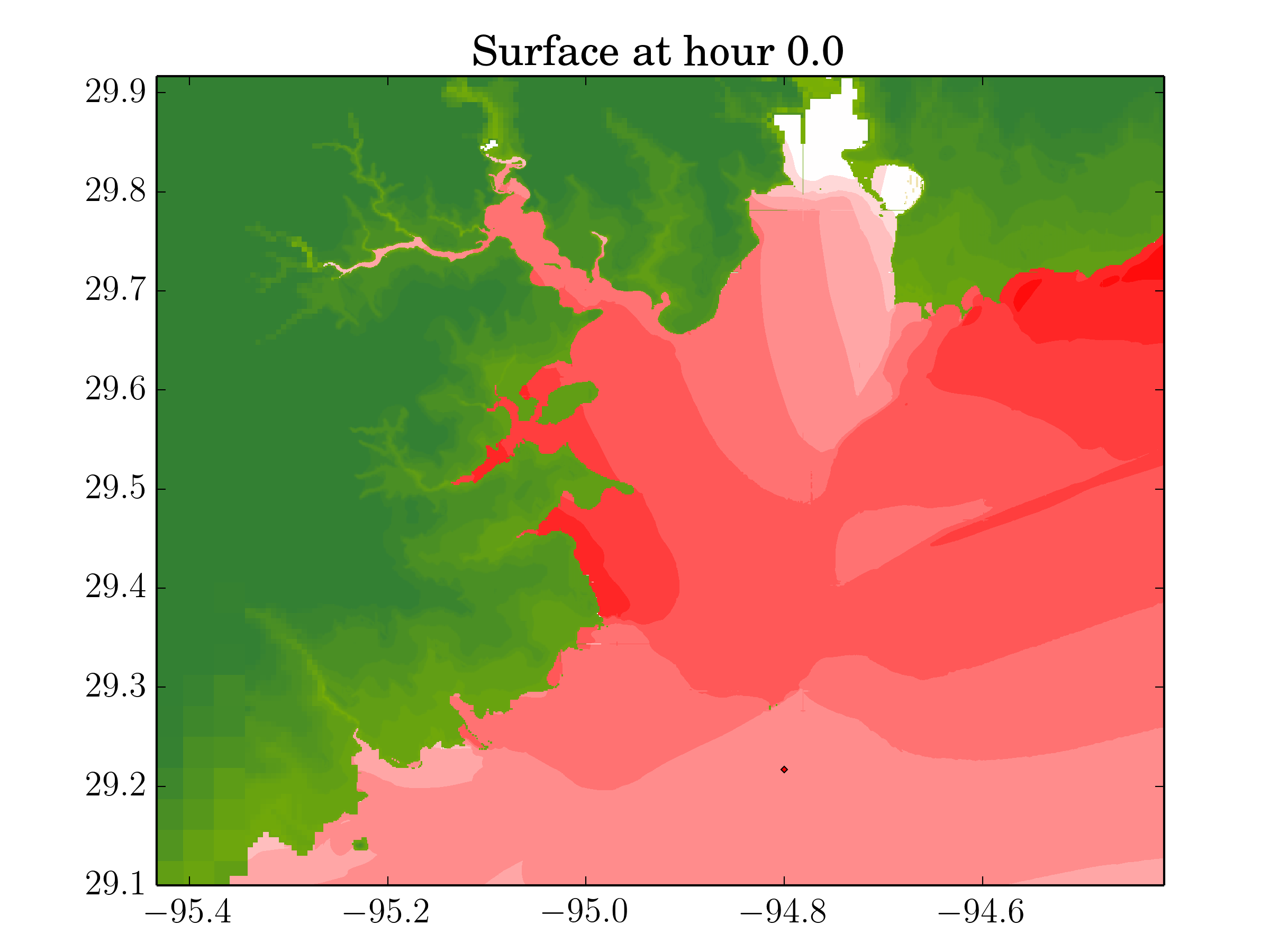} 
    \includegraphics[width=0.46\textwidth]{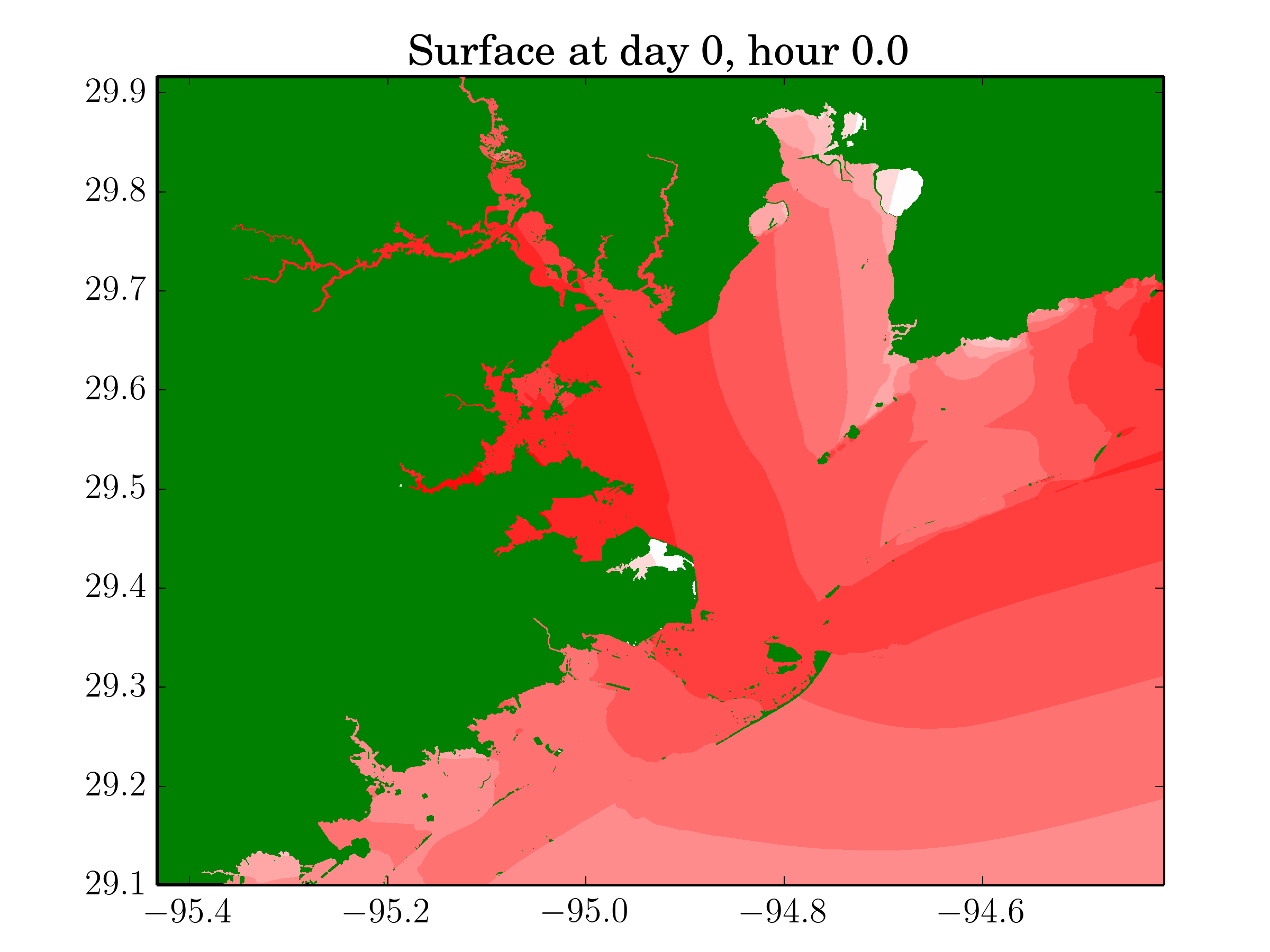} \\
    \includegraphics[width=0.46\textwidth]{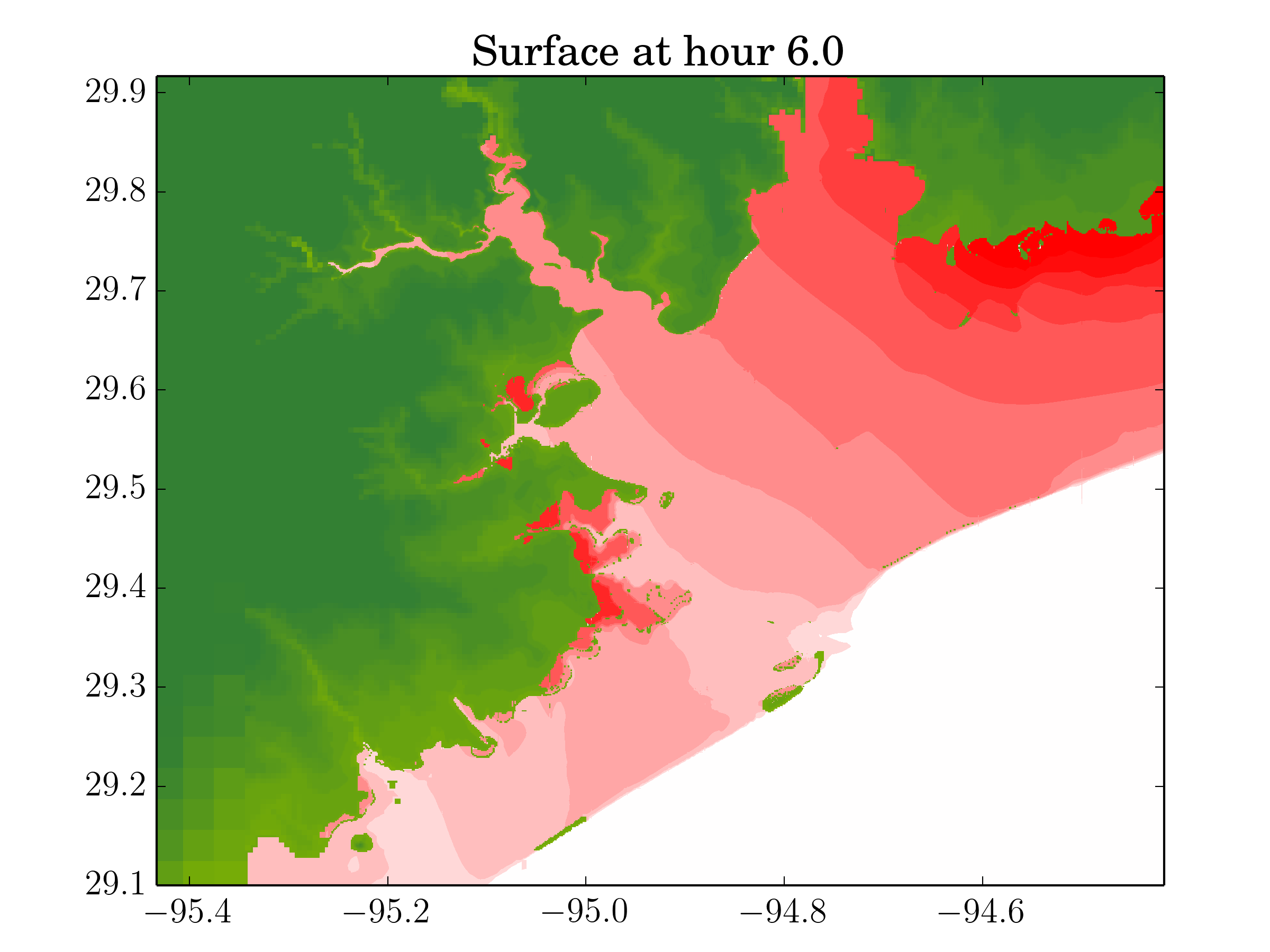} 
    \includegraphics[width=0.46\textwidth]{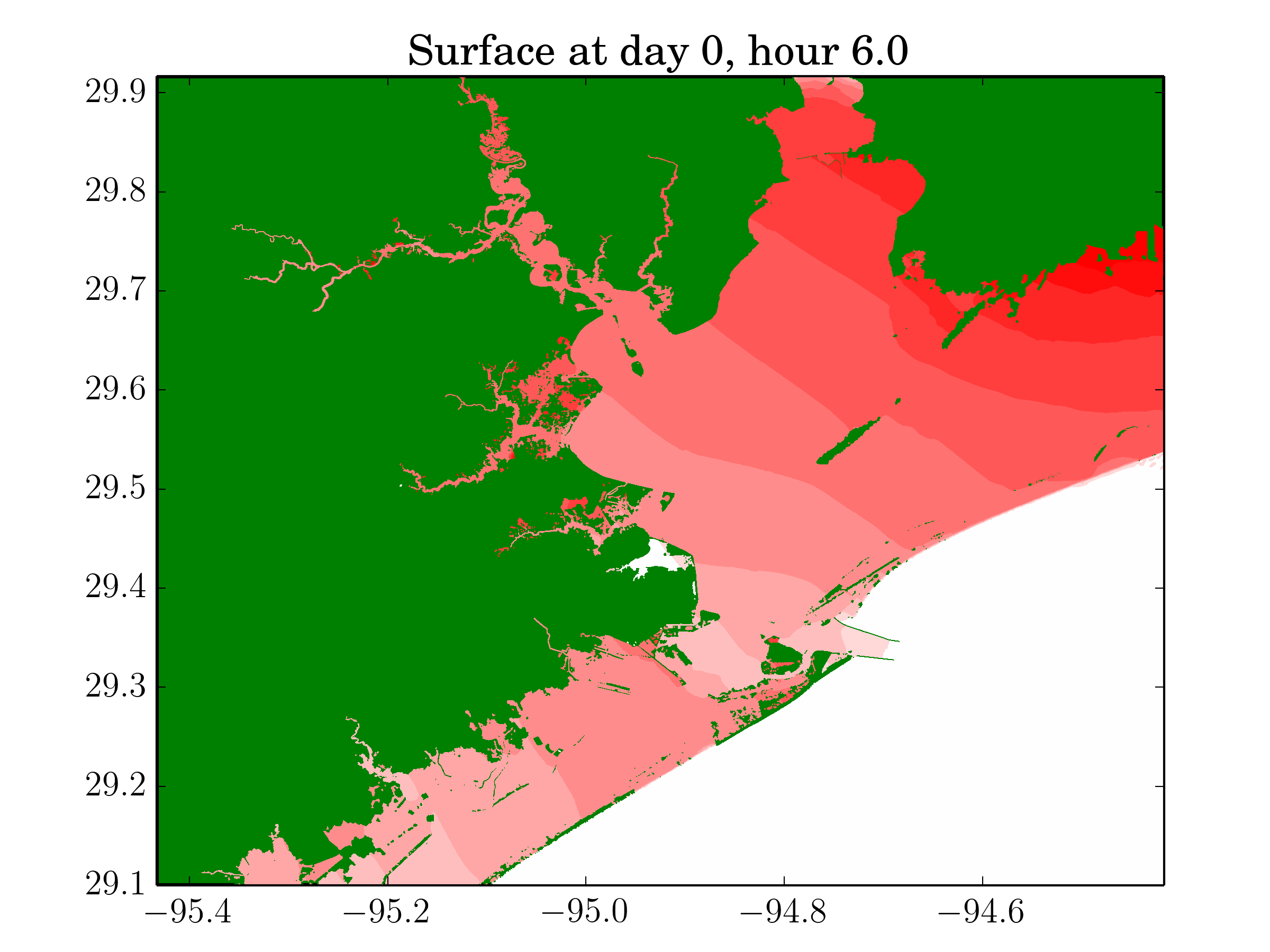} \\
    \includegraphics[width=0.46\textwidth]{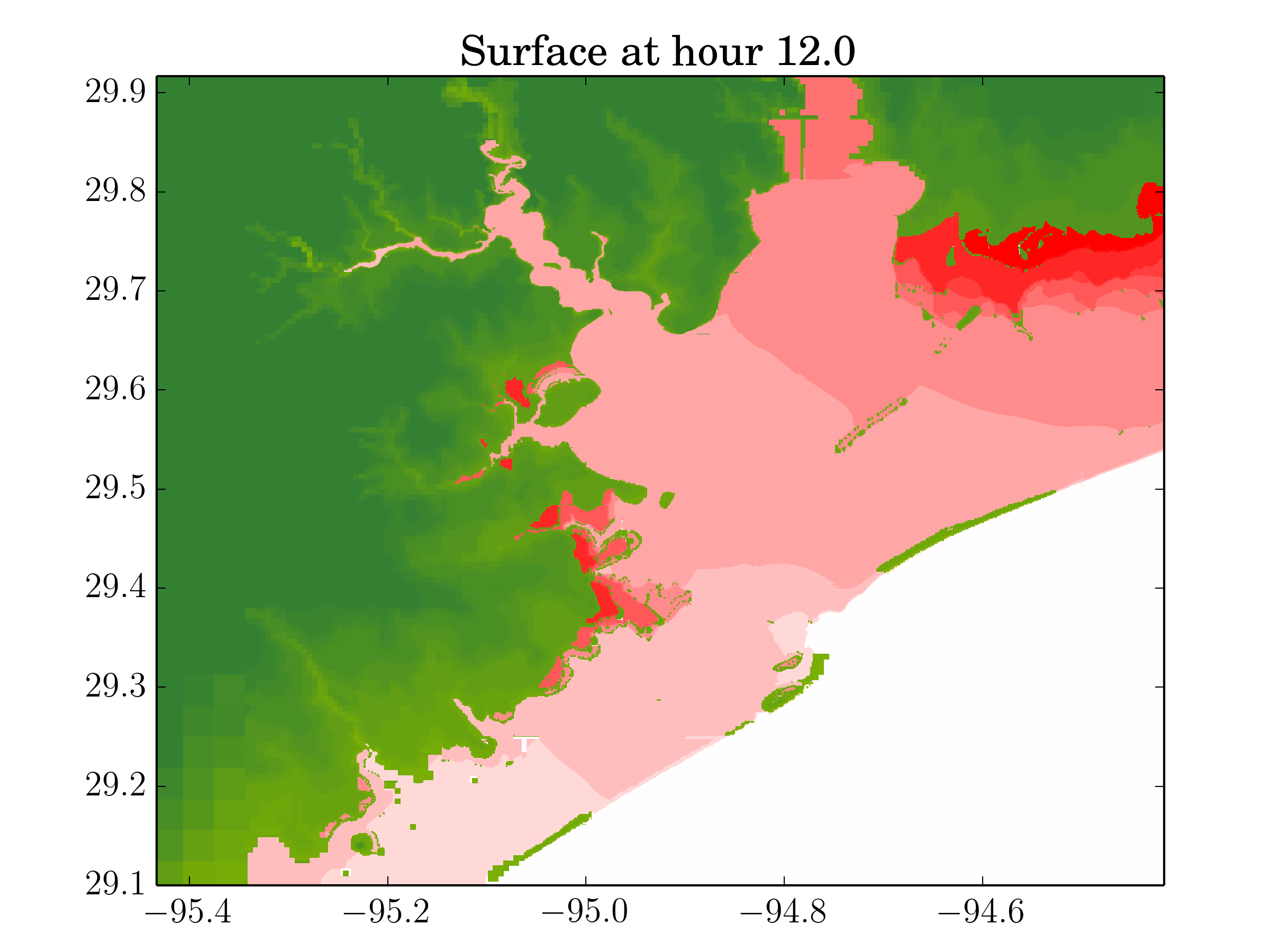} 
    \includegraphics[width=0.46\textwidth]{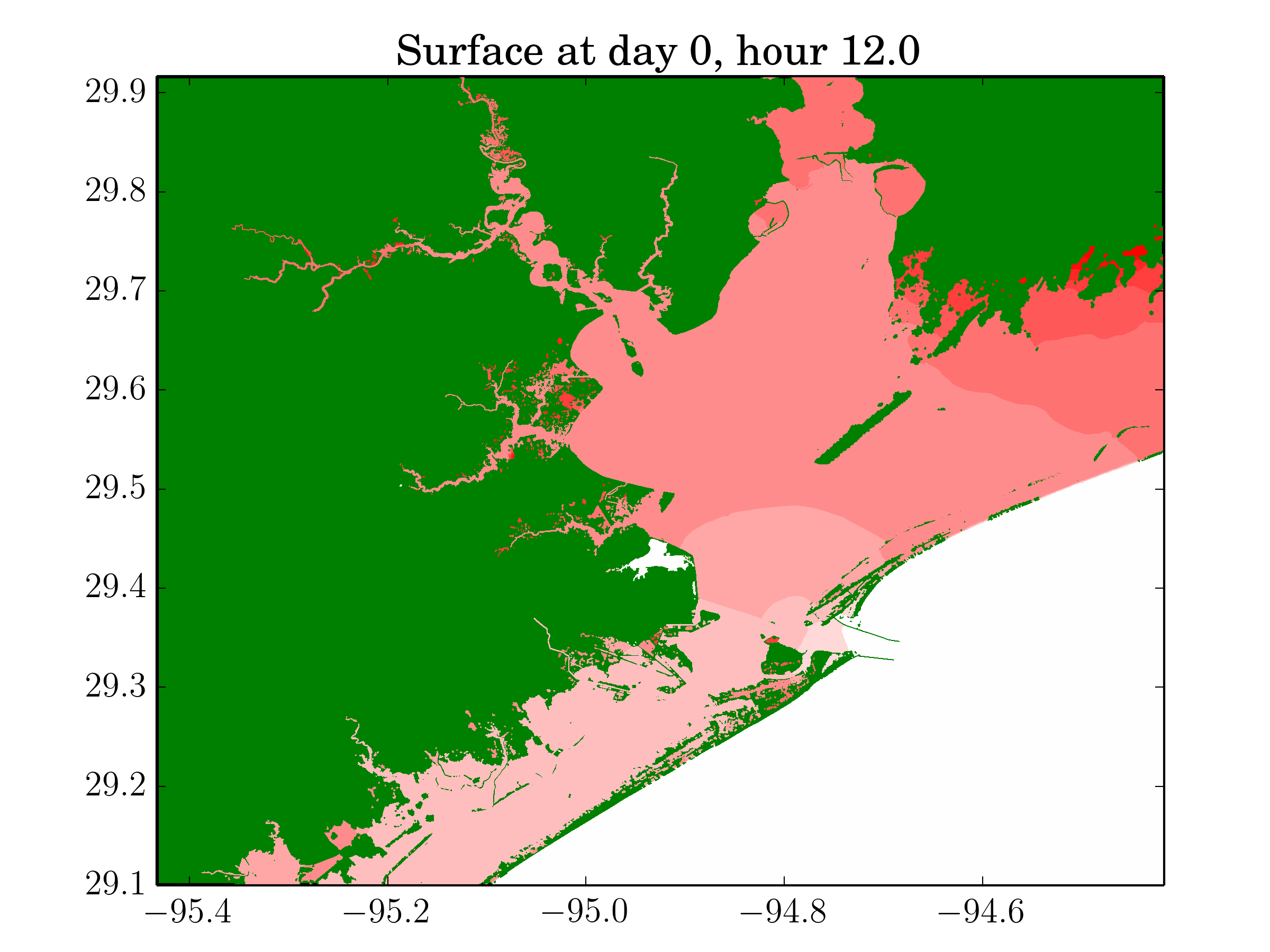} \\
    \includegraphics[width=0.9\textwidth]{surface_cmap.png}
    \caption{Surface elevations in the Galveston Bay region produced by \geoclaw (on the left) and \adcirc (on the right) after Hurricane Ike made landfall.}
    \label{fig:galveston_surface_after}
\end{figure}

\begin{figure}[htb]
    \centering
    \includegraphics[width=0.46\textwidth]{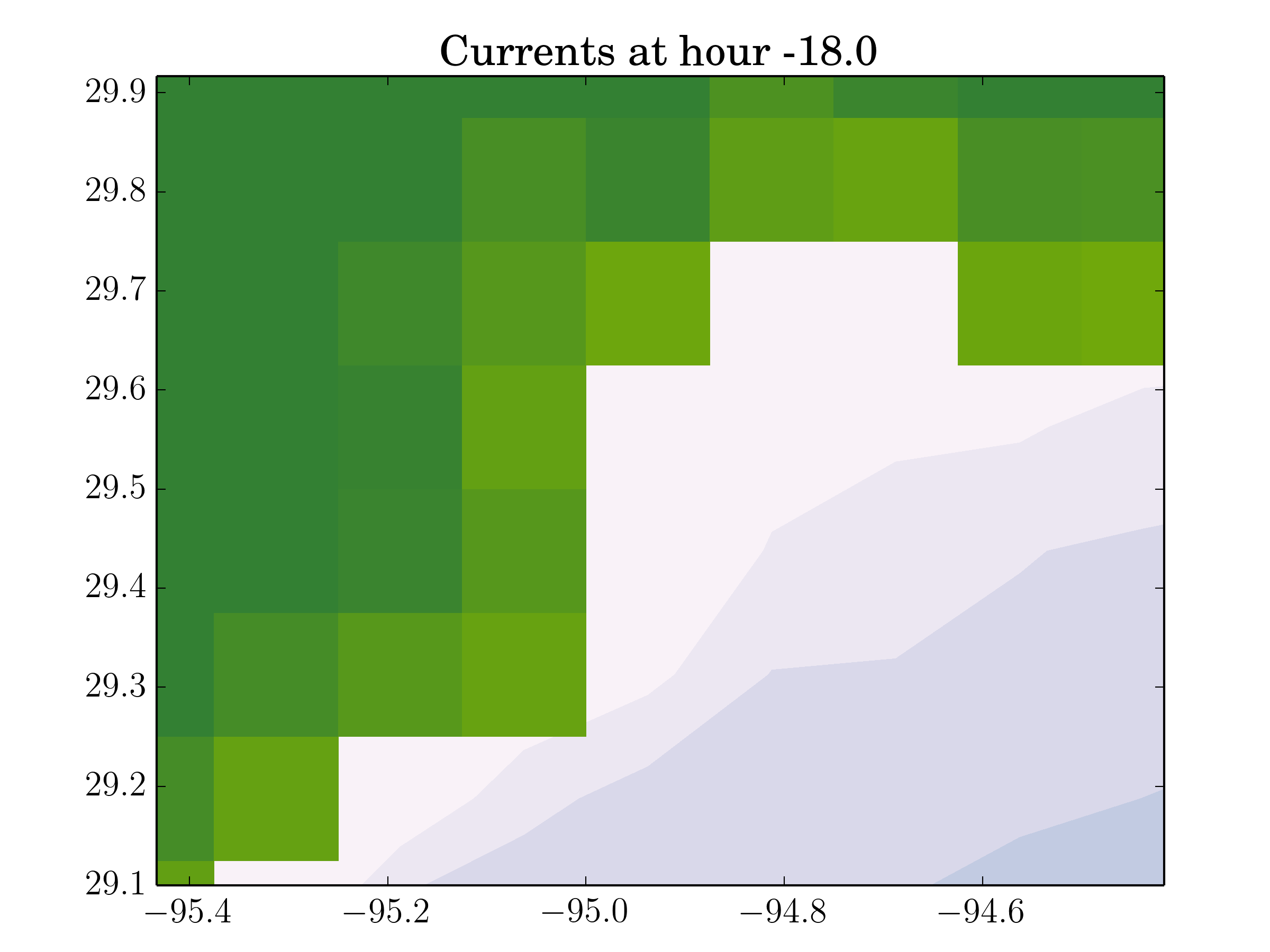}
    \includegraphics[width=0.46\textwidth]{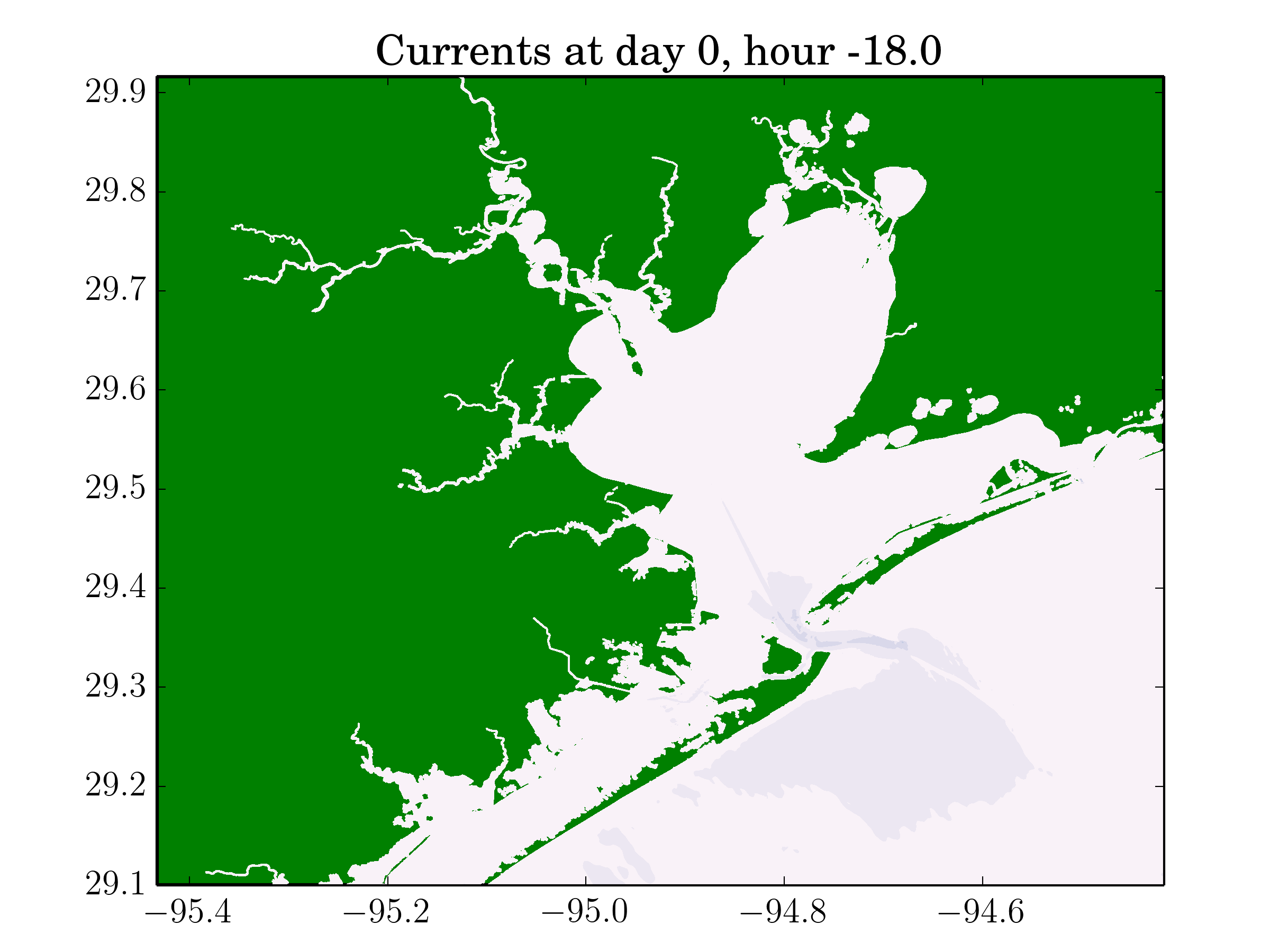} \\
    \includegraphics[width=0.46\textwidth]{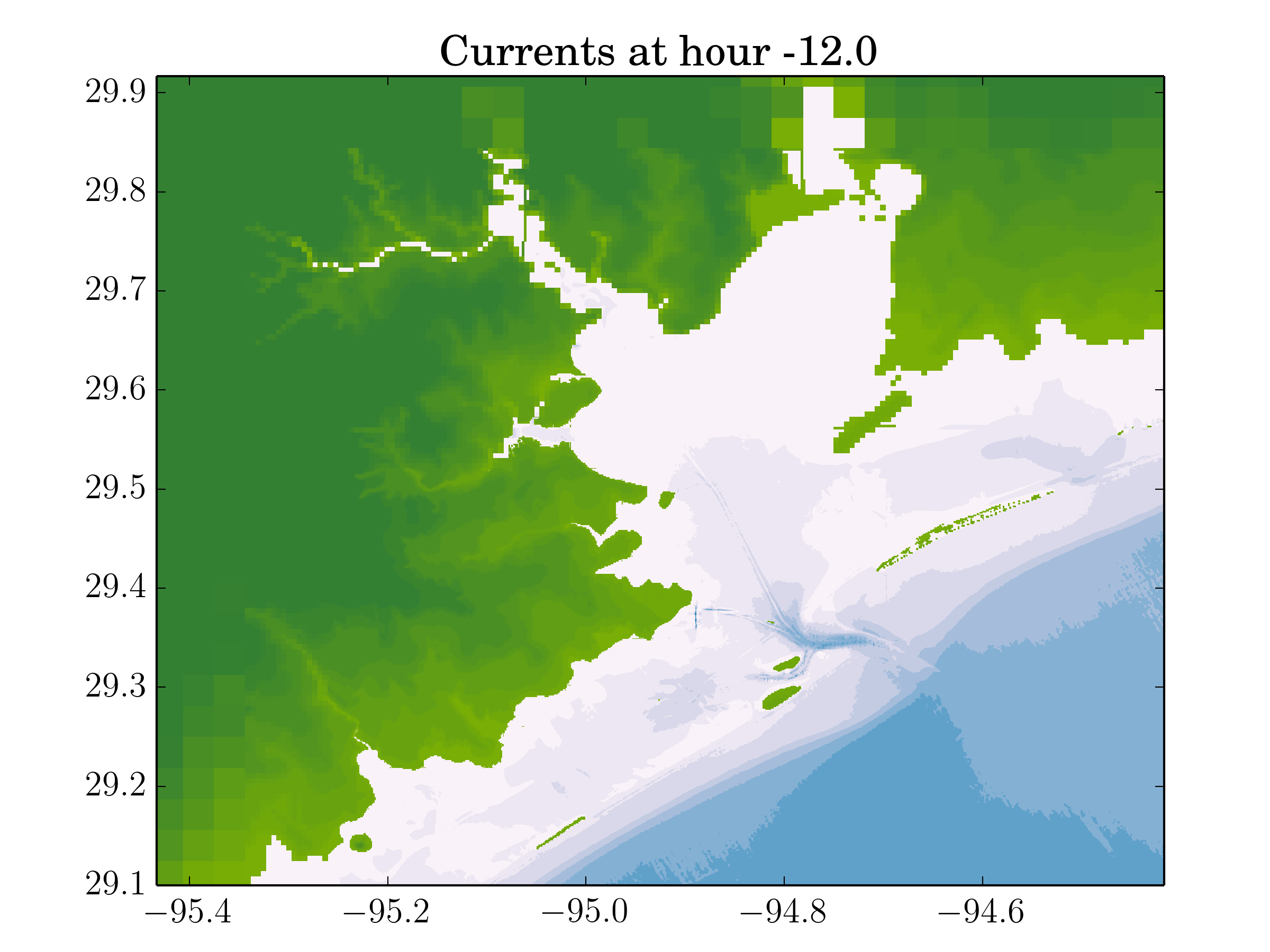}
    \includegraphics[width=0.46\textwidth]{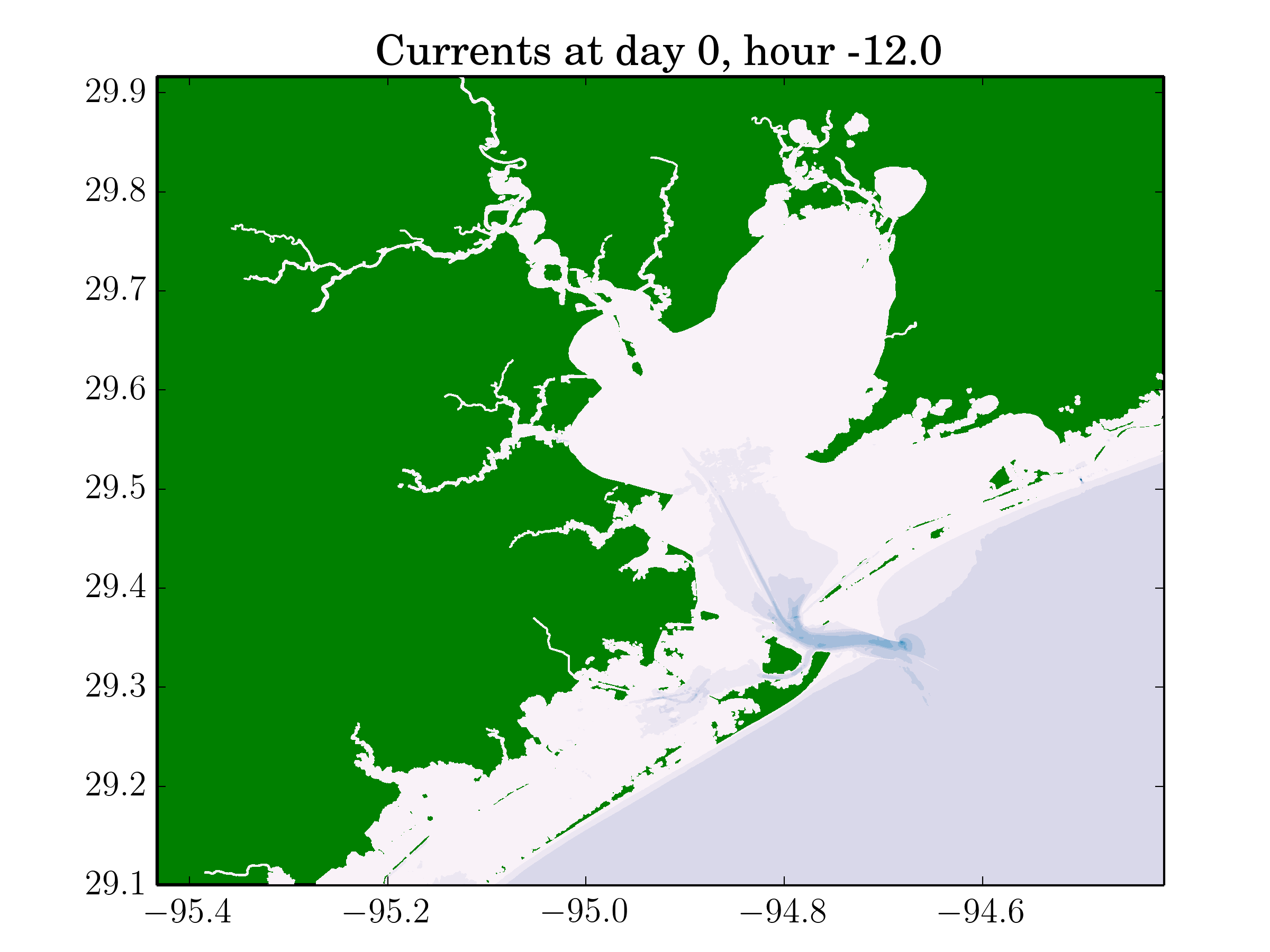} \\
    \includegraphics[width=0.46\textwidth]{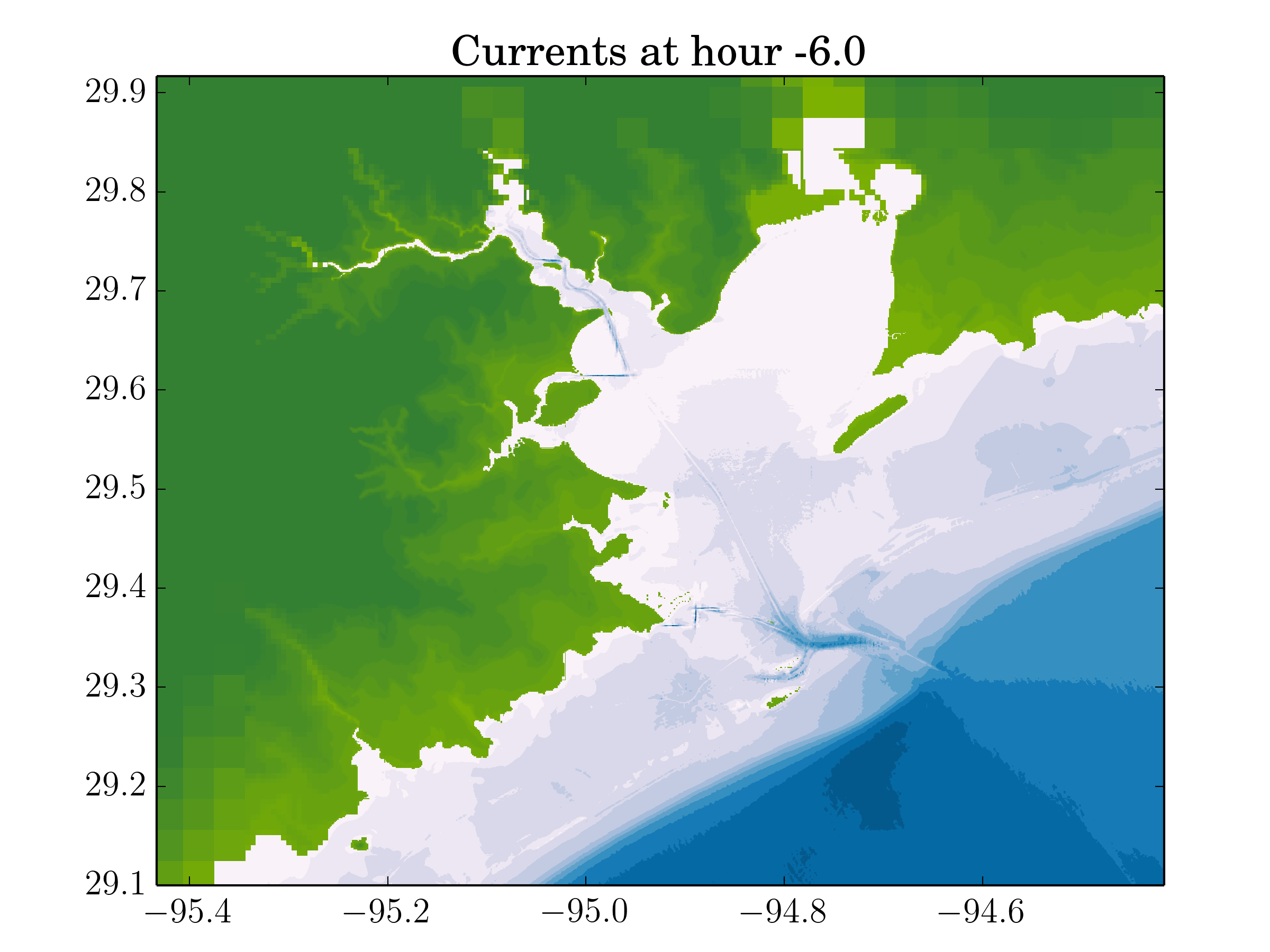}
    \includegraphics[width=0.46\textwidth]{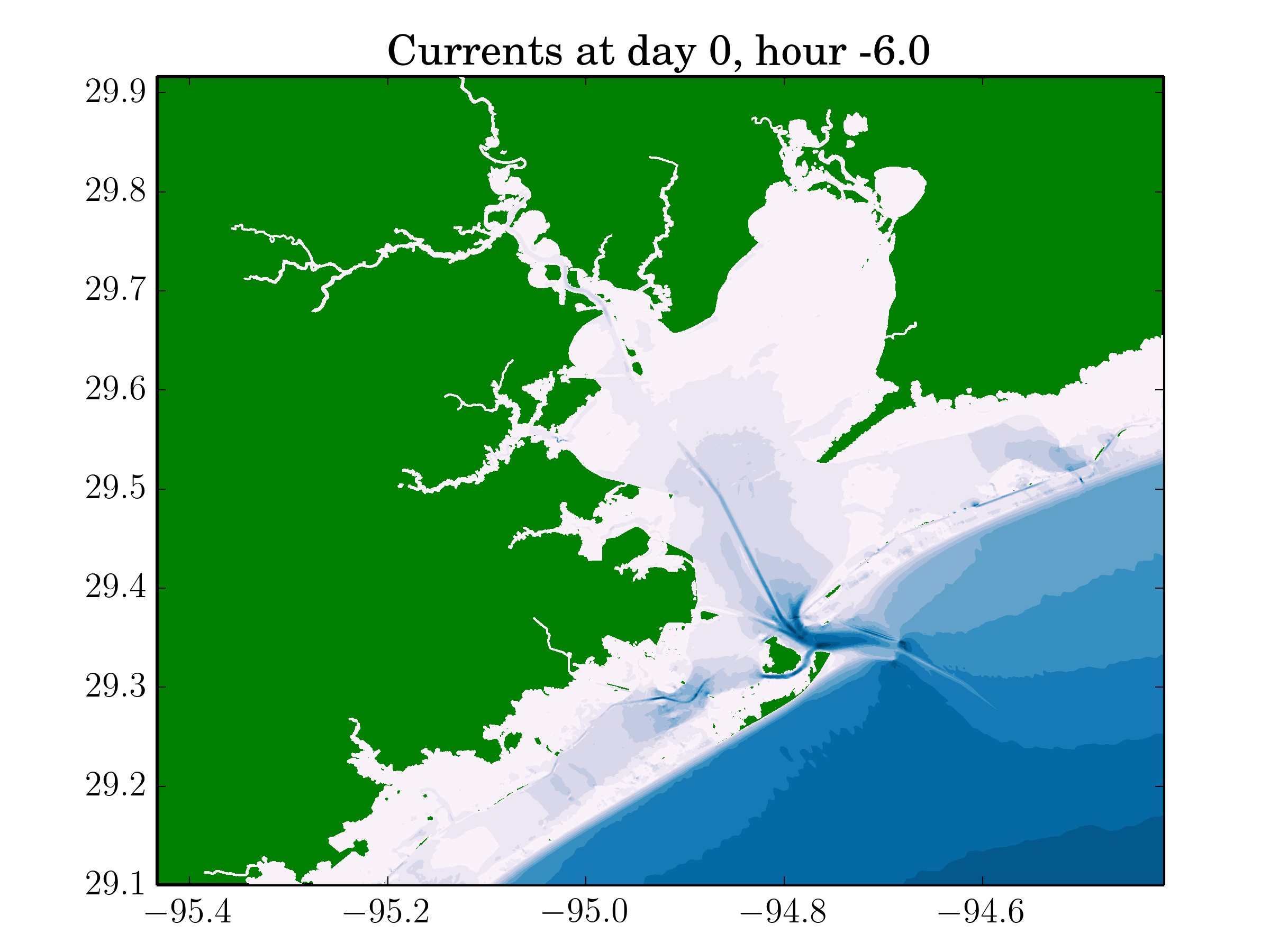} \\
    \includegraphics[width=0.9\textwidth]{speed_cmap.png}
    \caption{Currents in the Galveston Bay region produced by \geoclaw (on the left) and \adcirc (on the right) before Hurricane Ike makes landfall.}
    \label{fig:galveston_speed_before}
\end{figure}

\begin{figure}[htb]
    \centering
    \includegraphics[width=0.46\textwidth]{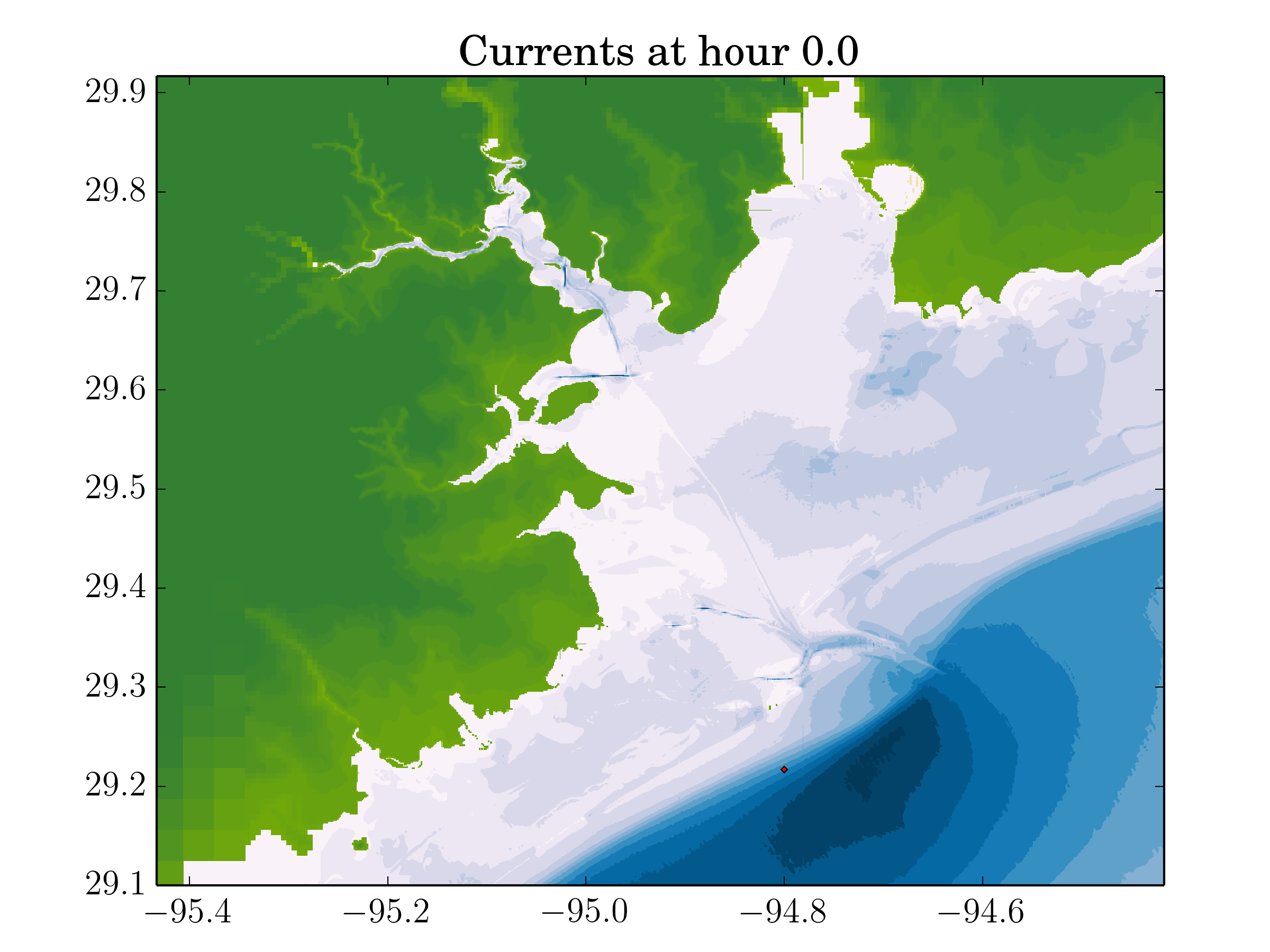}
    \includegraphics[width=0.46\textwidth]{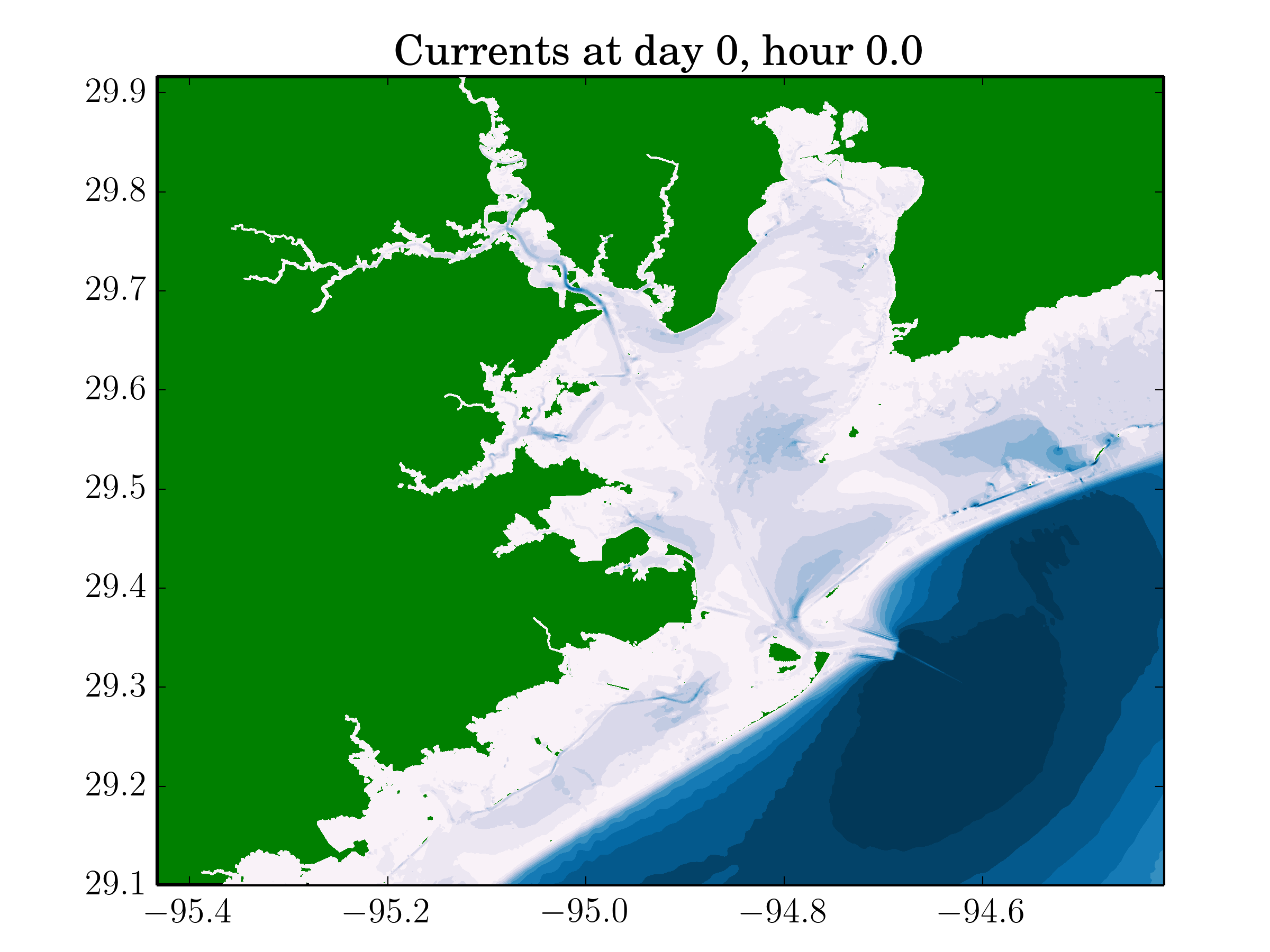} \\
    \includegraphics[width=0.46\textwidth]{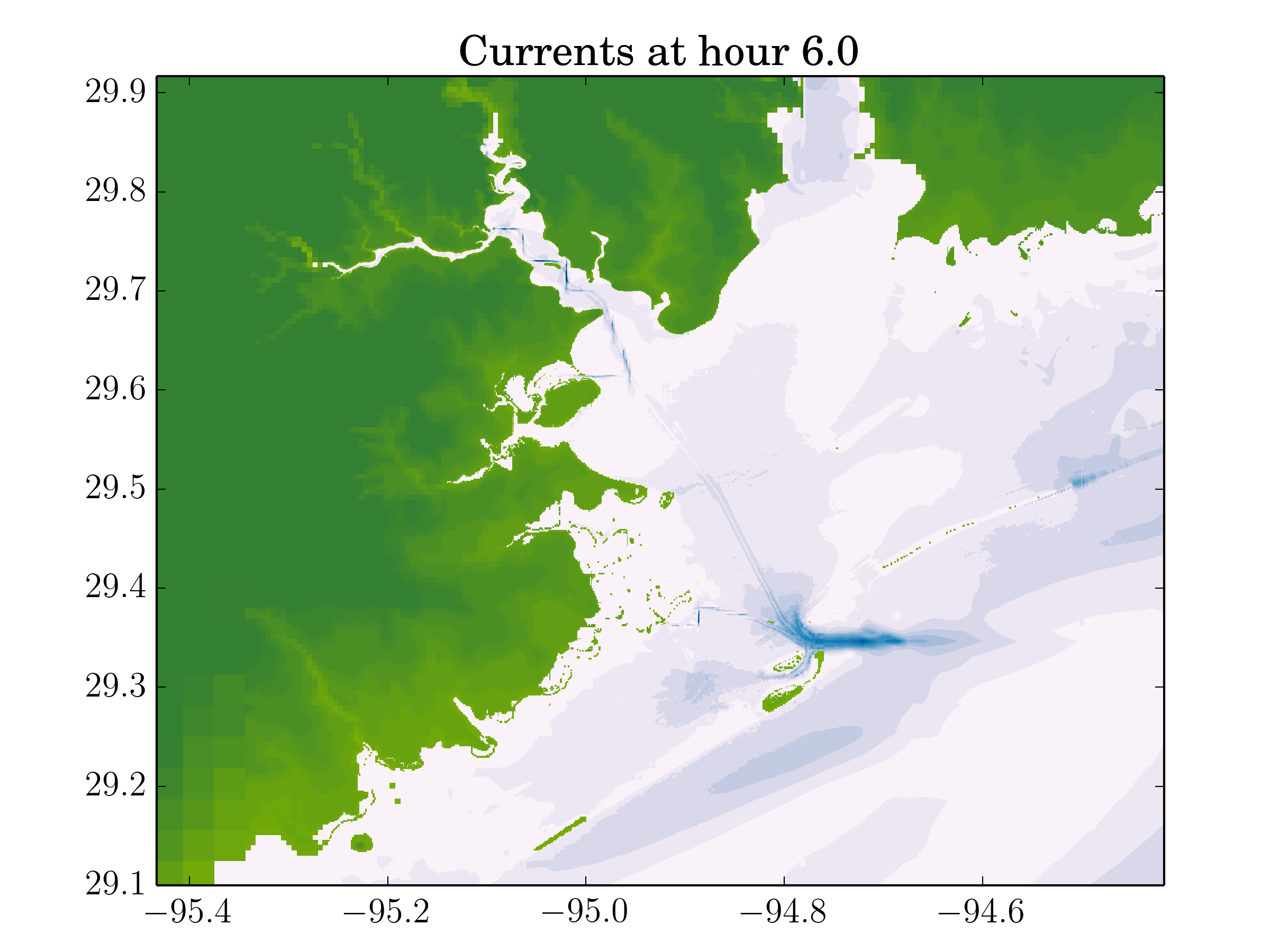}
    \includegraphics[width=0.46\textwidth]{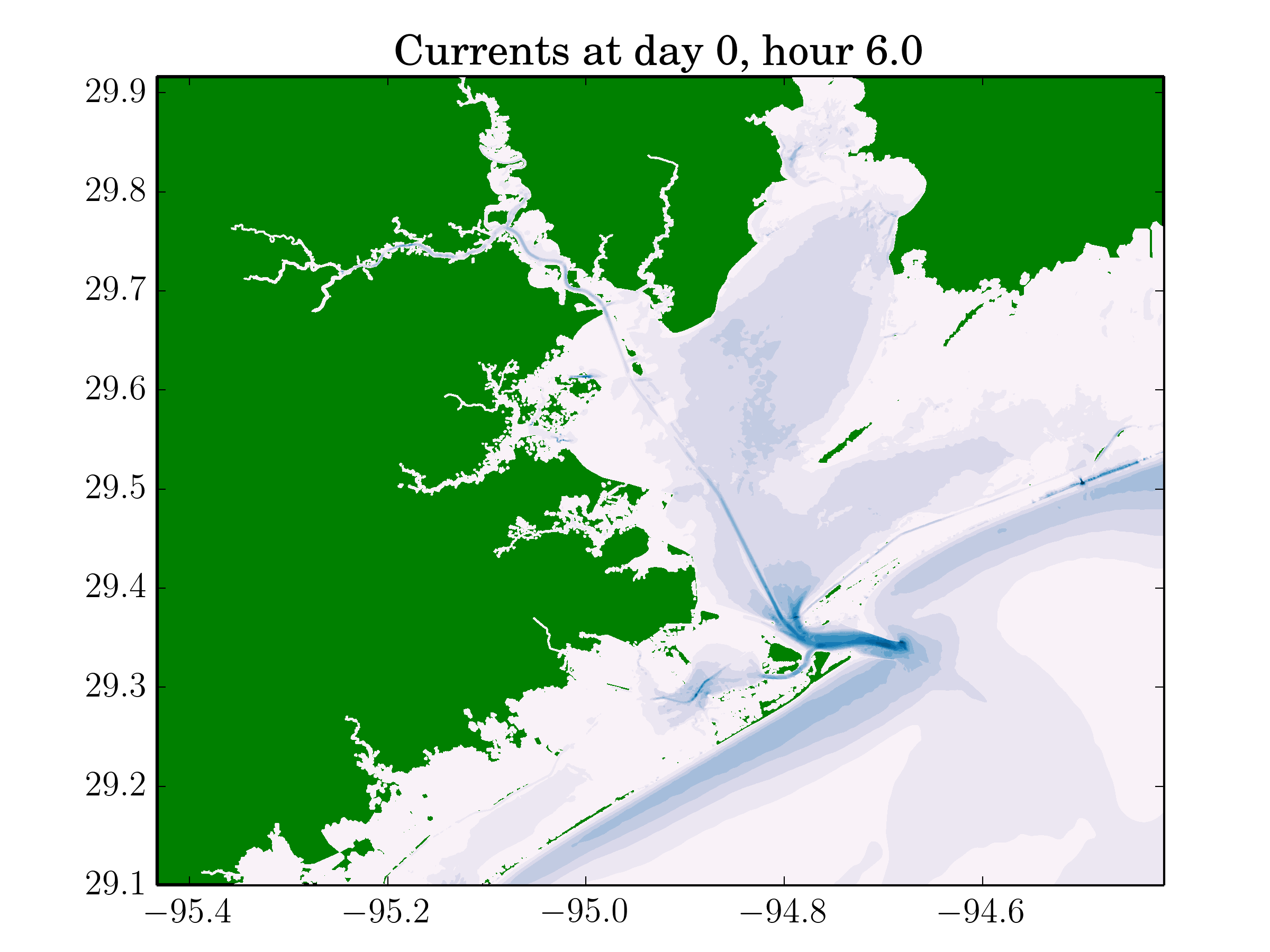} \\
    \includegraphics[width=0.46\textwidth]{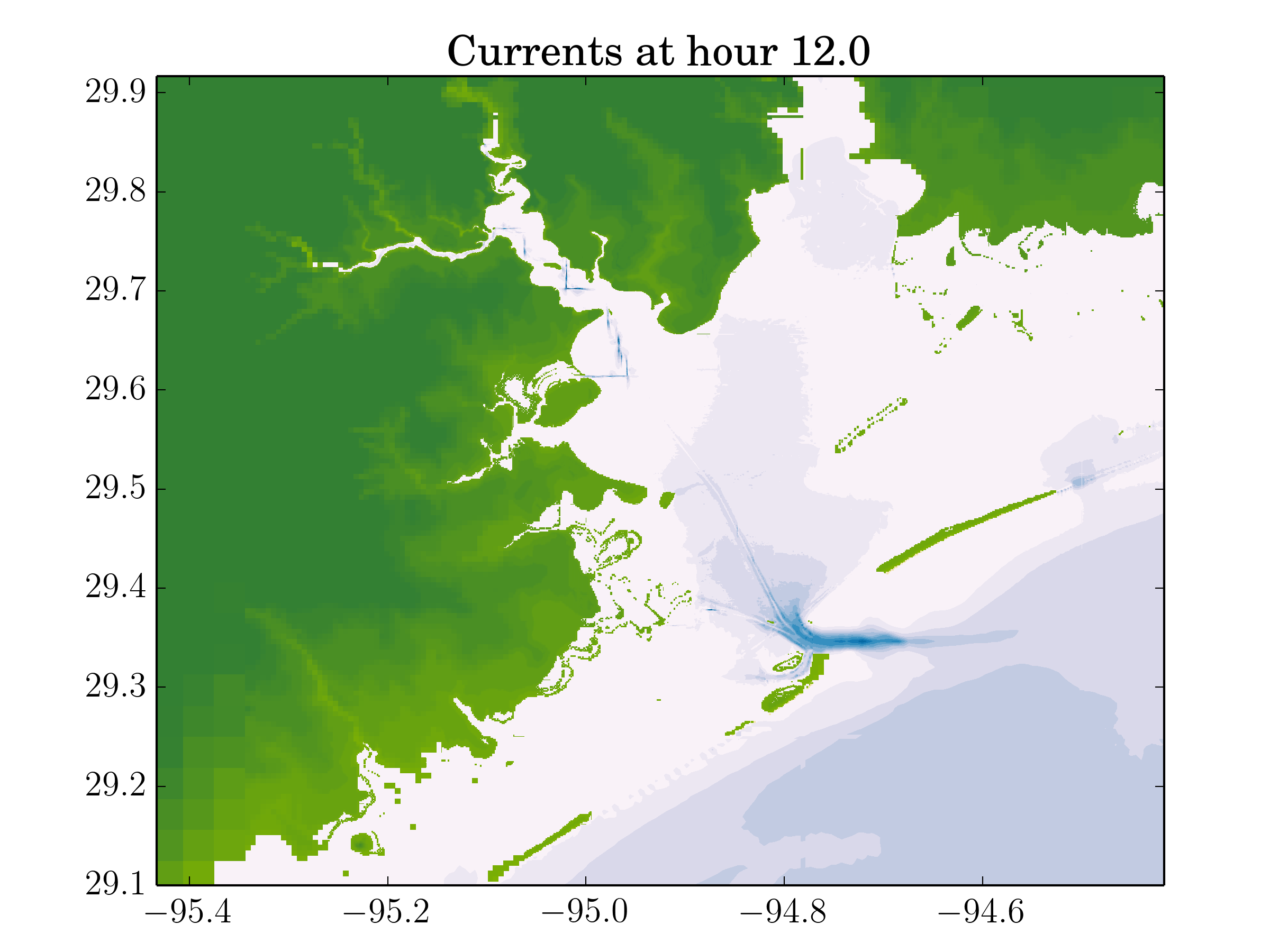}
    \includegraphics[width=0.46\textwidth]{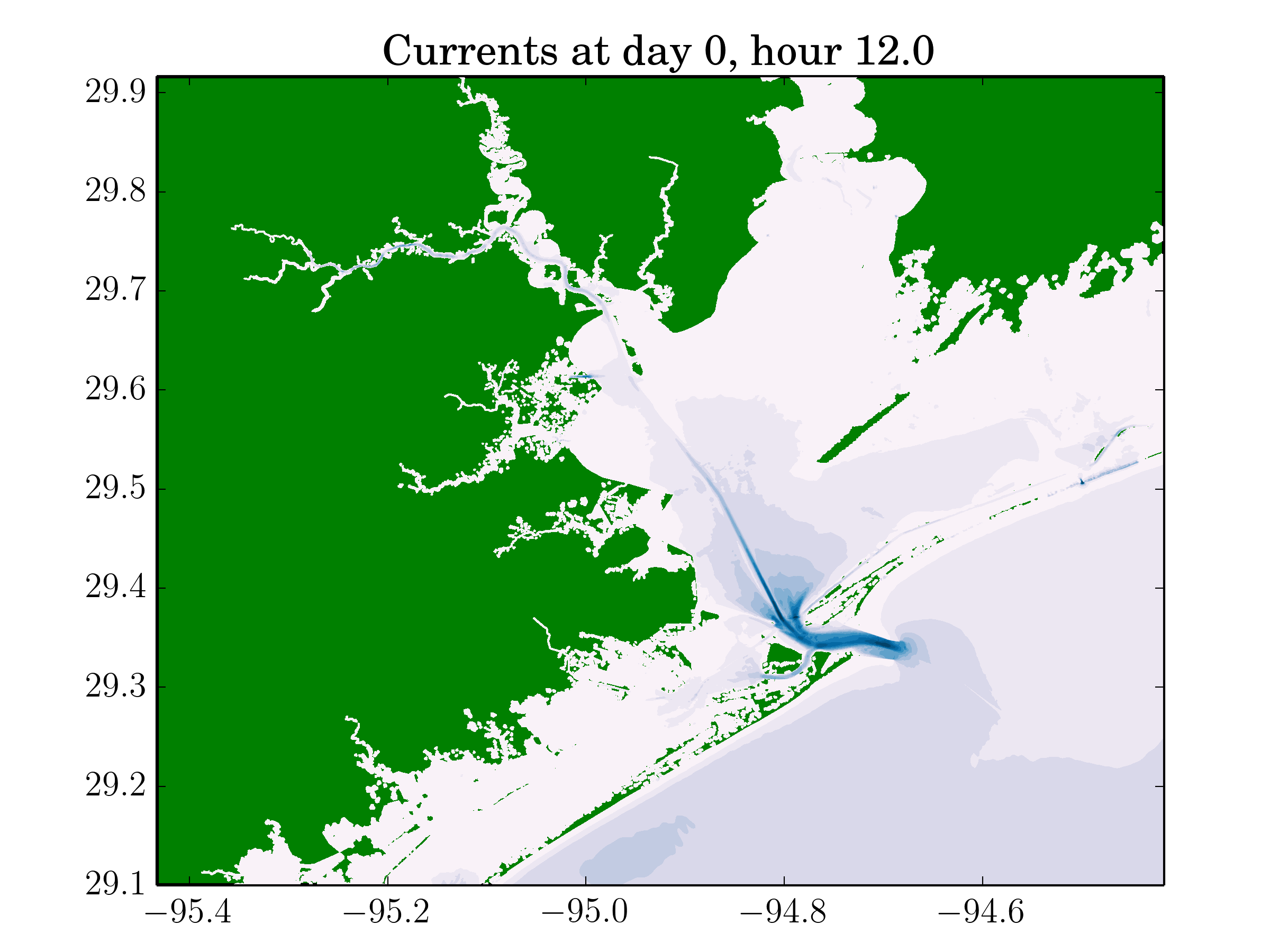} \\
    \includegraphics[width=0.9\textwidth]{speed_cmap.png}
    \caption{Currents in the Galveston Bay region produced by \geoclaw (on the left) and \adcirc (on the right) after Hurricane Ike makes landfall.}
    \label{fig:galveston_speed_after}
\end{figure}

Zooming in on the landfall region, Figures~\ref{fig:galveston_surface_before} and \ref{fig:galveston_surface_after} show surface elevations in the Galveston area before and after landfall respectively.  Here, at 18 hours before landfall, we can see that the refinement criteria for \geoclaw has not refined the region significantly.  This is primarily due to the fact that we are plotting from the datum's sea-level rather than the initial sea-level (0.28 meters) used in both simulations.  The real difference in sea-level is well below the refinement threshold of 1 meter of surge and therefore refinement has not been triggered.  Looking at Figures~\ref{fig:galveston_speed_before} and \ref{fig:galveston_speed_after} that show currents in the Galveston area before and after landfall, we can see that the velocity of the water in the area at 18 hours before landfall is near zero and again refinement is not triggered.  Beyond the first snap-shot, we can see that \geoclaw and \adcirc produce currents that are structurally very similar with \geoclaw producing slightly more diffusive results.  Again it should be noted, especially when comparing the current plots, that the \adcirc simulation has a much more complex friction field than the \geoclaw simulation.

\subsubsection{Gauge Data Comparisons} \label{ssub:gauge_data}

\begin{figure}[tb]
    \includegraphics[width=1.0\textwidth]{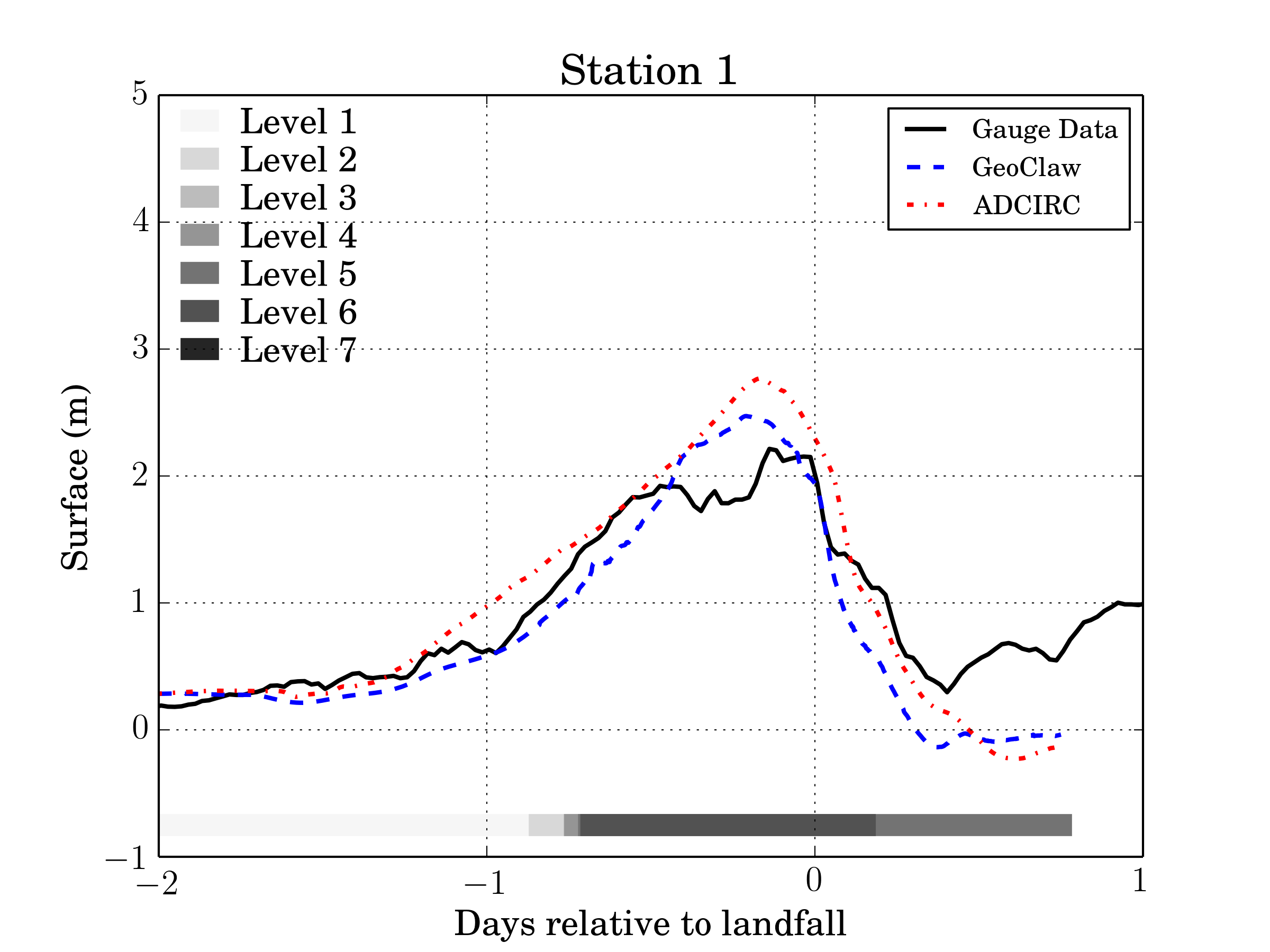} 
    \caption{Comparison of numerical gauge data from \geoclaw and \adcirc and collected data during Hurricane Ike at location 1 from Table~\ref{tab:gauge_locations}.  The shading below the gauge data represents the refinement level used to record the gauge.}
    \label{fig:gauge_1}
\end{figure}

\begin{figure}[tb]
    \includegraphics[width=1.0\textwidth]{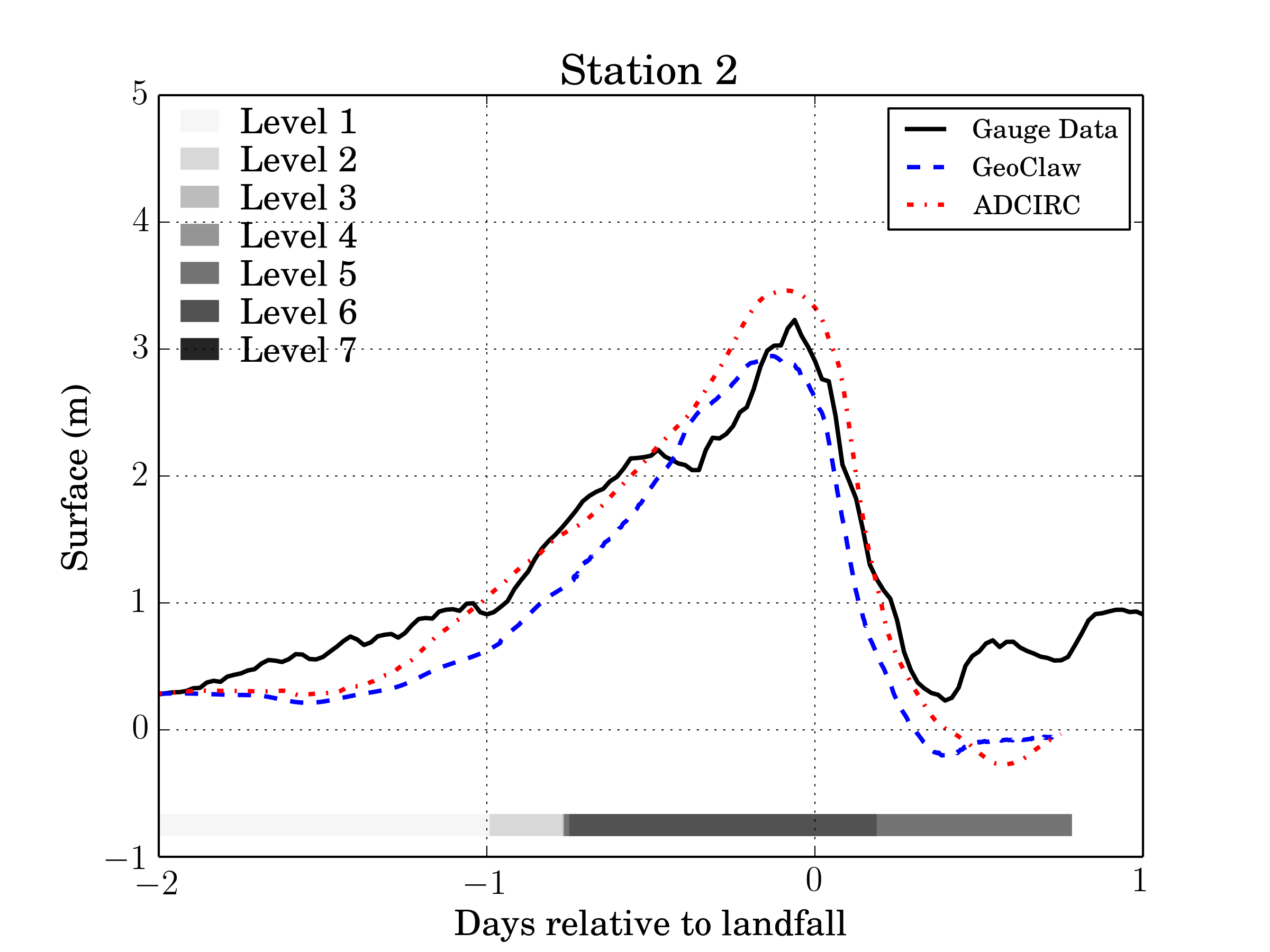}
    \caption{Comparison of numerical gauge data from \geoclaw and \adcirc and collected data during Hurricane Ike at location 2 from Table~\ref{tab:gauge_locations}.  The shading below the gauge data represents the refinement level used to record the gauge.}
    \label{fig:gauge_2}
\end{figure}

\begin{figure}[tb]
    \includegraphics[width=1.0\textwidth]{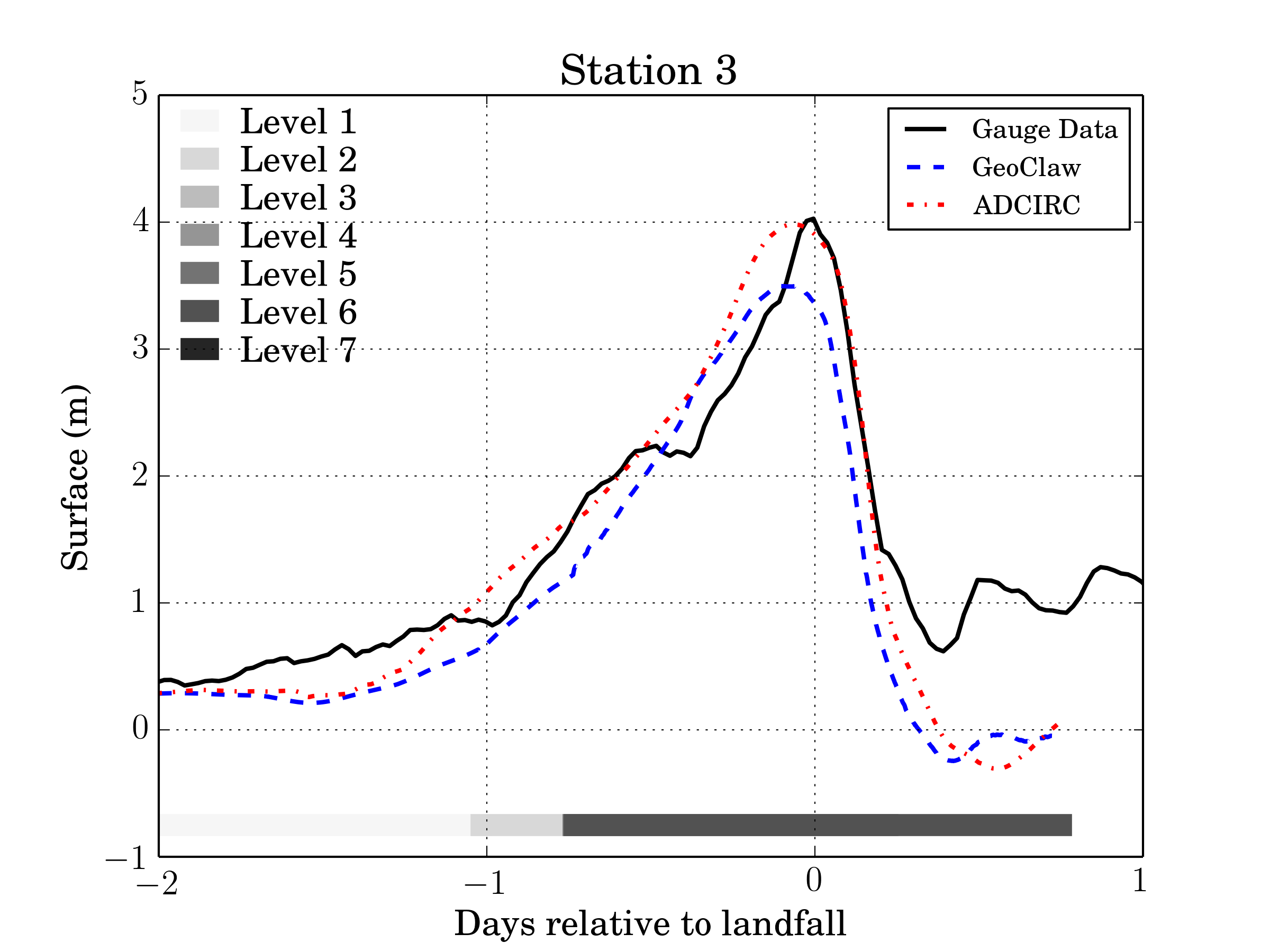} 
    \caption{Comparison of numerical gauge data from \geoclaw and \adcirc and collected data during Hurricane Ike at location 3 from Table~\ref{tab:gauge_locations}.  The shading below the gauge data represents the refinement level used to record the gauge.}
    \label{fig:gauge_3}
\end{figure}

\begin{figure}[tb]
    \includegraphics[width=1.0\textwidth]{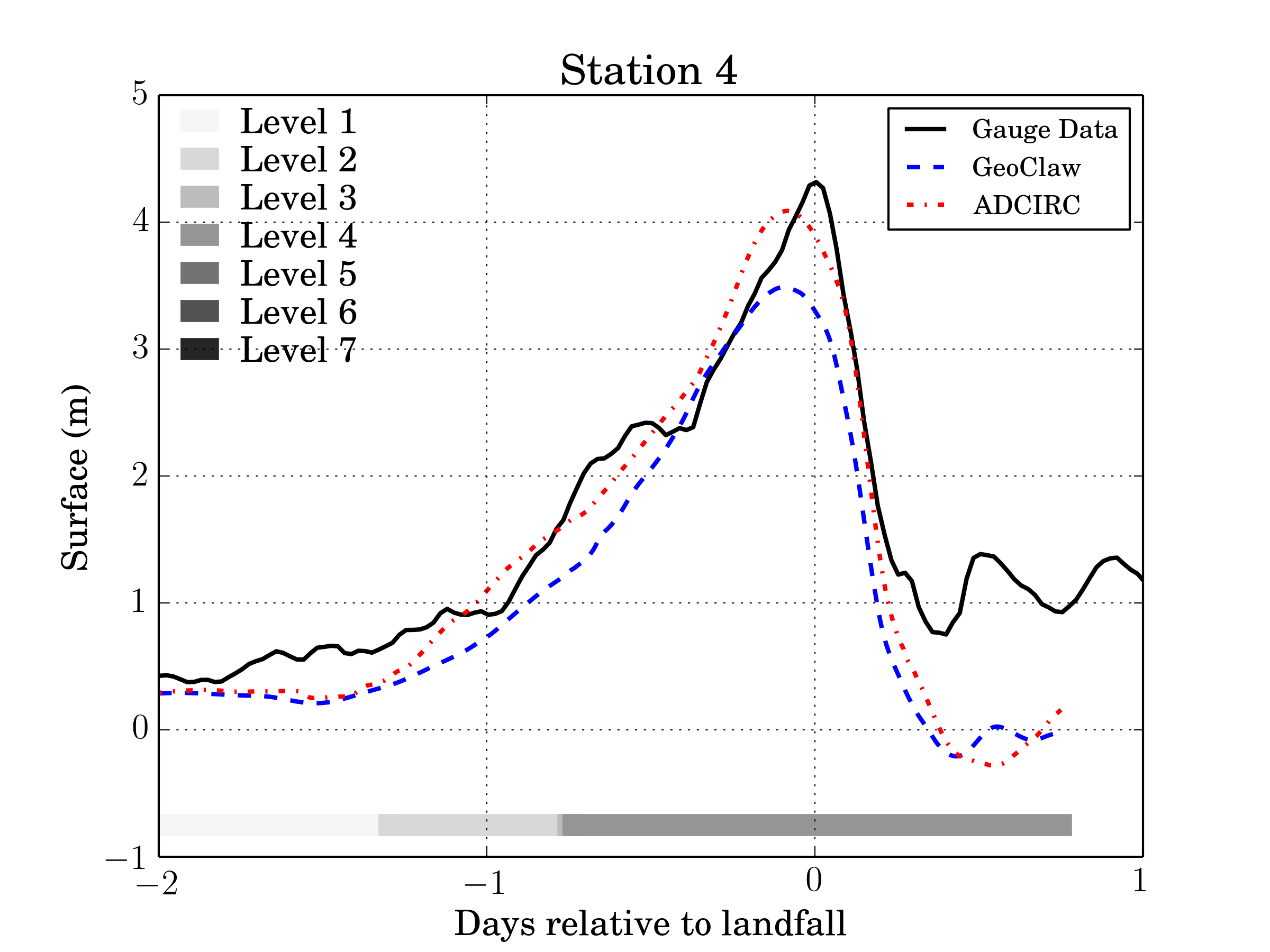}
    \caption{Comparison of numerical gauge data from \geoclaw and \adcirc and collected data during Hurricane Ike at location 4 from Table~\ref{tab:gauge_locations}.  The shading below the gauge data represents the refinement level used to record the gauge.}
    \label{fig:gauge_4}
\end{figure}

Gauge data was collected during Hurricane Ike off the Texas coast near Galveston \cite{Kennedy:2011kt}.  As another point of comparison, this gauge data was compared to both the \geoclaw and \adcirc numerical gauges.  In \geoclaw this data is obtained by nearest neighbor spatial interpolation to the location of the gauge and output in time every time step on the finest grid containing the gauge location.  In \adcirc spatial interpolation is also used and is output at a user defined frequency no smaller than the constant time step taken everywhere.  Table~\ref{tab:gauge_locations} describes the location of each gauge and Figures~\ref{fig:gauge_1}, \ref{fig:gauge_2}, \ref{fig:gauge_3}, and \ref{fig:gauge_4} provide a comparison between the data collected, and the \geoclaw and \adcirc gauge data.  

The first observation that can be made is that the numerical gauges in the \geoclaw simulation are consistently lower than the \adcirc gauges.  This could be due to the resolution near the gauges as the \geoclaw simulation never refines to the highest level of resolution in the gauge locations (as shown by the shading below each gauge figure) where as the \adcirc simulation's resolution at these gauge locations is approximately .  Apart from this, both simulations agree on the timing of the peak surge.  Neither the \geoclaw nor the \adcirc simulations capture the initial peak of surge (the forerunner).  This may be due to the use of the Holland based approximated storm fields or model inaccuracies.

\begin{table}[htb]
    \begin{center}
    \begin{tabular}{c|ccc}
    \hline
    \textbf{Gauge} & \textbf{Latitude} & \textbf{Longitude} & \textbf{Figure} \\
    \hline
    1  & 29.07 N & 95.04 W & \ref{fig:gauge_1}\\ 
    2  & 29.28 N & 94.71 W & \ref{fig:gauge_2}\\
    3  & 29.49 N & 94.39 W & \ref{fig:gauge_3} \\
    4  & 29.58 N & 94.13 W & \ref{fig:gauge_4}\\
    % 5  & 29.70 N & 95.00 W & Houston Ship Channel at SH 146 bridge \\
    % 6  & 29.74 N & 95.14 W & Houston Ship Channel at Port of Houston \\
    % 7  & 29.55 N & 95.08 W & Galveston Bay: Clear Creek entrance \\
    % 8  & 29.76 N & 94.75 W & Galveston Bay: Upper Trinity Bay \\
    % 9  & 29.72 N & 95.27 W & Manchester: Houston Ship Channel \\
    % 10 & 29.52 N & 94.51 W & Bolivar Peninsula, Rollover Pass \\
    % R  & 27.67 N & 97.12 W & ~ \\
    % S  & 28.21 N & 96.55 W & ~ \\
    % U  & 28.63 N & 95.75 W & ~ \\
    % V  & 28.87 N & 95.32 W & ~ \\
    \hline
    \end{tabular}
    \end{center}
    \caption{Locations of gauges presented.  The background field is of the sea-surface at landfall.
    \label{tab:gauge_locations}}
\end{table}

\begin{figure}[htb]
    \centering
    \includegraphics[width=1.0\textwidth]{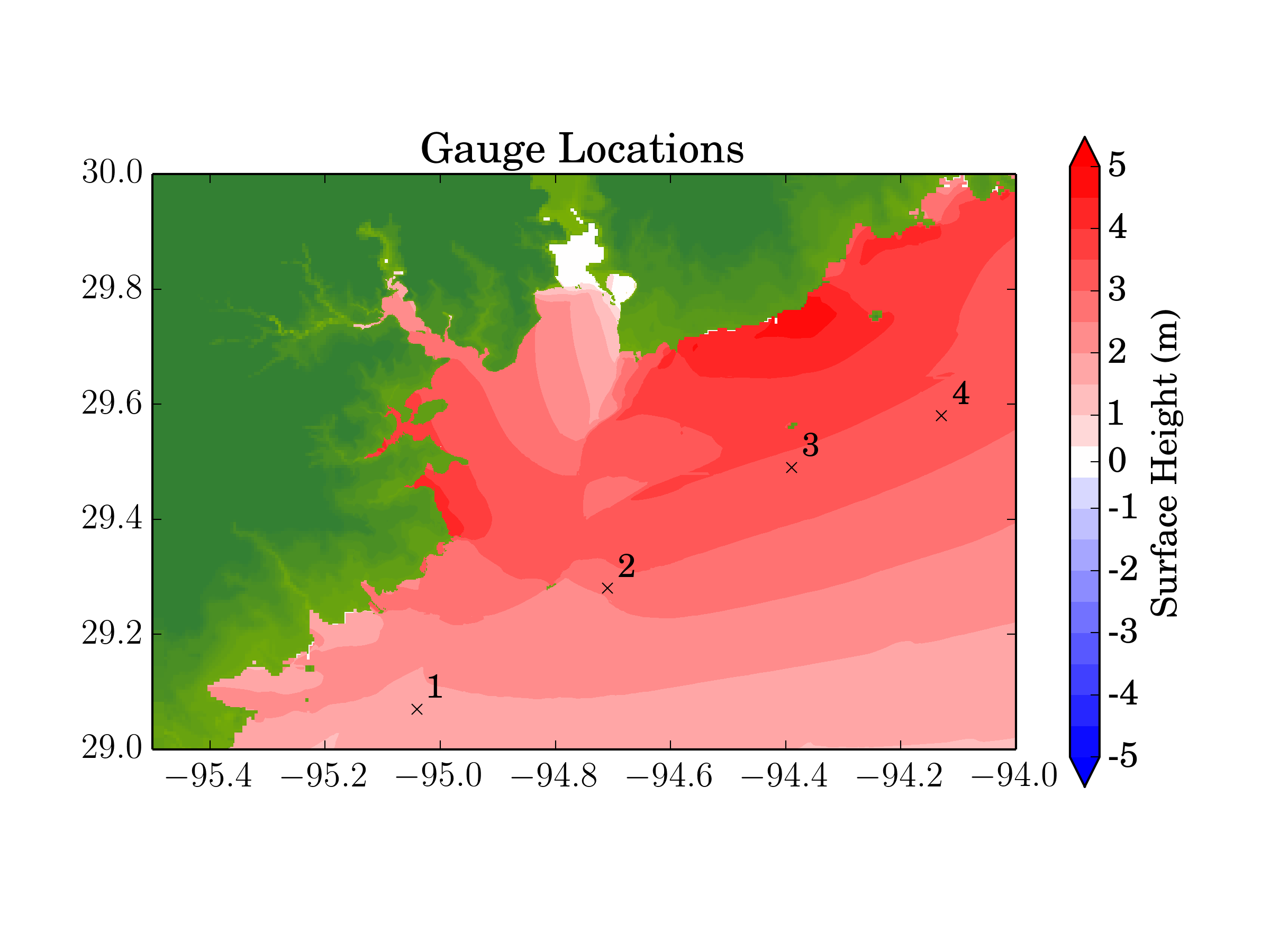}
    \caption{Locations of gauges presented.  The color is indicative of the surge heights at landfall (see Figure~\ref{fig:galveston_surface_after}) to give a perspective on where the gauges are relative to the surge.}
    \label{fig:figure1}
\end{figure}


\subsubsection{Computational Cost} \label{ssub:computation}

\begin{figure}[htb]
    \centering
    \begin{subfigure}[b]{0.49\textwidth}
        \includegraphics[width=\textwidth]{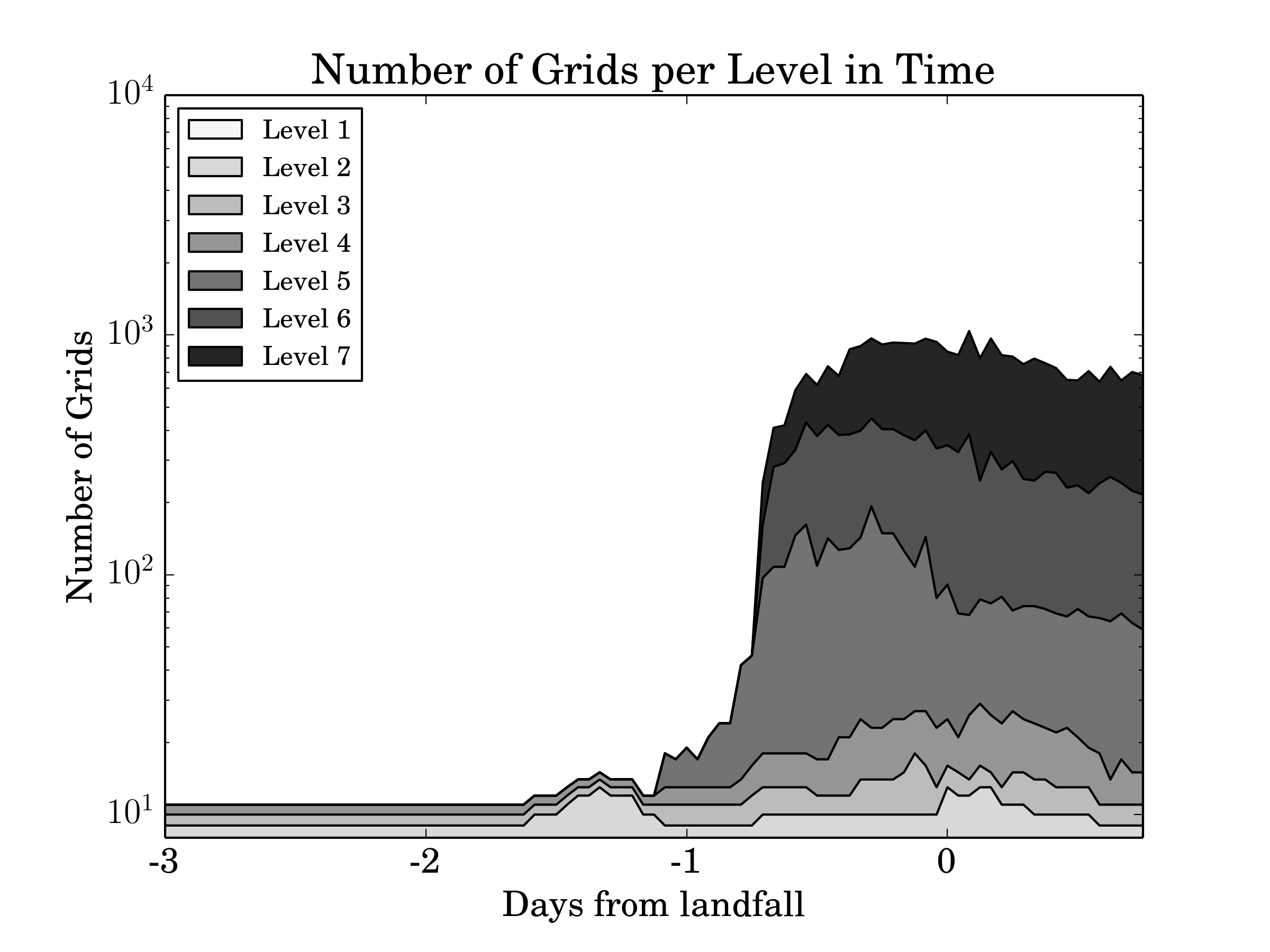}
        \caption{}
        \label{fig:num_grids}
    \end{subfigure}
    \begin{subfigure}[b]{0.49\textwidth}
        \includegraphics[width=\textwidth]{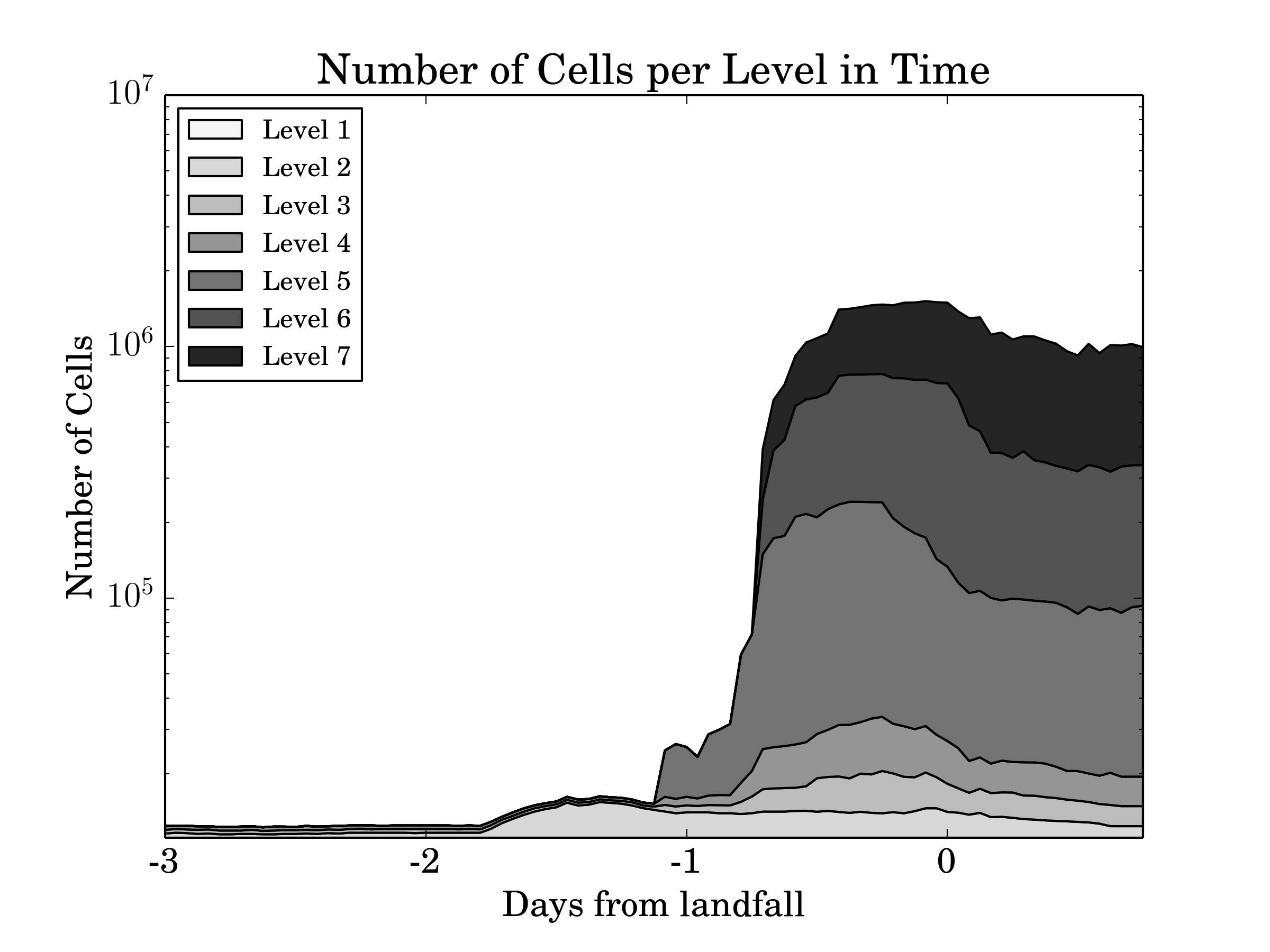}
        \caption{}
        \label{fig:num_cells}
    \end{subfigure}
    \caption{Domain resolution statistics for the \geoclaw simulation.  In Figure~\ref{fig:num_grids} a stack plot is shown representing the total number of grids per level and total.  Similarly in Figure~\ref{fig:num_cells} the number of grids cells is shown.}
    \label{fig:num_grids_cells}
\end{figure}

As one of the objectives of introducing AMR to storm surge forecasting was to reduce the computational cost, we report the overall wall clock time and computational cost of each of the compared simulations.  The \adcirc simulation was run on 4000 cores of the Stampede computer cluster at the Texas Advanced Computing Center.  Stampede is comprised of nodes containing two 8 core Xeon ES-2680 processors with 32GB of memory.  These nodes also contain Intel Xeon Phi SE10P co-processors but these were not utilized in the simulation.  The \geoclaw simulation was run on a MacBook Pro with an Intel Core i7 chip containing 4 cores and 16GB of memory.  As mentioned earlier, \adcirc utilizes MPI for parallelism while \geoclaw uses OpenMP.  Table~\ref{tab:timings} contains the raw performance and overall computational costs of the \adcirc and \geoclaw simulations.  As another measure of the computational cost in time of the \geoclaw simulation, Figure~\ref{fig:num_grids_cells} records the number of grids and grid cells used as a function of time.  Although \geoclaw has an overall wall clock time approximately 3.4 times that of the \adcirc simulation, the computational cost is significantly lessened.  It is difficult to fairly compare the equivalent number of degrees of freedom between the \geoclaw and \adcirc simulations but taking the maximum number of grid cells \geoclaw ever uses as a base, Figure~\ref{fig:num_grids_cells} shows the probable source of this decrease in computational cost, the number of grid cells used drastically changes over time.

\begin{table}[htb]
    \begin{center}
    \begin{tabular}{l|ccc}
    \hline
    \textbf{Package} & \textbf{Cores Used} & \textbf{Wall Clock Time} & \textbf{Core Hours} \\
    \hline
        \adcirc & 4000 & 35 minutes & 2333 hours \\
        \geoclaw & 4 & 2 hours & 8 hours \\
    \hline
    \end{tabular}
    \end{center}
    \caption{Timing and computational cost comparisons between \adcirc and \geoclaw. \label{tab:timings}}
\end{table}

\subsection{Discussion}

We close this section with some general observations and conclusions regarding the comparison presented.  First of all, the \geoclaw simulation compares favorably with the \adcirc simulation.  The peaks surge characteristics are similar although \geoclaw appears more diffusive in this comparison.  This may partly be due to the differences between the two simulations in the friction coefficient field but it is more likely to be a combination of insufficient resolution and the nature of the structured grid \geoclaw uses.  The \adcirc grid is designed to mesh waterways and coastlines well where as \geoclaw has to rely on sufficient resolution to represent a channel (unless said waterway is aligned with the chosen grid).  It is also interesting to consider Figure~\ref{fig:num_grids_cells} that provides some insight into the behavior of the presented simulation.  As the storm approaches shore, \geoclaw does not refine much until the storm is approximately 18 hours away from landfall.  At this time the resolution increases dramatically and stays relatively high, only decreasing slowly for the rest of the simulation time.  An improved refinement criteria may be able to address the apparent discrepancy in resolution while balancing that with decrease of resolution later in the simulation.

Neither simulation fully captures the true surge as was seen in Figures~\ref{fig:gauge_1}, \ref{fig:gauge_2}, \ref{fig:gauge_3}, and \ref{fig:gauge_4}.  In a forecasting scenario, capturing the forerunner would be a pertinent feature and it is unclear what caused both simulations to miss this feature.  We have also purposely excluded a number of forcing functions (\emph{e.g.} tidal constituents) which may have been important.  Currently tides have not been implemented in conjunction with AMR and may be an important issue to address in the future, especially where tidal forcing is more important such as the northern Atlantic coast.

One important issue that should be mentioned is that the performance of an AMR scheme, in this case \geoclaw, is highly dependent on the type of refinement criteria in conjunction with the base resolution and ratios used.  The setup here was picked after some experience running models with \geoclaw but should not be taken as the optimal settings.  In the course of writing this article these settings were changed and run multiple times leading to speed-up of a factor of 10 from the initial settings.  Clearly more work needs to be spent on what are optimal AMR settings for storm surge modeling.

% ==============================================================================
\section{Conclusions} \label{sec:conclusions}

At the outset of this article two capabilities were mentioned as essential for storm surge forecasting, ensemble based calculations and simulations containing resolution sufficient to capture multiple length scales.  The numerical models SLOSH and \adcirc address each of these capabilities independently but not simultaneously.  The question that is addressed here then is not if an AMR based code, such as \geoclaw, is better at either of these capabilities separately but if it can satisfy both without overly sacrificing either need.  The capability to calculate ensemble based forecasts is answered in Table~\ref{tab:timings}.  Given a roughly 2 hour forecasting window, on the same machine running 4 \adcirc simulations in this window, \geoclaw could be run over 1000 times.  It is also clear from the simulations presented that with AMR we are able to capture many of the fine-scale features near coastlines that \adcirc is able to capture.  Combined, these results strongly suggest that AMR is a compelling way to forecast storm surge with both ensemble calculation and high-resolution capabilities.

An ancillary impact of using AMR relates to the amount of effort needed to create the detailed unstructured grids that the \adcirc modelers actively improve upon.  This is an especially important consideration for regions at risk which do not have the type of resources required to create these detailed grids ahead of a storm, both in terms of personnel and software.  \geoclaw substantially reduces the effort needed and resources required in order to reach solutions, considering the growing ubiquity of bathymetry and land-use data.  Furthermore, \geoclaw also allows the grid to adapt to the storm being run, reducing computational cost.

Future improvements to the \geoclaw framework, as it pertains to storm surge modeling, primarily involve incorporating detailed variable friction fields and investigating possible improvements to the storm field representations.  The parallelism present in \geoclaw via OpenMP could be improved.  In particular, load-balancing and scalability on high numbers of threads will be critical in the future due to the current trend of higher core counts on consumer chips and the emerging many-core architectures.  More broadly, the use of AMR in storm surge forecasting will require additional work on the refinement criteria, balancing resolution with cost.  It may also be worthwhile to further delve into what aspect of the AMR implemented here provided the largest impact and adopt them into existing codes which are already in operational use.

% ==============================================================================
%  Acknowledgments
\vskip 10pt
{\bf Acknowledgments.}
The authors would like to thank Andrew Kennedy for the gauge observations and the reviewers for their time and effort in editing this article.  This research was supported in part by ONR Grant N00014-09-1-0649, the ICES postdoctoral fellowship program, the Gulf of Mexico Research Initiative Center for Advanced Research on Transport of Hydrocarbons in the Environment, and the King Abdullah University of Science and Technology Academic Excellence Alliance.

\bibliographystyle{elsarticle-num}
\bibliography{database}

\appendix

\section{Wave Propagation Methods} \label{app:wave_propagation_methods}

Wave-propagation methods start with the finite volume discretization of the domain into cells $\cell_i$ and the objective to evolve the grid cell average defined as $Q^n_i \equiv \frac{1}{\dx_i} \int_{\cell_i} q(x,t^n) dx$.  The goal then is to evolve these cell averages forward in time as the system of PDEs $q_t + f(q)_x = 0$ prescribes.  To accomplish this, the grid cell averages are used to reconstruct a piece-wise constant representation of the function $q$ throughout the domain.  This reconstruction leads naturally to the formulation of Riemann problems at each grid cell interface.  A Riemann problem consists of the original system of PDEs on an infinite domain with a jump-discontinuity located at $x=0$.  The solution to a Riemann problem generally consists of $m$ waves (often $m$ is equal to the number of equations in the system of PDEs but this is not a requirement) denoted by $\wave^p \in \R^m$ propagating away from the location of the jump-discontinuity traveling at speeds $s^p$.  These waves are related to the original jump-discontinuity via
\begin{equation} \label{eq:wave_definition}
    Q_i - Q_{i-1} = \sum^m_{p=1} \wave^p_{i-1/2}.
\end{equation}
From these waves, a first-order upwind method can be written as
\[
    Q^{n+1}_{i} = Q^{n}_i - \frac{\dt}{\dx} \left [ \apdq_{i-1/2} - \amdq_{i+1/2} \right]
\]
where $\apdq$ and $\amdq$ represent fluctuations coming from the right and left cell interfaces respectively and can be defined in terms of the waves as
\[
    \mathcal{A}^\pm \Delta Q_{i-1/2} = \sum^m_{p=1} (s^p_{i-1/2})^\pm \wave^p_{i-1/2}
\]
where $s^+_{i-1/2} = \text{max}(s^p_{i-1/2},0)$ and $s^-_{i-1/2} = \text{min}(s^p_{i-1/2},0)$.  Extensions to $2^\text{nd}$-order are possible using limiters applied to each wave such that the update becomes
\[
    Q^{n+1}_{i} = Q^n_i - \frac{\dt}{\dx} (\amdq_{i+1/2} + \apdq_{i-1/2}) - \frac{\dt}{\dx} (\tilde{F}_{i+1/2} - \tilde{F}_{i-1/2})
\]
where 
\[
    \tilde{F}_{i-1/2} = \frac{1}{2} \sum^m_{p=1} |s^p_{i-1/2}| \left (1 - \frac{\dt}{\dx}|s^p_{i-1/2}| \right ) \widetilde{\wave}^p_{i-1/2}
\]
and $\widetilde{\wave}^p_{i-1/2}$ are limited versions of $\wave^p_{i-1/2}$.

In the case where the system being considered is linear, the waves $\wave^p$ can be written as scalar multiples $\alpha^p$ of the right eigenvectors $r^p$ of the flux where $f(q)_x = A q_x$ such that
\[
    \wave^p_{i-1/2} = (\ell^p_{i-1/2})^T(Q_i - Q_{i-1}) r^p_{i-1/2}
\]
where $\ell^p$ are the left eigenvectors of $A$.  If the system is non-linear often local linearizations are employed to find an approximate flux Jacobian $\hat{A}_{i-1/2}$ whose eigenvectors are again used to find $\wave^p$.

An alternative formulation to the splitting in \eqref{eq:wave_definition} is to instead split the jump in the fluxes at each cell interface as
\[
    f(Q_i) - f(Q_{i-1}) = \sum^m_{p=1} \fwave^p_{i-1/2}
\]
where now $\fwave^p_{i-1/2}$ waves carrying a jump in the fluxes.  This formulation, called the f-wave approach, has the advantage that spatially dependent flux functions and source terms can be incorporated directly into the Riemann solver.  In the case of an f-wave splitting, the fluctuations $\apdq$ and $\amdq$ can be written as
\[
    \mathcal{A}^\pm \Delta Q_{i-1/2} = \sum^m_{p=1} \text{sgn}(s^p_{i-1/2}) \fwave^p_{i-1/2}.
\]
The accuracy again can be improved as before by limiting the f-waves.

\section{Augmented Riemann Solver} \label{app:riemann_solver}

The augmented Riemann solver used in \geoclaw was originally proposed in \cite{George:2008aa} where the novel approach of splitting both the conserved quantities $Q$ and the fluxes $f(Q)$ into waves was explored.  Here the augmented set of waves take the form
\begin{equation} \label{eq:augmented_wave_vector}
    \begin{bmatrix}
        Q_i - Q_{i-1} \\
        f(Q_i) - f(Q_{i-1})
    \end{bmatrix} = 
    \begin{bmatrix}
        \sum^m_{p=1} \wave^p_{i-1/2} \\
        \sum^m_{p=1} \fwave^p_{i-1/2}
    \end{bmatrix}
\end{equation}
leading to $2m$ waves.  Since only $m$ waves are needed to update the $m$ equations, these extra waves can be used to provide many desirable properties to the Riemann solution.

Another primary feature of the Riemann solver in \geoclaw is the direct inclusion of the bathymetry source term $-ghb_x$ into the Riemann solver.  One important consequence of this incorporation in conjunction with the augmented approach is the genesis of a ``steady state'' wave.  As was shown in \cite{George:2008aa} this additional wave preserves steady states where $(hu)_x = 0$, a much more general class of steady states than is usually considered.  Additionally this extra wave can be used to better represent large rarefaction waves which are commonly found at the wet-dry front.

\end{document}